\documentclass[
preprint, 3p, 
number, 
sort&compress,
]{elsarticle}
\pdfoutput=1

\usepackage[utf8]{luainputenc}
\usepackage[english]{babel}
\usepackage{csquotes}

\usepackage[plainpages=false,pdfpagelabels,hidelinks,unicode]{hyperref}

\usepackage{lipsum}
\usepackage{amsfonts}
\usepackage{graphicx}
\usepackage{epstopdf}

\usepackage{amsmath}
\allowdisplaybreaks
\usepackage{amssymb}
\usepackage{commath}
\usepackage{mathtools}
\usepackage{bbm}

\usepackage{siunitx}
\usepackage{adjustbox}

\usepackage{algorithm}
\usepackage[noend]{algpseudocode}

\usepackage[normalem]{ulem}
\newcommand{\rulesep}{\unskip\ \hrule\ }

\usepackage{amsthm}
\theoremstyle{plain}
  \newtheorem{theorem}{Theorem}

\theoremstyle{definition}
  
  \newtheorem{remark}[theorem]{Remark}

\usepackage{color}
\usepackage{graphicx}
\usepackage[small]{caption}
\usepackage{subcaption}

\ifx\useTikzForPlotting\undefined
\else
  \usepackage{pgfplots}
  \pgfplotsset{compat=1.11}
  \usetikzlibrary{external}
\fi

\usepackage{booktabs}
\usepackage{rotating}
\usepackage{multirow}

\usepackage{multicol}
\usepackage{enumitem}

\usepackage{calc}
\usepackage{xparse}

\DeclareMathOperator*{\argmin}{arg\,min}
\renewcommand{\vec}[1]{\underline{#1}}
\NewDocumentCommand{\mat}{mo}{%
  \IfValueTF{#2}{%
    \underline{\underline{#1}}{#2}
  }{%
    \underline{\underline{#1}}\,
  }%
}
\newcommand{\diag}[1]{\operatorname{diag}\left(#1\right)}

\newcommand{\sign}{\operatorname{sign}}
\newcommand{\scp}[2]{\left\langle{#1,\, #2}\right\rangle}

\renewcommand{\d}{\mathrm{d}}

\newcommand{\fnum}{f^{\mathrm{num}}}

\newcommand{\vecfnum}{\vec{f}^{\mathrm{num}}}

\newcommand{\Oref}{\Omega_{\mathrm{ref}}}

\renewcommand{\epsilon}{\varepsilon}
\renewcommand{\phi}{\varphi}
\renewcommand{\rho}{\varrho}

\newcommand{\R}{\mathbb{R}}

\newsavebox{\DelimiterBox}
\newlength{\DelimiterHeight}
\newlength{\DelimiterDepth}
\newsavebox{\ArgumentBox}
\newlength{\ArgumentHeight}
\newlength{\ArgumentDepth}
\newlength{\ResizedDelimiterHeight}
\newlength{\ResizedDelimiterDepth}

\usepackage{lipsum}
\makeatletter
\def\ps@pprintTitle{%
 \let\@oddhead\@empty
 \let\@evenhead\@empty
 \def\@oddfoot{}%
 \let\@evenfoot\@oddfoot}
\makeatother

\begin{document}

\begin{frontmatter}

\title{High Order Edge Sensors with $\ell^1$ Regularization for Enhanced Discontinuous Galerkin Methods}

\author[label1]{Jan Glaubitz\corref{cor1}}
\ead{j.glaubitz@tu-bs.de}

\author[label2]{Anne Gelb}
\ead{annegelb@math.dartmouth.edu}

\cortext[cor1]{Corresponding author: Jan Glaubitz}
\address[label1]{Institute for Computational Mathematics, TU Braunschweig, 38106 Braunschweig, Germany.} 
\address[label2]{Department of Mathematics, Dartmouth College, Hanover, NH 03755, United States.}

\begin{abstract}
  This paper investigates the use of $\ell^1$ regularization for solving hyperbolic conservation laws 
based on high order discontinuous Galerkin (DG) approximations. 
We first use the polynomial annihilation method to construct a high order edge sensor which 
enables us to flag ``troubled'' elements. 
The DG approximation is enhanced in these troubled regions by activating $\ell^1$ regularization to 
promote sparsity in the corresponding jump function of the numerical solution.  
The resulting $\ell^1$ optimization problem is efficiently implemented using the alternating 
direction method of multipliers. 
By enacting $\ell^1$ regularization only in troubled cells, our method remains accurate and 
efficient, as no additional regularization or expensive iterative procedures are needed in 
smooth regions. 
We present results for the inviscid Burgers' equation as well as for a nonlinear system of 
conservation laws using a nodal collocation-type DG method as a solver.

\end{abstract}


\end{frontmatter}

\section{Introduction}
\label{sec:introduction}

This work is about a novel shock capturing procedure in spectral element (SE) type of methods for solving 
time-dependent hyperbolic conservation laws 
\begin{equation}\label{eq:cons-law}
    \partial_t u + \partial_x f(u) = 0, 
\end{equation} 
with smooth and discontinuous solutions. 
Such methods include the spectral difference (SD) \cite{liu2006spectral,wang2007spectral}, discontinuous Galerkin (DG) 
\cite{hesthaven2007nodal}, and flux reconstruction \cite{huynh2007flux} methods. 
While all these methods provide fairly high orders of accuracy for smooth problems, they often lack desired stability and robustness properties, especially in the presence of (shock) discontinuities.  

Many shock capturing techniques have therefore been developed over the last few decades.
Such efforts date back more than $60$ years to the pioneering work 
of von Neumann and Richtmyer, \cite{vonneumann1950method} in which they add artificial viscosity terms to 
\eqref{eq:cons-law} in order to construct stable finite difference 
schemes for the equations of hydrodynamics.
Since then, artificial viscosity has been added to a variety of  algorithms
\cite{jameson1981numerical,persson2006sub,klockner2011viscous,nordstrom2006conservative,ranocha2018stability,glaubitz2019smooth}. 
However, augmenting \eqref{eq:cons-law} with additional (higher) viscosity terms requires care about 
their design and size. 
Otherwise, new time stepping constraints for explicit methods can considerably decrease computational efficiency; 
see, e.g. \cite[equation (2.1)]{klockner2011viscous}.
Other interesting alternatives are based on (modal) filters \cite{hesthaven2008filtering,glaubitz2018application,ranocha2018stability} 
or applying viscosity to the different spectral scales \cite{tadmor1990shock,tadmor1993super}.
Finally, we mention those methods based on order reduction 
\cite{cockburn1989tvb}, 
mesh adaptation \cite{dervieux2003theoretical}, 
and weighted essentially nonoscillatory (WENO) concepts \cite{shu1988efficient,shu1989efficient}.  
Yet, a number of issues regarding the fundamental convergence properties for these methods still 
remain unresolved. 
Moreover, even when the extension to multiple dimensions is straightforward, these schemes may be 
too computationally expensive for practical usage. 

In this work we propose $\ell^1$ regularization as a novel tool to capture shocks in SE methods by promoting sparsity in 
the jump function of the approximate solution.  
$\ell^1$ regularization methods are frequently encountered in signal processing and imaging applications.
They are still of limited use in solving partial differential equations numerically, however, and only a few studies, 
(see, e.g., 
\cite{schaeffer2013sparse,hou2015sparse+,lavery1989solution,lavery1991solution,scarnati2017using,guermond2008fast})
have considered sparsity or $\ell^1$ regularization of the numerical solution. 
A brief discussion of these investigations can be found in \cite{scarnati2017using}.
We note that while problems with discontinuous initial conditions were studied in 
\cite{lavery1989solution,lavery1991solution,guermond2008fast}, problems that form shocks were not. 
The technique developed in  \cite{schaeffer2013sparse} was designed to promote sparsity in the 
frequency domain, making it less amenable to problems emitting shocks, where the frequency domain is 
not sparse. 
Further, with the exception of \cite{scarnati2017using}, 
each of these investigations applied $\ell^1$ regularization  directly to numerical solution or to 
its residual, rather than incorporating it directly in the time stepping evolution. 

In this investigation we follow the approach in \cite{scarnati2017using}, which incorporates 
$\ell^1$ regularization directly into the time dependent solver. 
Specifically, we promote {\em the sparsity of the jump function} that corresponds to the 
discontinuous solution. 
The jump function approximation is performed using the (high order) polynomial annihilation (PA) 
operator \cite{archibald2016image,archibald2005polynomial,wasserman2015image}, which eliminates the 
unwanted ``staircasing'' effect, a common degradation of detail of the piecewise smooth solution 
arising from the classical total variation (TV) regularization. 
More specifically, the high order PA operator allows the resulting solution to be comprised of 
piecewise polynomials instead of piecewise constants.
We solve the resulting $\ell^1$ optimization problem by the alternating direction method of 
multipliers (ADMM) \cite{li2013efficient,sanders2016matlab,sanders2017composite}. 
A similar application of $\ell_1$ regularization was used in \cite{scarnati2017using} to numerically solve hyperbolic 
conservation laws, though only for the Lax--Wendroff scheme and 
Chebyshev and Fourier spectral methods. 
It was concluded in \cite{scarnati2017using} that 
although the Lax--Wendroff scheme yielded sufficient accuracy for relatively simple problems, its 
lower order convergence properties made it difficult to resolve more complicated ones. 
The new technique fared better using Chebyshev polynomials, although their global construction made 
it difficult to resolve the local structures without oscillations or excessive smoothing.

One possible solution is to use SE methods as the underlying mechanism for 
solving the hyperbolic conservation law.  SE methods have the advantage of being more localized, 
for instance, and allow element-to-element variations in the optimization problem. 
In particular for the method we develop here, $\ell^1$ regularization is only activated in troubled elements, i.e., in 
elements where discontinuities are detected. 
This further enhances efficiency of the method. 
In the process, a novel discontinuity sensor based on PA operators of increasing orders is proposed, which is able to flag troubled elements. 
The discontinuity sensor steers the optimization problem and thus locally calibrates the method with respect to the 
smoothness of the solution.  
Numerical tests are performed for a nodal collocation-type DG method and the inviscid Burgers' equation as well as for a 
nonlinear system. 
It should be stressed that the proposed procedure also carries over to other classes of methods, 
with the obvious extension to SE type methods.
The extension to other types of methods, such as finite volume methods, is also possible under slight modifications of the procedure. 
We should also stress that in our new development of the $\ell^1$ regularization method for 
solving PDEs that admit shocks in their solutions that we {\em do not} require different methods to 
be used in smooth and nonsmooth regions. 
Such methods have been developed in, for instance, \cite{don2016hybrid} and have been shown to be effective. 
Here we demonstrate that it is possible to avoid such additional complexities.

The rest of this paper is organized as follows: 
\S \ref{sec:preliminaries} briefly reviews the nodal collocation-type DG method, $\ell^1$ regularization, and PA 
operators which are needed for the development of our method. 
In \S \ref{sec:l1_limiter} we describe the application of $\ell^1$ regularization by higher-order 
edge detectors to SE type methods. 
Further, a novel discontinuity sensor based on PA operators of increasing orders is proposed. 
Numerical tests for the inviscid Burgers' equation, the linear advection equation, 
and a nonlinear system of conservation laws are presented in  \S \ref{sec:tests}. 
The tests demonstrate that we are able to better resolve numerical solutions when $\ell^1$ regularization is utilized. 
We close this work with concluding thoughts in  \S \ref{sec:conclusion}.

\section{Preliminaries}
\label{sec:preliminaries}

In this section we briefly review all necessary concepts in order to introduce $\ell^1$ regularisation into the framework of discontinuous Galerkin methods in the subsequent section.

\subsection{A nodal discontinuous Galerkin method}
\label{sub:DG}

Let us consider a hyperbolic conservation law 
\begin{equation}\label{eq:cl}
    \partial_t u + \partial_x f\left(u\right) = 0, \quad x \in \Omega,
\end{equation}
with suitable initial and boundary conditions. 
The domain $\Omega \subset \R$ is decomposed into $I$ disjoint, face-conforming elements $\Omega_i$, 
$\Omega = \bigcup_{i=1}^I \Omega_i$. 
All elements are mapped to a reference element, typically $\Oref = [-1,1]$, where all computations are performed. 

In this work we consider a nodal collocation-type discontinuous Galerkin (DG) method 
\cite{hesthaven2007nodal} on the reference element. 
The solution $u$ as well as the flux function $f$ are approximated by interpolation polynomials of 
the same degree, giving the advantage of highly efficient operators. 
We further collocate the flux approximation based on interpolation with the numerical 
quadrature used for the evaluation of the inner products \cite{kopriva2009implementing}.

The first step is to introduce a nodal polynomial approximation 
\begin{equation}\label{eq:approx-u} 
        u\left(t,\xi\right) \approx u_p\left(t,\xi\right) = \sum_{k=0}^p u_k\left(t\right) \ell_k\left(\xi\right), 
\end{equation}
where $p$ is the polynomial degree and $\{ u_k \}_{k=0}^p$ are the $p+1$ time dependent nodal degrees of freedom at the element grid nodes $-1 \leq \xi_0 < \dots < \xi_p \leq 1$. 
Common choices are either the Gauss-Legendre or Gauss-Lobatto nodes \cite{davis2007methods}. 
We use Gauss-Lobatto points in the latter numerical tests, as they include the boundary points $\xi_0 = -1, \xi_p=1$ and thus render the method more robust 
\cite{gassner2013skew,gassner2014kinetic,kopriva2014energy}.  
For the procedure of $\ell^1$ regularisation proposed in this work the specific choice of grid nodes is not crucial however. 

Further, the Lagrange basis functions of degree $p$ are given by 
\begin{equation}
    \ell_k\left(\xi\right) = \prod_{\begin{smallmatrix}i=0 \\ i\neq k\end{smallmatrix}}^p \frac{\xi - \xi_i}{\xi_k - \xi_i}
\end{equation} 
and satisfy the cardinal property $\ell_k\left(\xi_i\right) = \delta_{ki}$. 
The flux function is approximated in the same way, i.e. 
\begin{equation}\label{eq:approx-f} 
    f\left(u\right) \approx f_p\left(t,\xi\right) = \sum_{k=0}^p f_k\left(t\right) \ell_k\left(\xi\right),
\end{equation} 
where, collocating the nodes for both approximations, the nodal degrees of freedom $\{ f_k \}_{k=0}^p$ are given by $f_k\left(t\right) = f\left(u_k\left(t\right)\right)$. 

We now obtain the formulation of the nodal DG method by inserting the polynomial approximations \eqref{eq:approx-u} and \eqref{eq:approx-f} into the conservation law \eqref{eq:cl}, multiplying by a test function $\ell \in \{ \ell_k \}_{k=0}^p$, integrating over the reference element, and applying integration by parts, resulting in 
\begin{equation}\label{eq:DG-int}
    \int_{-1}^1 \dot{u}_p \ell_i \d \xi + \Bigl. \left( \fnum \ell_i \right)\Bigr|_{-1}^1 
        - \int_{-1}^1 f_p \ell_i' \d \xi, 
        \quad i=0,\dots,p.
\end{equation}
Here, $\fnum$ is a suitably chosen numerical flux, providing a mechanism to couple the
solutions across elements \cite{toro2013riemann}. 
Further, $\dot{u}_p$ denotes the time derivative of the approximation while $\ell_i'$ denotes the spatial
derivative of the basis element with respect to $\xi$. 

Next, the integrals are approximated by an (interpolatory) quadrature rule using the same nodes, 
\begin{equation}
    \int_{-1}^1 g\left(\xi\right) \d \xi 
        \approx \int_{-1}^1 \sum_{k=0}^p g\left(\xi_k\right) \ell_k\left(\xi\right) \d \xi 
        = \sum_{k=0}^p \omega_k g\left(\xi_k\right). 
\end{equation}
From this quadrature rule, we introduce a discrete inner product  
\begin{equation}\label{eq:disc_inner_prod}
    \scp{u}{v}_{\mat{M}} 
        = \sum_{k=0}^p \omega_k u\left(\xi_k\right) v\left(\xi_k\right) 
        = \vec{u}^T \mat{M} \vec{v}
\end{equation} 
with mass matrix 
\begin{equation}
    \mat{M} = \operatorname{diag}\left( [\omega_0,\dots,\omega_p] \right)
\end{equation} 
and vectors of nodal degrees of freedom 
\begin{equation}
    \vec{u} = [ u_0, \dots, u_p ]^T, \quad \vec{v} = [ v_0, \dots, v_p ]^T. 
\end{equation}
Using the discrete inner product, the spacial approximation \eqref{eq:DG-int} becomes 
\begin{equation}
    \scp{\dot{u}_p}{\ell_i}_{\mat{M}} = \scp{f_p}{\ell_i'}_{\mat{M}} - \Bigl. \left( \fnum \ell_i \right)\Bigr|_{-1}^1, \quad i=0,\dots,p. 
\end{equation}
Finally going over to a matrix vector representation and utilising the cardinal property of the Lagrange basis functions, the DG approximation can be compactly rewritten in its \emph{weak form} as 
\begin{equation}\label{eq:DG-weak}
    \mat{M} \dot{\vec{u}} = 
        \mat{D}^T \mat{M} \vec{f} - \mat{R}^T \mat{B} \vecfnum, 
\end{equation} 
where 
\begin{equation}
    \mat{D} = \left( \ell_i'\left(\xi_k\right) \right)_{k,i=0}^p, \quad
    \mat{R} = 
        \begin{pmatrix} 
            0 & \dots & 0 & 1 \\ 
            1 & 0 & \dots & 0
        \end{pmatrix}, \quad 
    \mat{B} = 
        \begin{pmatrix} 
            -1 & 0 \\ 
            0 & 1
        \end{pmatrix}
\end{equation} 
are the differentiation matrix, restriction matrix, and boundary matrix, respectively. 
Let $\vecfnum$ denote the vector containing the values of the numerical flux at the element boundaries. 
Note that applying integration by parts a second time to \eqref{eq:DG-int} would result in the \emph{strong form} 
\begin{equation}\label{eq:DG-strong}
    \mat{M} \dot{\vec{u}} = 
        - \mat{M} \mat{D} \vec{f} - \mat{R}^T \mat{B} \left( \vecfnum - \mat{R} \vec{f} \right)
\end{equation} 
of the DG approximations. 
Both forms are equivalent when using summation by parts operators, which satisfy  
${\mat{M} \mat{D} + \mat{D}^T \mat{M} = \mat{R}^T \mat{B} \mat{R}}$. 
The strong form \eqref{eq:DG-strong} can further be recovered as a special case of flux reconstruction schemes \cite{allaneau2011connections,yu2013connection,de2014connections}.

\subsection{\texorpdfstring{$\ell^1$ regularisation}{l1 regularisation}}
\label{sub:l1_reg}

Let $u\left(\xi\right) = u\left(t,\xi\right)$ be the unknown solution on an element $\Omega_i$ transformed into the reference element $\Oref$ and $u_p \in \mathbb{P}_p\left(\Oref\right)$ a spacial polynomial approximation at fixed time $t$. 
Assume that some measurable features of $u$ have sparse representation. 
Consequently, the approximation $u_p$ is desired to have this sparse representation as well. 

Let $H$ be a regularisation functional which measures sparsity. 
The objective is to then solve the constrained optimisation problem 
\begin{equation}\label{eq:const_problem}
    \argmin_{v \in \mathbb{P}_p\left(\Oref\right)} H\left(v\right) 
    \quad \text{s.t.} \quad 
    \norm{v-u_p} = 0.
\end{equation}
The equality constraint, referred to as the \emph{data fidelity term}, measures how accurately the reconstructed approximation fits the given data with respect to some seminorm $\norm{\cdot}$. 
Typically, the continuous $L^2$-norm $\norm{f}_2^2 = \int |f|^2$ or some discrete counterpart is used. 
The \emph{regularisation term} $H\left(v\right)$ enforces the known sparsity present in the underlying solution $u$ by penalising missing sparsity in the approximation. 
The regularisation functional $H$ further restricts the approximation space to a desired class of functions, here $\mathbb{P}_p\left(\Oref\right)$. 
Note that any $p$-norm with $p \leq 1$ will enforce sparsity in the approximation. 
In this work, we choose $H$ to be the $\ell^1$-norm of certain transformations of $v$. 

It should be stressed that if $\norm{\cdot}$ is not just a seminorm but a strictly convex norm, 
for instance induced by an inner product, the equality constraint immediately and uniquely 
determines the approximation.
Thus, instead of \eqref{eq:const_problem}, typically the related denoising problem 
\begin{equation}\label{eq:relaxed_problem}
    \argmin_{v \in \mathbb{P}_p\left(\Oref\right)} H\left(v\right) 
    \quad \text{s.t.} \quad 
    \norm{v-u_p} < \sigma
\end{equation}
with $\sigma > 0$ is solved, which relaxes the equality constraint on the data fidelity term. 
Equivalently, the denoising problem \eqref{eq:relaxed_problem} can also be formulated as the \emph{unconstrained (or penalised) problem} 
\begin{equation}\label{eq:unconst_problem}
    \argmin_{v \in \mathbb{P}_p\left(\Oref\right)} \left( \norm{v-u_p}_2^2 + \lambda H\left(v\right) \right)
\end{equation}
by introducing a non-negative \emph{regularisation parameter} $\lambda \geq 0$. 
$\lambda$ represents the trade-off between fidelity to the original approximation and sparsity. 

The unconstrained problem \eqref{eq:unconst_problem} is often solved with $H\left(v\right) = TV\left(v\right)$, where $TV(v)$ is the total variation of $v$. 
Following \cite{scarnati2017using} however, in this work we solve \eqref{eq:unconst_problem} with 
\begin{equation}
    H\left(v\right) = \norm{ L_m v }_1,
\end{equation}
where $L_m$ is a polynomial annihilation operator introduced in the next subsection. 
Using higher order polynomial annihilation operators will help to eliminate the staircase effect that occurs when using the TV operator (polynomial annihilation for $m=1$) for $H$, see \cite{stefan2010improved}.

\subsection{Polynomial annihilation}
\label{sub:poly_anni}

Polynomial annihilation (PA) operators were originally proposed in \cite{archibald2005polynomial}.  One 
main advantage in using PA operators of higher orders ($m>1$) for the regularisation functional 
$H$ is that they allow distinction between jump discontinuities and steep gradients, 
which is critical in the numerical treatment of nonlinear conservation laws. 
PA regularisation is also preferable to TV regularisation when the resolution is poor, 
even when the underlying solution is piecewise constant.

Let $u\left(\xi^-\right)$ and $u\left(\xi^+\right)$ respectively denote the left and right hand 
side limits of $u:\Oref=[a,b] \to \R$ at $\xi$. 
We define the \emph{jump function of $u$} as 
\begin{equation}\label{eq:jump_fun}
    [u]\left(\xi\right) = u\left(\xi^+\right) - u\left(\xi^-\right)
\end{equation} 
and note that $[u]\left(\xi\right) = 0$ at every $\xi \in \Oref$ where $u$ has no jump. 
We thus say that the jump function $[u]$ has a sparse representation. 
The \emph{polynomial annihilation operator of order $m$}, 
\begin{equation}\label{eq:PA-operator}
    L_m[u]\left(\xi\right) = \frac{1}{q_m\left(\xi\right)} 
    \sum_{x_j \in S_\xi} c_j\left(\xi\right) u\left(x_j\right), 
\end{equation}
is designed in order to approximate the jump function $\left[ u \right]$. 
Here
\begin{equation}
    S_\xi = \{ x_0\left(\xi\right), \dots, x_m\left(\xi\right) \} \subset \Oref
\end{equation}
is a set of $m+1$ local grid points around $\xi$, the \emph{annihilation coefficients} $c_j:\Oref \to \R$ are given by 
\begin{equation}\label{eq:syst-anni-coef}
    \sum_{x_j \in S_\xi} c_j\left(\xi\right) p_l\left(x_j\right) 
      = p_l^{\left(m\right)}\left(\xi\right), \quad j=0,\dots,m,
\end{equation}
and $\{ p_l \}_{l=0}^m$ is a basis of $\mathbb{P}_m\left(\Oref\right)$. 
An explicit formula for the annihilation coefficients utilising Newton's divided differences 
is given by (\cite{archibald2005polynomial}):
\begin{equation}\label{eq:anni-coef}
    c_j\left(\xi\right) = \frac{m!}{\omega_j\left(S_\xi \right)} 
    \quad \text{with} \quad 
    \omega_j\left(S_\xi \right) = \prod_{\begin{smallmatrix} x_i \in S_\xi \\i\neq j\end{smallmatrix}} \left(x_j - x_i\right)
\end{equation}
for $j=0,\dots,m$.
Finally, the \emph{normalisation factor} $q_m$, calculated as
\begin{equation}\label{eq:norm-fact}
    q_m\left(\xi\right) = \sum_{x_j \in S_\xi^+} c_j\left(\xi\right),
\end{equation}
ensures convergence to the right jump strength at every discontinuity. 
Here $S_\xi^+ $ denotes the set $\{ x_j \in S_\xi \ | \ x_j \geq \xi \}$ of all local grid points to the right of $\xi$. 

In this work, the PA operator is applied to the reference element $\Oref = [-1,1]$ of the 
underlying nodal DG method using $p+1$ collocation points $\{ \xi_k \}_{k=0}^p$, 
typically Gauss-Lobatto points including the boundaries. 
We can thus construct polynomial annihilation operators up to order $p$ by allowing 
the sets of local grid points $S_\xi$ to be certain subsets of the $p+1$ collocation points. 

In \cite{archibald2005polynomial} it was shown that
\begin{equation}\label{eq:convergence}
    L_m[u]\left(\xi\right) = 
    \left\{ 
    \begin{array}{ccl}
        [u]\left( x \right) + \mathcal{O}\left(h\left( \xi \right)\right) & , & 
	  \text{if } x_{j-1} \leq \xi, x \leq x_j, \\ 
        \mathcal{O}\left(h^{\min\left(m,k\right)}\left(\xi\right)\right) & , & 
	  \text{if } u \in C^k\left(I_{\xi}\right),
    \end{array} \right.
\end{equation}
where $h\left(\xi\right)= \max\{ |x_i - x_{i-1} \ | \ x_{i-1}, x_i \in S_{\xi} \}$ and $I_{\xi}$ is 
the smallest closed interval such that $S_{\xi} \subset I_{\xi}$. 
Note that $h(\xi)$ depends on the density of set of local grid points around $\xi$.
Thus, $L_m$ provides $m$th order convergence in regions where $u \in C^m$, and yields a first order approximation to the jump.
It should be stressed that oscillations develop around points of discontinuity as $m$ increases. 
The impact of the oscillation could be reduced by post-processing methods, 
such as the minmod limiter, as was done in \cite{archibald2005polynomial}. 
However, as long as there is sufficient resolution between two shock locations,
such oscillations do not directly impact our method. 
This is because we use the PA operator not to detect the precise location of jump discontinuities, 
but rather to enforce sparsity. 

\begin{figure}[!htb]
  \centering
  \begin{subfigure}[b]{0.32\textwidth}
    \includegraphics[width=\textwidth]{%
      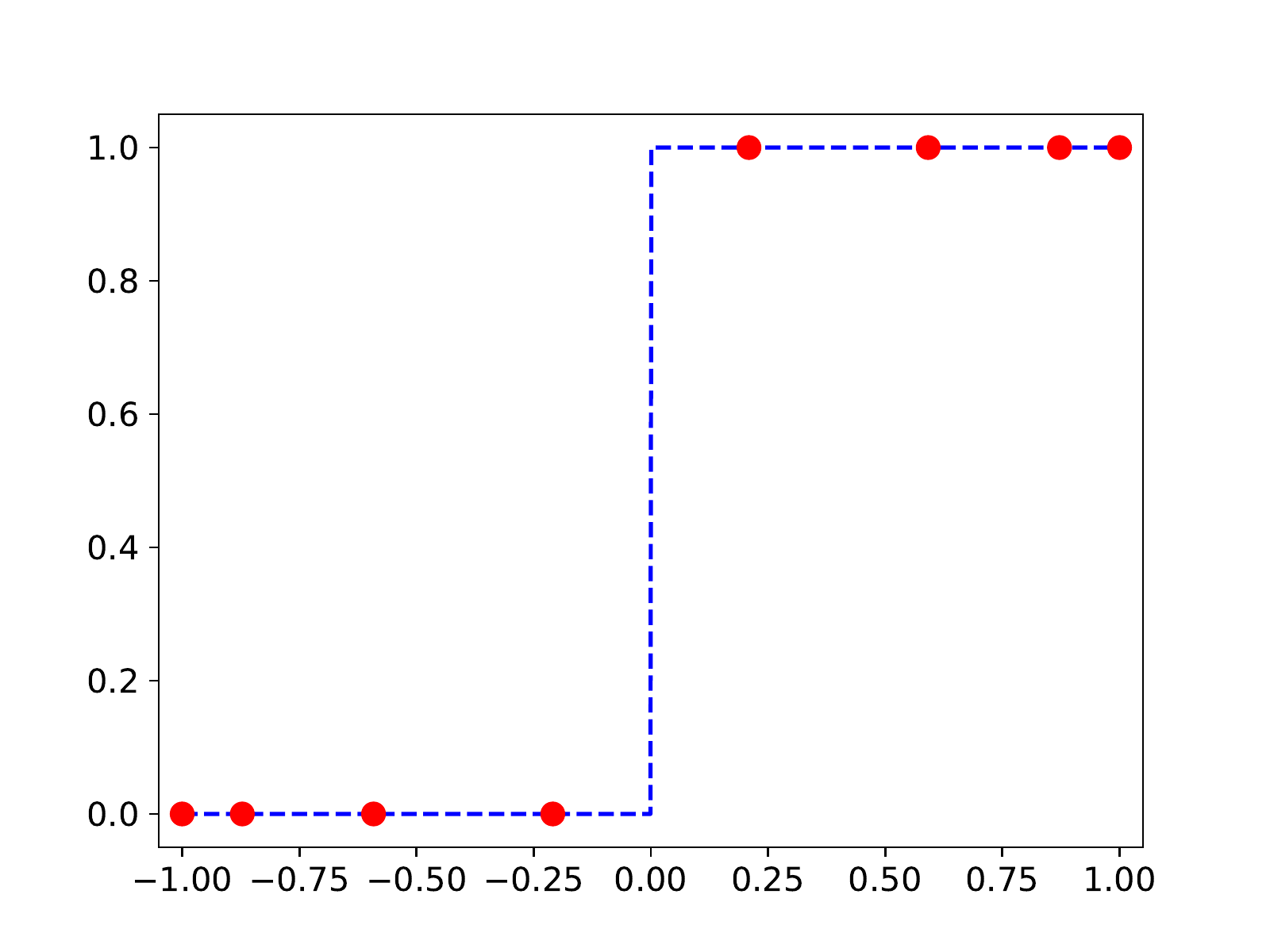}
    \caption{$f_1(\xi)$}
    \label{fig:step_0_fun}
  \end{subfigure}%
  ~
  \begin{subfigure}[b]{0.32\textwidth}
    \includegraphics[width=\textwidth]{%
      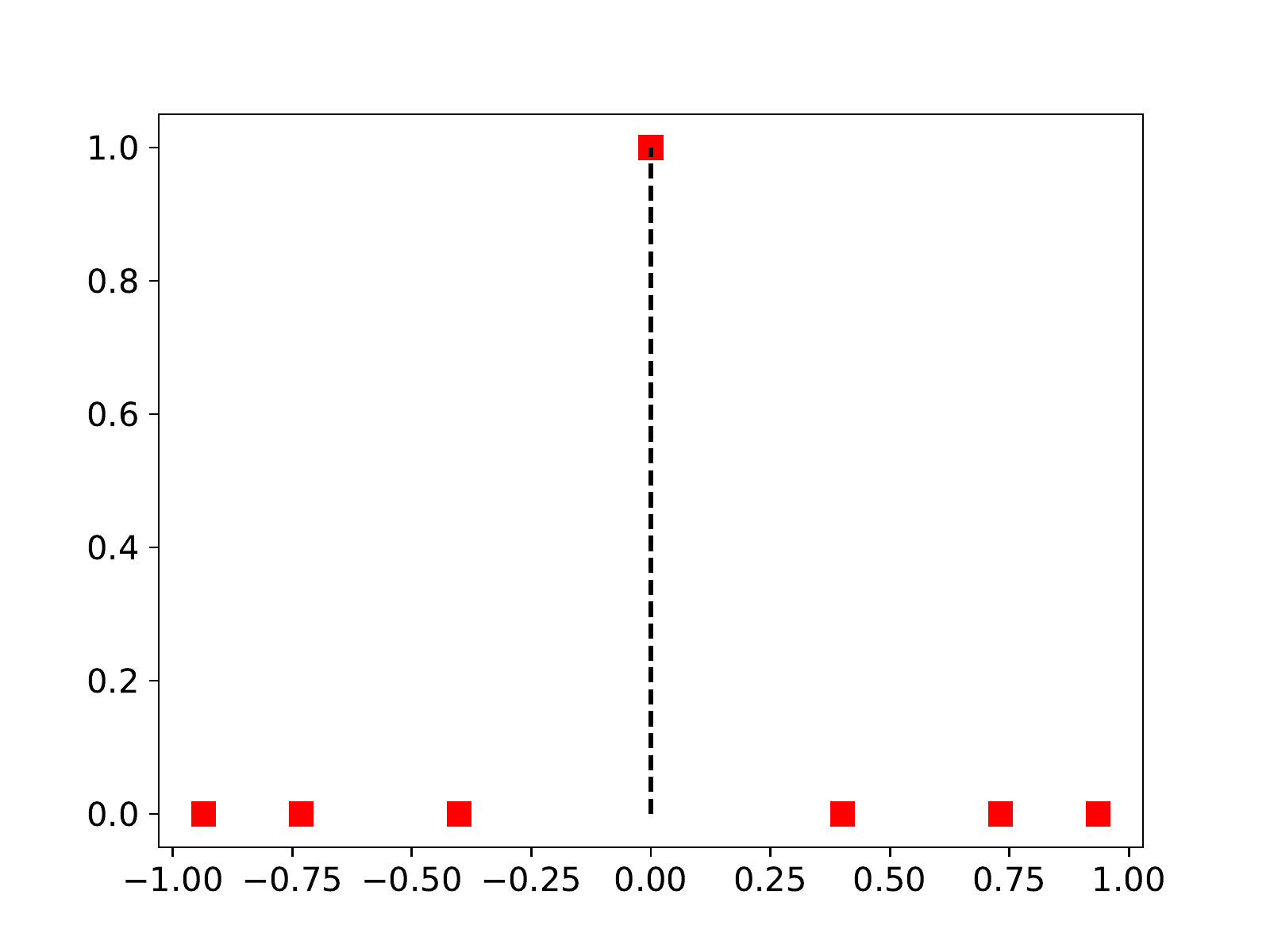}
    \caption{$f_1(\xi)$, PA for $m=1$}
    \label{fig:step_0_m1}
  \end{subfigure}%
  ~
  \begin{subfigure}[b]{0.32\textwidth}
    \includegraphics[width=\textwidth]{%
      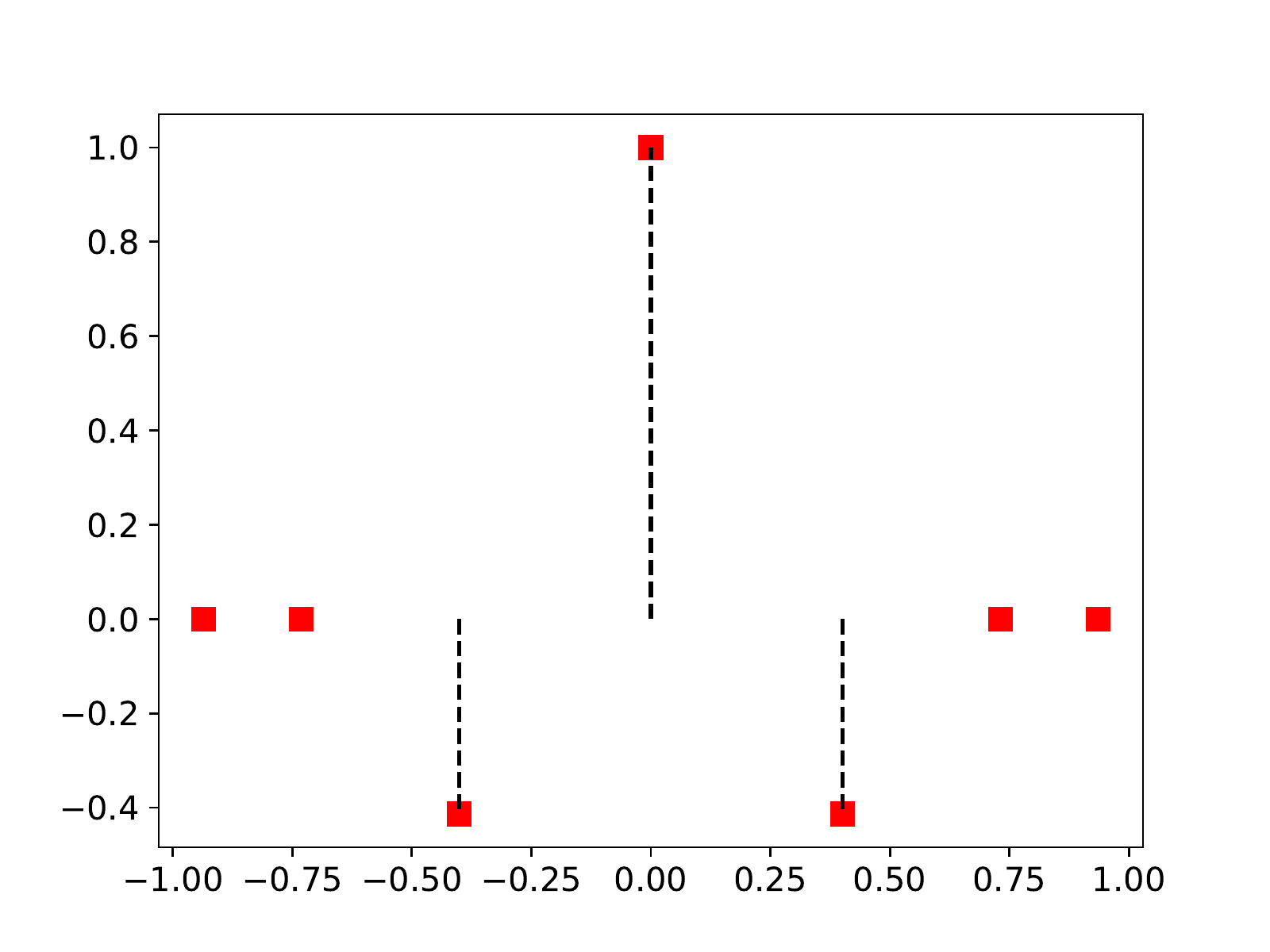}
    \caption{$f_1(\xi)$, PA for $m=3$}
    \label{fig:step_0_m3}
  \end{subfigure}%
  \rulesep
  \begin{subfigure}[b]{0.32\textwidth}
    \includegraphics[width=\textwidth]{%
      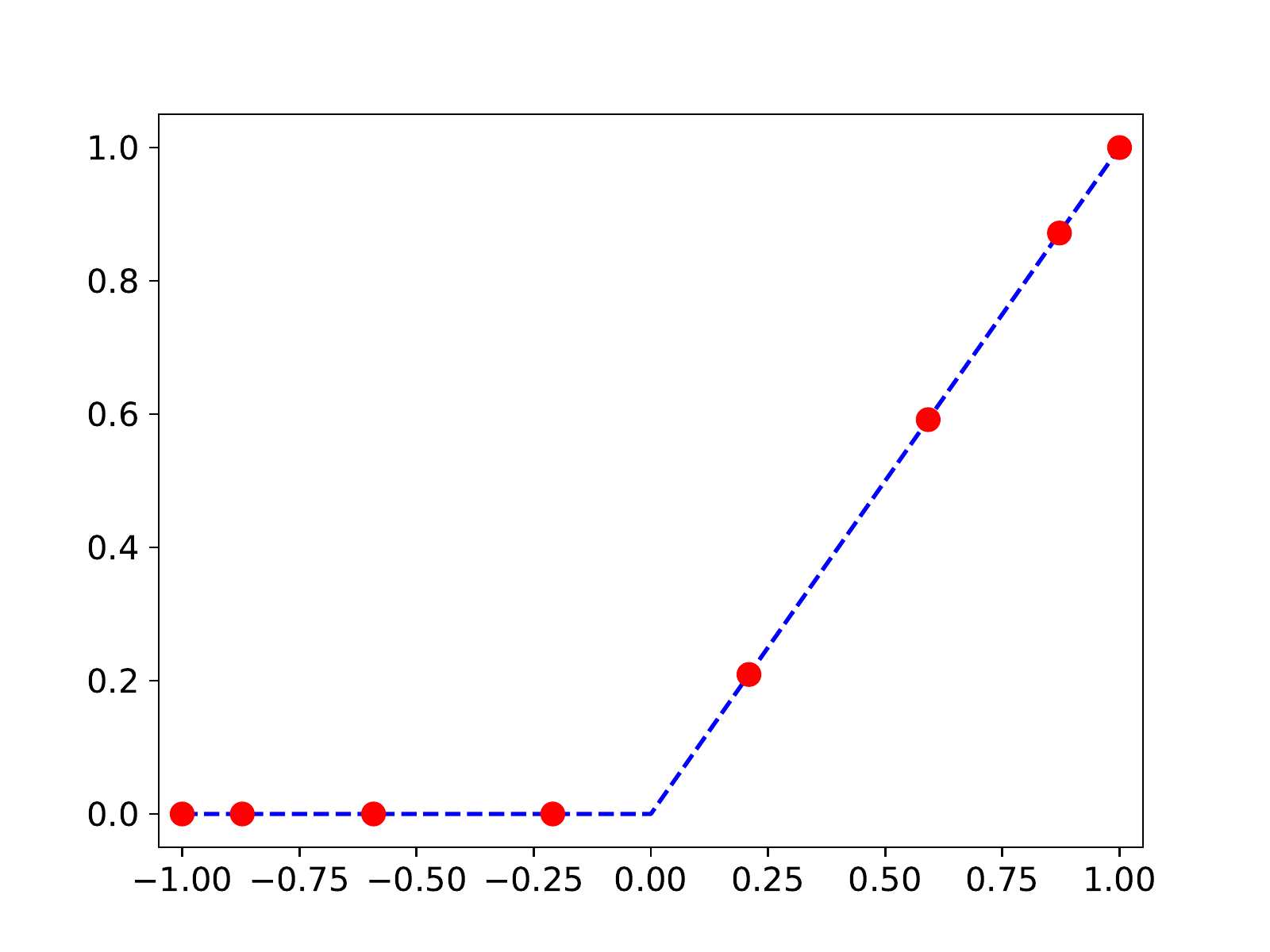}
    \caption{$f_2(\xi)$}
    \label{fig:c0-kink_fun}
  \end{subfigure}%
  ~
  \begin{subfigure}[b]{0.32\textwidth}
    \includegraphics[width=\textwidth]{%
      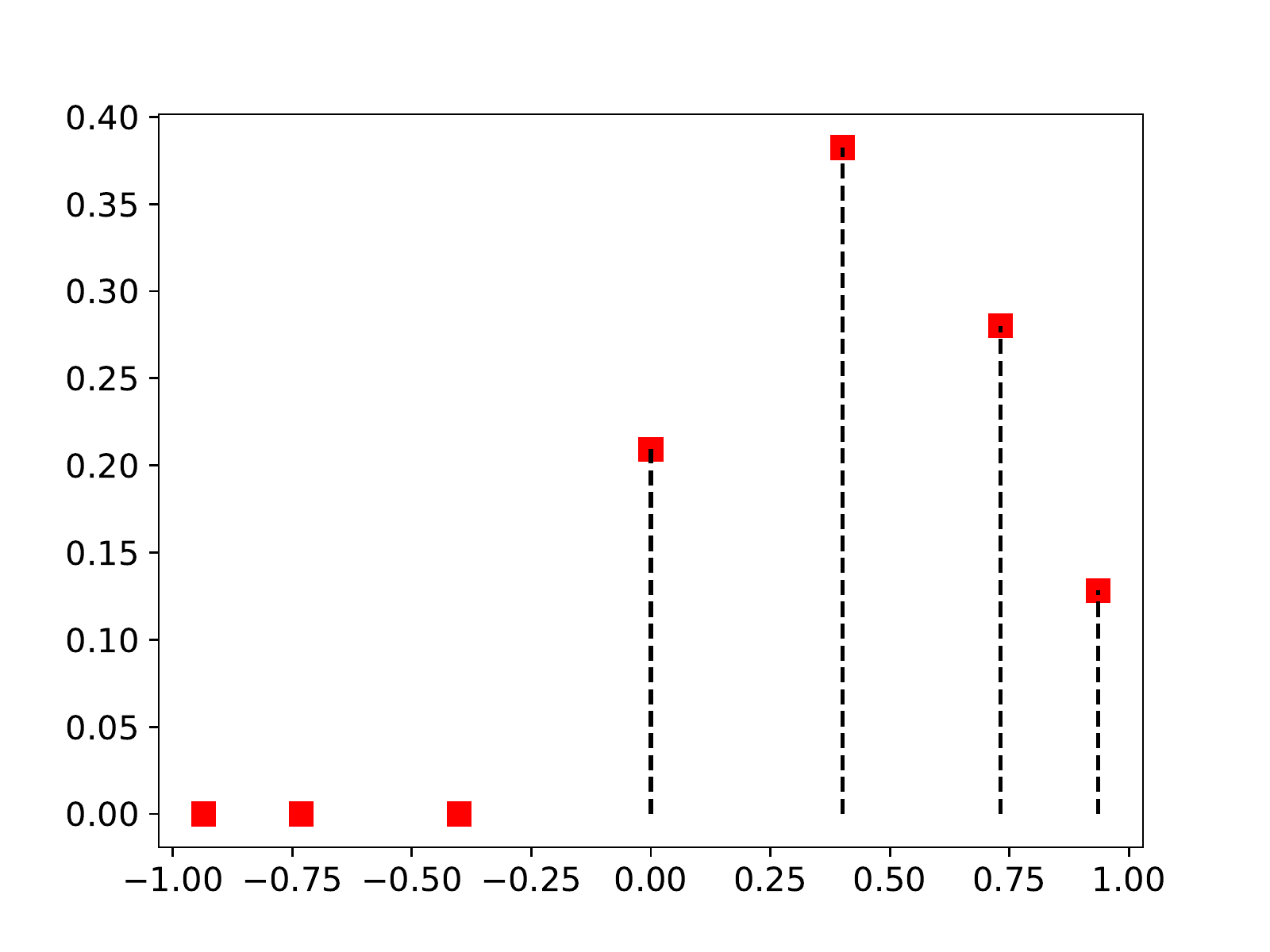}
    \caption{$f_2(\xi)$, PA for $m=1$}
    \label{fig:c0-kink_m1}
  \end{subfigure}%
  ~
  \begin{subfigure}[b]{0.32\textwidth}
    \includegraphics[width=\textwidth]{%
      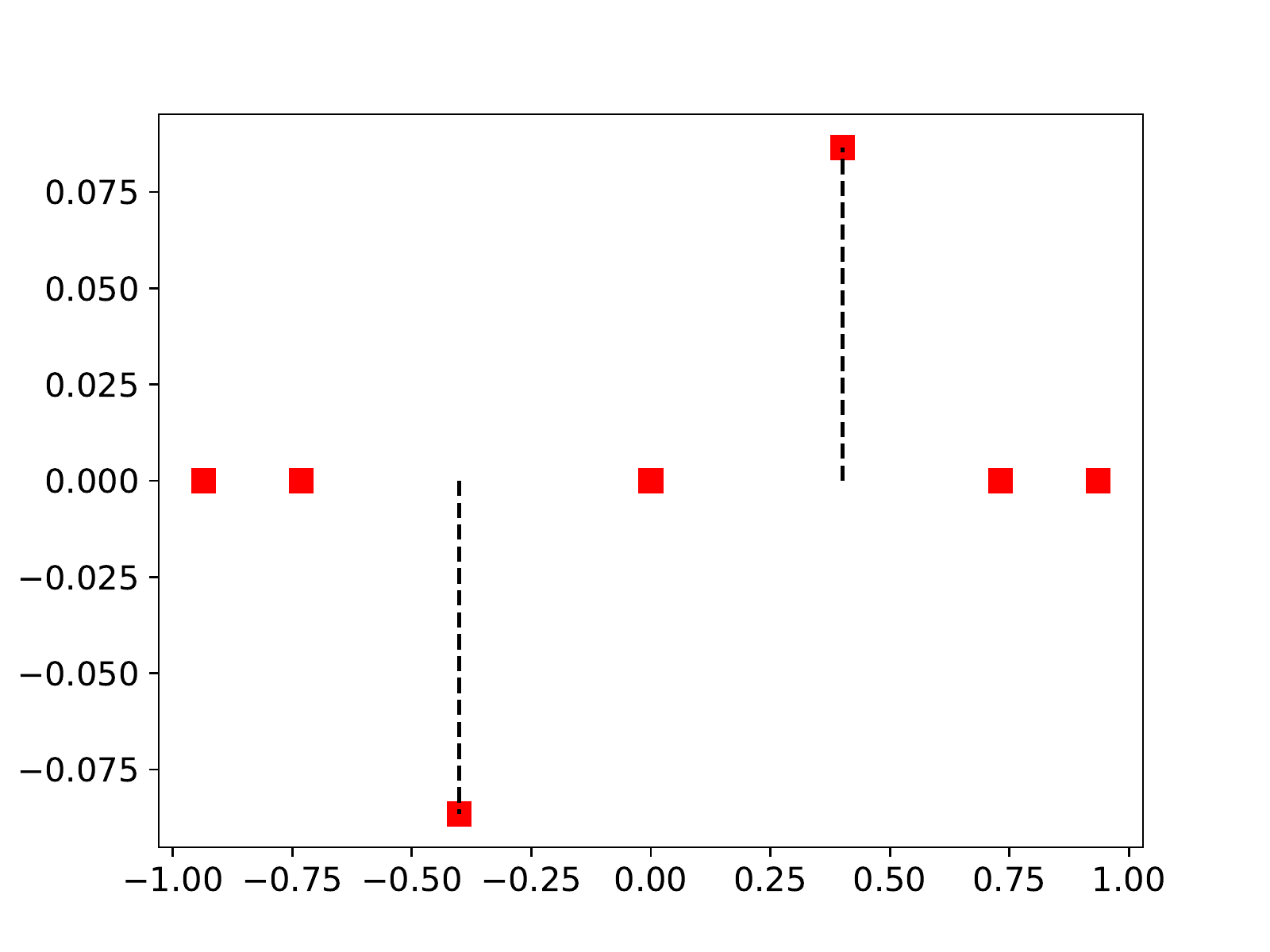}
    \caption{$f_2(\xi)$, PA for $m=3$}
    \label{fig:c0-kink_m3}
  \end{subfigure}%
  \rulesep
  \begin{subfigure}[b]{0.32\textwidth}
    \includegraphics[width=\textwidth]{%
      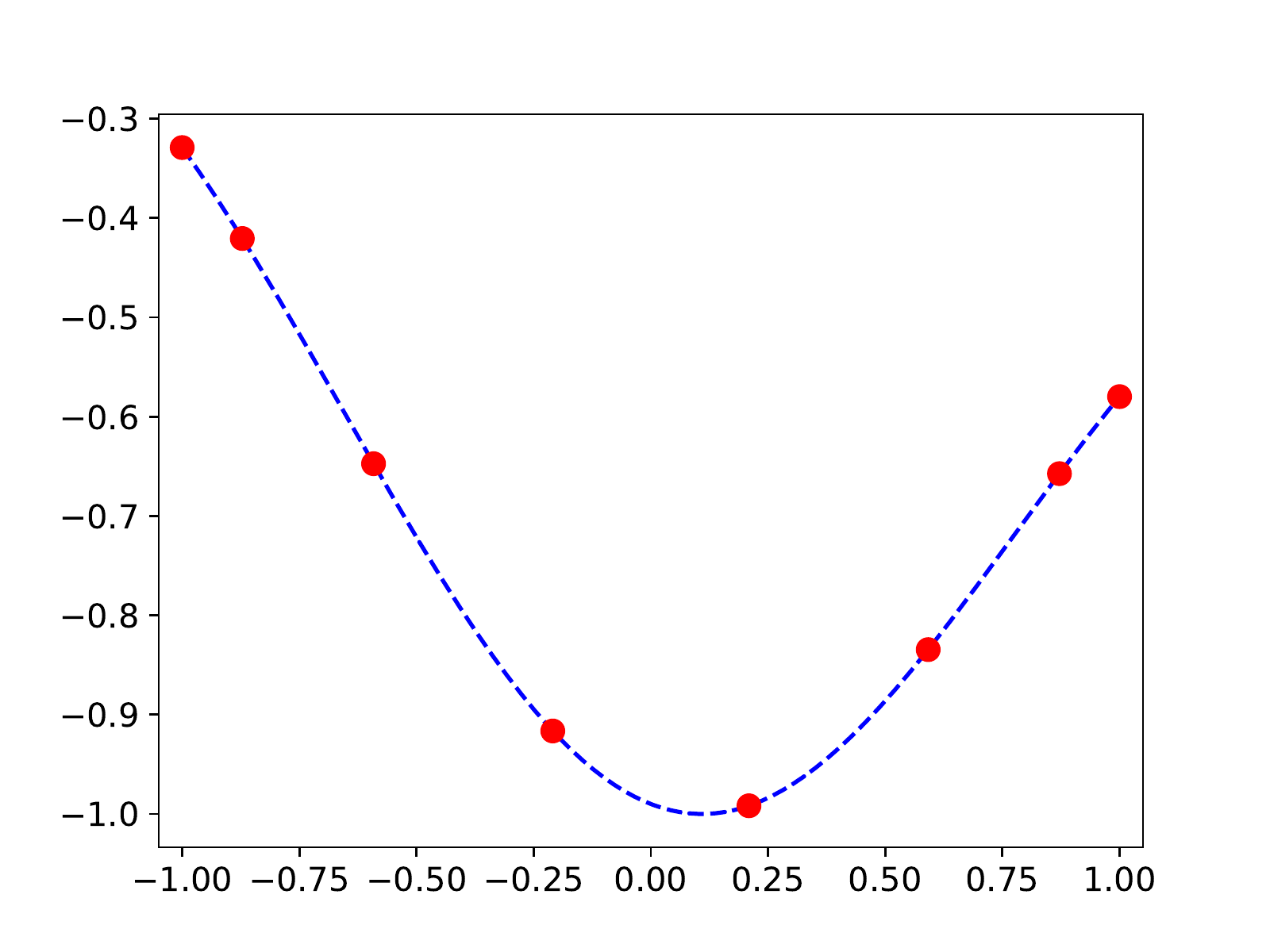}
    \caption{$f_3(\xi)$}
    \label{fig:cos-sin_fun}
  \end{subfigure}%
  ~
  \begin{subfigure}[b]{0.32\textwidth}
    \includegraphics[width=\textwidth]{%
      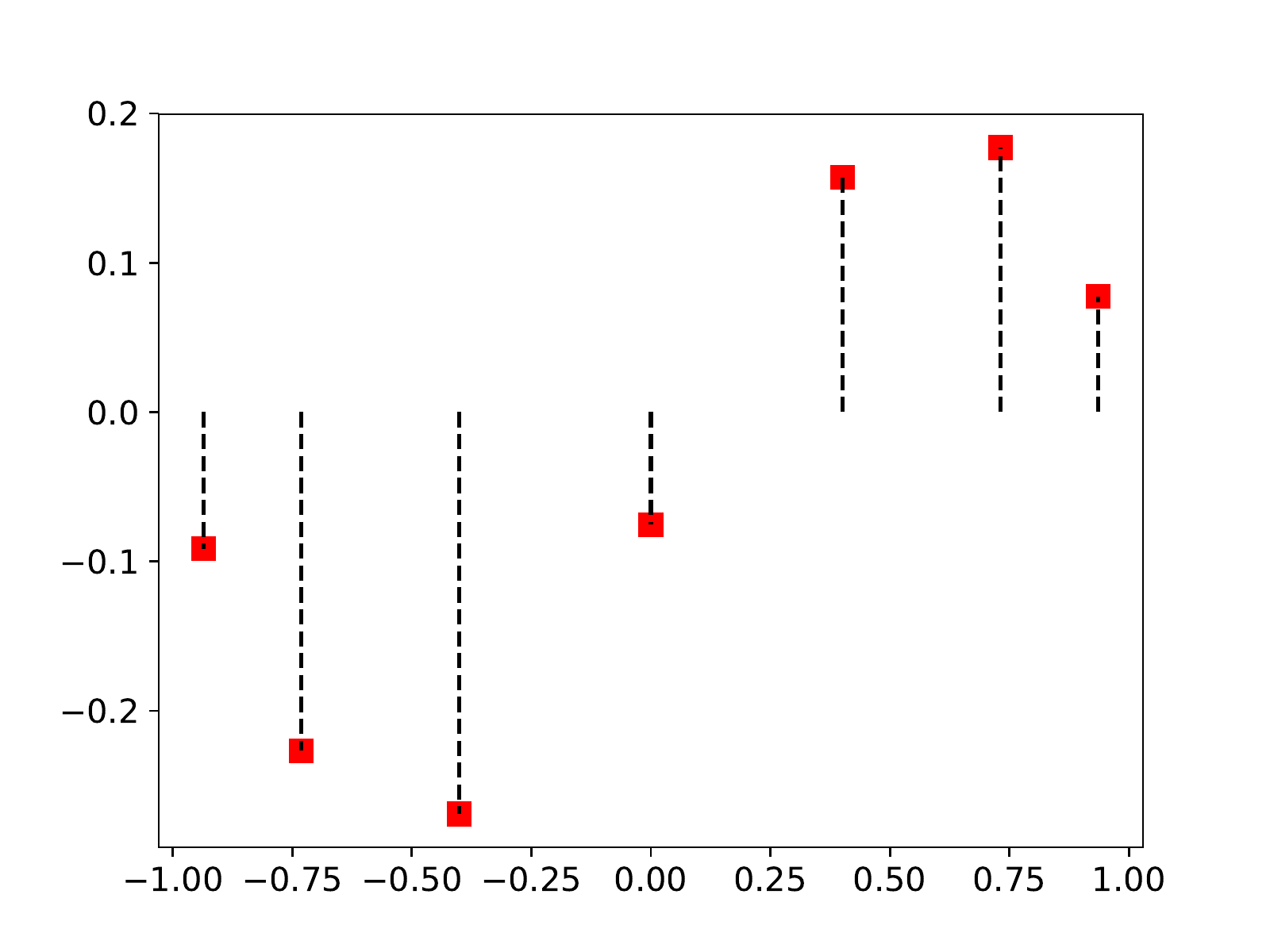}
    \caption{$f_3(\xi)$, PA for $m=1$}
    \label{fig:cos-sin_m1}
  \end{subfigure}%
  ~
  \begin{subfigure}[b]{0.32\textwidth}
    \includegraphics[width=\textwidth]{%
      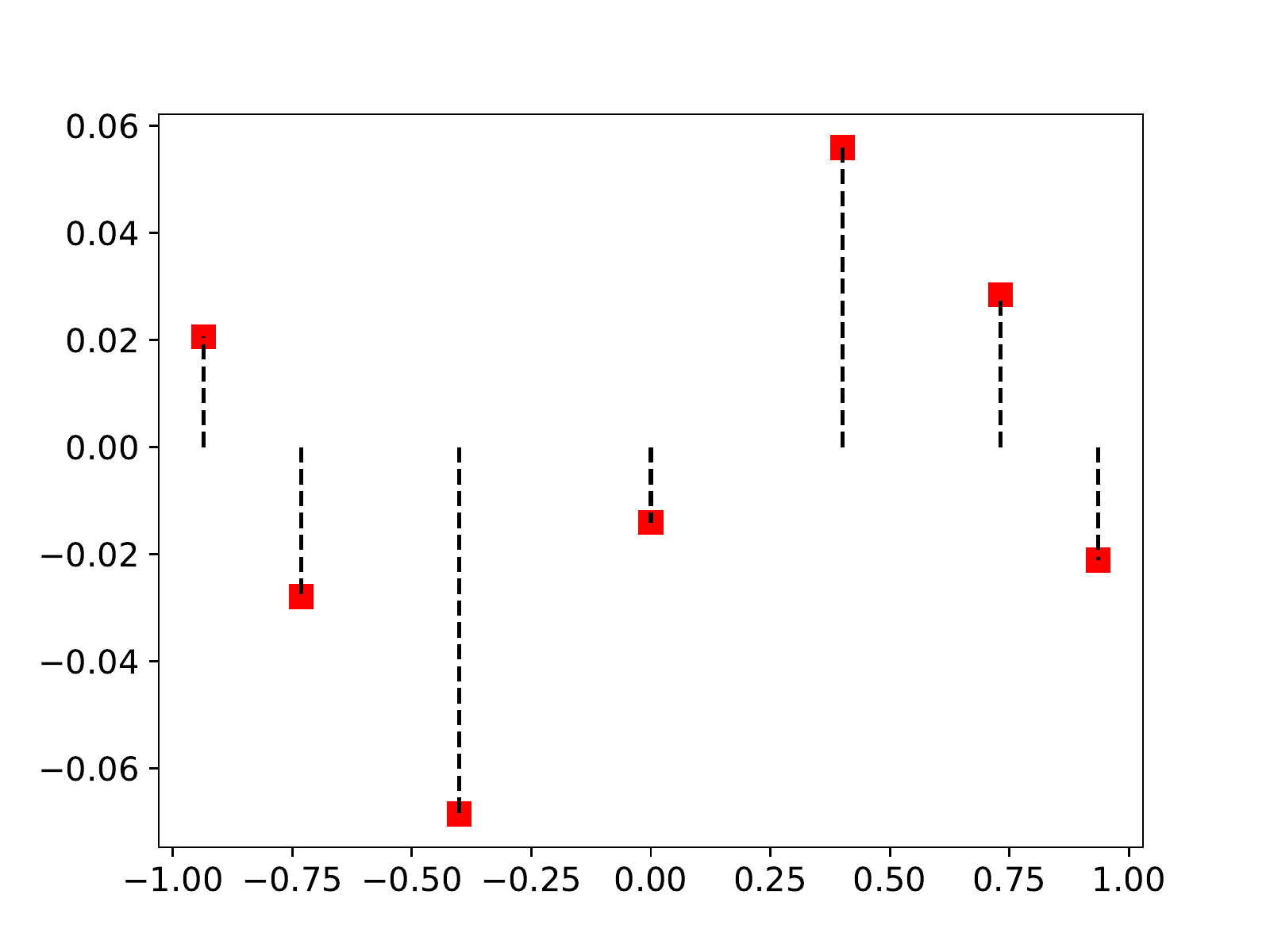}
    \caption{$f_3(\xi)$, PA for $m=3$}
    \label{fig:cos-sin_m3}
  \end{subfigure}%
  \caption{Three different functions $f$ at $8$ Gauss-Lobatto points. 
    The PA operator is evaluated at their mid points.}
  \label{fig:PA}
\end{figure}

Figures \ref{fig:step_0_m1} and \ref{fig:step_0_m3} demonstrate the PA operator for the 
discontinuous function illustrated in Figure \ref{fig:step_0_fun}, 
while Figures \ref{fig:c0-kink_m1} and \ref{fig:c0-kink_m3} do so for the 
continuous but not differentiable function shown in Figure \ref{fig:cos-sin_fun}, 
and Figures \ref{fig:cos-sin_m1} and \ref{fig:cos-sin_m3} illustrate the PA operator for 
the smooth function in Figure \ref{fig:cos-sin_fun}. 

For each function we used eight Gauss-Lobatto points to compute the PA operators for $m=1,3$ at the mid points $\xi_{j+1/2} = \left(\xi_{j-1}+\xi_{j}\right)/2$. 
As illustrated in Figure \ref{fig:step_0_m1} and \ref{fig:step_0_m3}, the jump function 
approximation becomes oscillatory as $m$ is increased from $1$ to $3$, especially in the region of 
the discontinuity. 
The maximal absolute value of the PA operator is also not decreasing in this case. 
On the other hand, for the continuous and smooth functions displayed in Figures 
\ref{fig:c0-kink_fun} and \ref{fig:cos-sin_fun}, the maximal absolute value decreases significantly when going from $m = 1$ to $m = 3$.
This is consistent with the results in \eqref{eq:convergence}, and will be 
used to construct the discontinuity sensor in \S \ref{sub:sensor}.
For a discussion on the convergence of the PA operator see \cite{archibald2005polynomial}.

\begin{remark}
  In this work, we only consider one dimensional conservation laws. 
  It should be stressed, however, that polynomial annihilation can be extended to multivariate 
irregular data in any domain. 
  It was demonstrated in \cite{archibald2005polynomial} that polynomial annihilation 
is numerically cost efficient and entirely independent of any specific shape or complexity of 
boundaries. 
  In particular, in \cite{archibald2009discontinuity} and \cite{jakeman2011characterization} the 
method was applied to high dimensional functions that arise when solving stochastic partial 
differential equations, which reside in a high dimensional space which includes the original space 
and time domains as well as additional random dimensions. 
\end{remark}
 
\section{\texorpdfstring{Application of $\ell_1$ regularization}{Application of l1 regularization}}
\label{sec:l1_limiter}

In this section, we describe how the proposed $\ell^1$ regularization using PA operators, 
i.e., $H\left(v\right) = \norm{ L_m v }_1$ in \eqref{eq:unconst_problem}, can be incorporated into a DG method. 
While this kind of regularization functional was already investigated in \cite{scarnati2017using}, this work is the 
first to extend these ideas to an SE method and thus to allow element-to-element variations in the 
optimization problem.
It should be stressed that the subsequent procedure relies on a piecewise polynomial approximation in space. 
Yet, by appropriate modifications of the procedure, it is also possible to apply $\ell^1$ regularization (with 
PA) to any other method.

\subsection{Procedure}
\label{sub:procedure}

One of the main challenges in solving nonlinear conservation laws \eqref{eq:cl} is balancing high 
resolution properties and the amount of viscosity introduced to maintain stability, 
especially near shocks \cite{vonneumann1950method,persson2006sub,glaubitz2019smooth}. 
Applying the techniques presented in 
\S \ref{sub:l1_reg} and \S \ref{sub:poly_anni}, 
we are now able to adapt the nodal DG method described in \S \ref{sub:DG} to include $\ell^1$ 
regularization. 

Our procedure consists of replacing the usual polynomial approximation $u_p$ by a \emph{sparse reconstruction}  
\begin{equation}\label{eq:op-problem}
    u_p^{\mathrm{spar}} = \argmin_{v \in \mathbb{P}_p\left(\Oref\right)} \left( \frac{1}{2}
\norm{v-u_p}_2^2 + \lambda  \norm{L_m v}_1 \right)
\end{equation}
with \emph{regularization parameter} $\lambda$ in troubled elements after every (or every $k$th) time step by an 
explicit time integrator. 

For the ADMM described in \S \ref{sub:ADMM}, 
it is advantageous to rewrite \eqref{eq:op-problem} in the usual form of an $\ell^1$ regularized 
problem as
\begin{equation}\label{eq:op-problem2}
    u_p^{\mathrm{spar}} = \argmin_{v \in \mathbb{P}_p\left(\Oref\right)} \left( \norm{L_m v}_1 + \frac{\mu}{2} \norm{v-u_p}_2^2 \right), 
\end{equation} 
where $\mu = \frac{2}{\lambda}$ is referred to as the \emph{data fidelity parameter.}
Note that (\ref{eq:op-problem}) and (\ref{eq:op-problem2}) are equivalent.
In the later numerical tests, the data fidelity  parameter $\mu$ and the regularization 
parameter $\lambda$ will be steered by a discontinuity sensor proposed in 
\S \ref{sub:sensor}.

\begin{remark}
\label{rem2}
One of the main drawbacks in using $\ell^1$ regularization for solving numerical partial 
differential equations, as well as for image restoration or sparse signal recovery, is in choosing 
the regularization parameter $\lambda$ (or $\mu$). 
Ideally, one would want to balance the terms in \eqref{eq:op-problem} or \eqref{eq:op-problem2}, 
but this is difficult to do without knowing their comparative size. 
Indeed, the $\ell^1$ regularization term $\norm{L_m v}_1$ heavily depends on the magnitudes of 
nonzero values in the sparsity domain, in this case the jumps. 
Larger jumps are penalized significantly more in the $\ell^1$ norm than smaller values.

Iterative {\em spacially varying} weighted $\ell^1$ regularization techniques (see, e.g., 
\cite{candes2008enhancing,Gelb2018reducing}) are designed to help reduce the size of the norm, since the 
remaining values should be close to zero in magnitude. 
Specifically, the jump discontinuities which are meant to be in the solution can ``pass through'' 
the minimization. 
In this way, with some underlying assumptions made on the accuracy of the fidelity term, one could 
argue that both terms are close to zero. 
Consequently, the choice of $\lambda$ (or $\mu$) should not have as much impact on the results, 
leading to greater robustness overall. 
For the numerical experiments in this investigation, we simply chose regularization parameters that 
worked well. 
We did not attempt to optimize our results and leave parameter selection to future work.
\end{remark}

\subsection{\texorpdfstring{Selection of the regularization parameter $\lambda$}{Selection of the regularization 
parameter}}
\label{sub:selection}

The $\ell^1$ regularization should only be activated in troubled elements. 
In particular, we do not want to unnecessarily degrade the accuracy in the smooth regions of the 
solution.
We thus propose to adapt the regularization parameter $\lambda$ in \eqref{eq:op-problem} to 
appropriately capture different discontinuities and regions of smoothness. 
As a result, the optimization problem will be able to calibrate the resulting sparse reconstruction to the smoothness of 
the solution. 
More specifically, to avoid unnecessary regularization, we choose $\lambda = 0$ in elements 
corresponding to smooth regions. 
Note that this also renders the proposed method more efficient.

On the other hand, when a discontinuity is detected in an element, $\ell^1$ regularization will be fully activated by 
choosing $\lambda = \lambda_{\text{max}}$ in \eqref{eq:op-problem}, 
which corresponds to the amount of regularization necessary to reconstruct sharp shock profiles. 
While no effort was made to optimize or even adapt this parameter, we found that using 
$\lambda_{\text{max}} = 4 \cdot 10^2$ in all of our numerical experiments yielded  good results. 
A heuristic explanation for choosing $\lambda_{\text{max}}$ in this way stems from the goal of 
balancing the size of $\norm{L_m v}_1$ with the expected size of the fidelity term, which in this case means to be
consistent with the order of accuracy of the underlying numerical PDE solver. 
As mentioned previously, choosing an appropriate $\lambda$ will be the subject of future work. 

Between these extreme cases, i.e., $\lambda = 0$ and $\lambda = \lambda_{\text{max}}$, 
we allow the regularization parameter to linearly vary and choose $\lambda$ as a function of 
the discontinuity sensor proposed in \S \ref{sub:sensor}.
As a consequence, we obtain more accurate sparse reconstructions while still maintaining stability 
in regions around discontinuities.

\subsection{Discontinuity sensor}
\label{sub:sensor}

We now describe the discontinuity sensor which is used to activate the $\ell^1$ regularization and to steer the 
regularization parameter $\lambda$ in \eqref{eq:op-problem}. 
The sensor is based on comparing PA operators of increasing orders. 
To the best of our knowledge this is the first time the PA operator is utilized for shock 
(discontinuity) detection in a PDE solver.\footnote{Of course (W)ENO schemes 
\cite{qiu2005comparison,shu1988efficient,shu1989efficient} compare slope magnitudes for determining 
troubled elements and choosing approximation stencils, so in this regard our method was inspired by 
WENO type methods.}

At least for smooth solutions, discontinuous Galerkin methods are capable of spectral 
orders of accuracy.
$\ell^1$ regularization as well as any other shock capturing procedure 
\cite{vonneumann1950method,persson2006sub,glaubitz2018application,ranocha2018stability,glaubitz2019smooth} 
should thus be just applied in (and near) elements where discontinuities are present. 
We refer to those elements as \emph{troubled} elements. 

Many shock and discontinuity sensors have been proposed over the last $20$ years for the selective
application of shock capturing methods. 
Some of them use information about the $L^2$-norm of the residual of the variational form 
\cite{bassi1995accurate,jaffre1995convergence}, the primary orientation of the
discontinuity \cite{hartmann2006adaptive}, the magnitude of the facial interelement jumps 
\cite{barter2010shock,feistauer2007robust}, or entropy pairs \cite{guermond2008entropy} to detect troubled elements. 
Others are not just able to detect troubled elements, but even the location of discontinuities in the 
corresponding element, such as the concentration method in 
\cite{gelb1999detection,gelb2000detection,gelb2006adaptive,gelb2008detection,offner2013detecting,don2016hybrid,TadmorWaagan2012}. 
The PA operator is also capable of detecting the location of a discontinuity 
\cite{archibald2005polynomial,archibald2008determining,archibald2009discontinuity} 
up to the resolution size, as in (\ref{eq:convergence}), using nonlinear postprocessing 
techniques.
It should be stressed, however, that in this work we do not fully make use of this feature. 

\begin{figure}[!htb]
  \centering
  \begin{subfigure}[b]{0.95\textwidth}
    \includegraphics[width=\textwidth]{%
      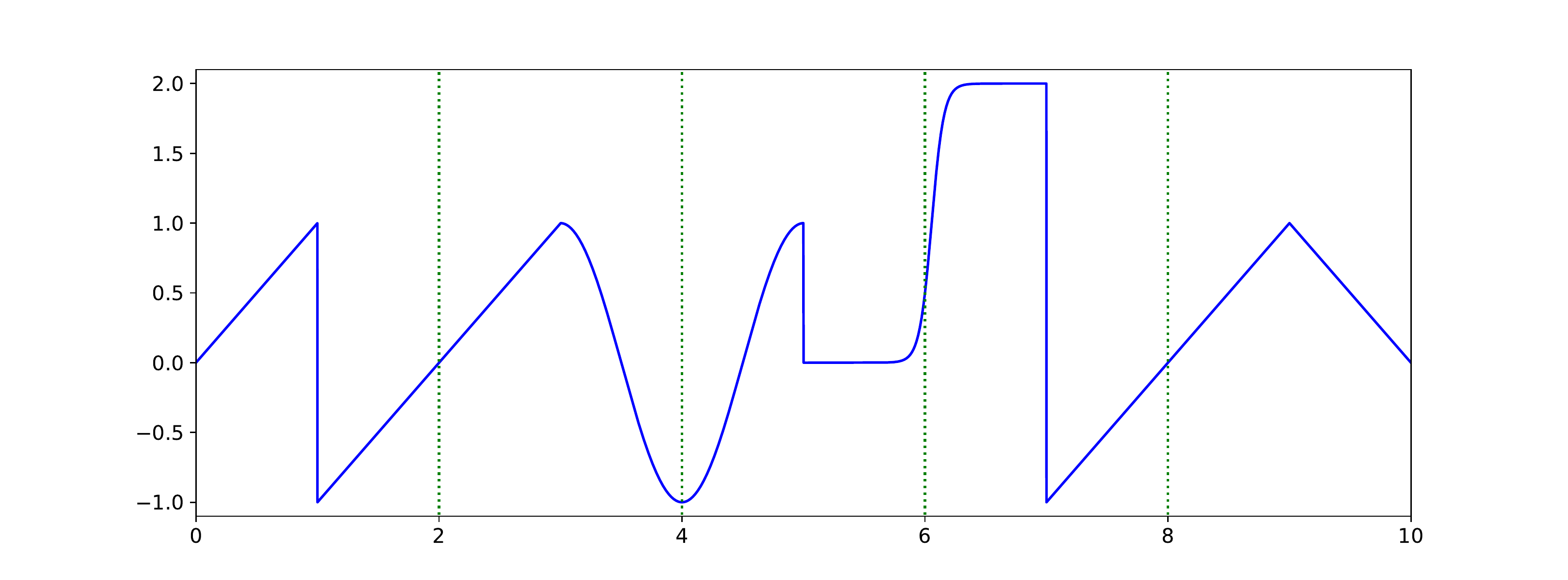}
  \end{subfigure}%
  \caption{Test function with various different discontinuities as well as continuous but not differentiable points.}
  \label{fig:test-fun}
\end{figure}

In what follows we present a novel discontinuity sensor based on the PA operator. 
Let the \emph{sensor value of order $m$} be
\begin{equation}\label{eq:switchS}
    S_m = \max_{0 \leq k \leq p-1} \left| L_m[u_p]\left(\xi_{k+\frac{1}{2}}\right) \right|,
\end{equation} 
i.e., the greatest absolute value of the PA operator at the midpoints 
\begin{equation}
    \xi_{k+\frac{1}{2}} = \frac{\xi_k + \xi_{k+1}}{2}, 
    \quad \text{for } k=0,\dots,p-1,
\end{equation}
of the collocation points. 
Since $L_m$ provides convergence to $0$ of order $m$ in elements where $u$ has $m$ 
continuous derivatives, we expect $S_3 < S_1$ to hold for an at least continuous function. 
Thus, $\ell^1$ regularization should just be fully activated if 
\begin{equation}\label{eq:sensor}
    S(u) := \frac{S_3}{S_1} \geq 1,
\end{equation} 
i.e., if the sensor value does not decrease as $m$ increases from $1$ to $3$.
Sensor values of order $1$ and $3$ (rather than $2$) are chosen since having symmetry of the grid points in 
$S_\xi$ surrounding the point $\xi$ yields a simpler form for implementation \cite{archibald2005polynomial}. 
Various modifications of the \emph{PA sensor} \eqref{eq:sensor} are possible and will be the topic of future research. 
In the following, the resulting PA sensor is demonstrated for the function displayed in Figure \ref{fig:test-fun} on 
the interval $[0,10]$.  We decompose the interval into $I=5$ elements and apply the PA operator and resulting PA sensor 
separately on each element.

Table \ref{tab:sensor} lists the sensor values of order $m=1$ and $m=3$ on each element in Figure \ref{fig:test-fun}. 
The last row further shows if the PA sensor \eqref{eq:sensor} reacts. 
As can be seen in Table \ref{tab:sensor}, the discontinuity sensor exactly identifies the troubled elements, simply by 
comparing the sensor values of order $m=1$ and $m=3$.

\begin{table}[!htb] 
\resizebox{\columnwidth}{!}{%
    \centering
    \begin{tabular}{ c | c | c | c | c | c }
        \hline
        Element & $\Omega_1=\left(0,2\right)$ & $\Omega_2=\left(2,4\right)$ & $\Omega_3=\left(4,6\right)$ & 
$\Omega_4=\left(6,8\right)$ & $\Omega_5=\left(8,10\right)$ \\
        \hline
        $S_1$ & 1.58 & 1.07 & 1.07 & 2.79 & 0.38 \\
        $S_3$ & 1.99 & 0.60 & 1.38 & 3.00 & 0.17 \\ 
        \hline
        Discontinuity & yes & no & yes & yes & no \\
        \hline
    \end{tabular} 
}
    \caption{Shock sensor for function displayed in Figure \ref{fig:test-fun}.}
    \label{tab:sensor}
\end{table}

Even though we use this sensor only to detect troubled elements, in principle, 
this new sensor can also be used to more precisely determine the location and strength of the 
discontinuities. 
Such information can be used for refined domain decomposition 
\cite{archibald2009discontinuity}. 
Finally, we decide for the regularization parameter $\lambda$ to linearly vary between $\lambda=0$ and 
$\lambda=\lambda_{\text{max}}$ and thus utilize the parameter function 
\begin{equation}\label{eq:param-fun}
    \lambda(S) = 
    \left\{ 
    \begin{array}{ll}
        0 & \text{if } S \leq \kappa, \\ 
        \lambda_{\text{max}} (S-\kappa)/(1-\kappa) & \text{if } \kappa < S < 1, \\ 
        \lambda_{\text{max}} & \text{if } 1 \leq S, \\ 
    \end{array} \right.
\end{equation} 
where $\kappa \in [0,1)$ is a problem dependent \emph{ramp parameter}. 
Observe that using $\lambda(S)$ is comparable to employing the weighted $\ell_1$ 
regularization as discussed in Remark \ref{rem2}.

For the later numerical tests we also considered other parameter functions, some as discussed in \cite{huerta2012simple}. 
Yet the best results were obtained with \eqref{eq:param-fun}. 
The same holds for other discontinuity sensors, such as the modal-decay based sensor of Persson and Peraire \cite{persson2006sub,huerta2012simple} and its refinements \cite{barter2010shock,klockner2011viscous} as well as the KXRCF sensor \cite{krivodonova2004shock,qiu2005comparison} of Krivodonova et al., which is built up on a strong superconvergence phenomenon of the DG method at outflow boundaries.  For brevity, those results are omitted here.

\begin{remark}
\label{rem3}
  We note that the PA sensor might produce false positive or false negative misidentifications
in certain cases.
  A false negative misidentification might arise from a discontinuity where the solution is 
detected to be smooth. 
  This is encountered by the ramp parameter $\kappa$, which is observed to work robustly for 
$\kappa = 0.8$ or $\kappa = 0.9$ in all later numerical tests.
  A false positive misidentification might arise from a smooth solution which is detected to be 
nonsmooth (possibly discontinuous). 
  In this case smooth parts of the solution with steep gradients will result in significantly 
greater values of $L_mu$ than parts of the solution with less steep gradients. 
  As a result, the standard $\ell^1$ regularization would heavily penalize these features of the 
solution, yielding inappropriate smearing of steep gradients in smooth regions. 
  By making $\lambda = \lambda(S)$ dependent on the sensor value, this problem can be somewhat 
alleviated. 
  Using a weighted $\ell^1$ regularization, as suggested in Remark \ref{rem2}, should also 
reduce the unwanted smearing effect.
  Failure to detect a discontinuity would, after a number of time steps, yield instability. 
  However it is unlikely that this would occur as the growing oscillations would more likely be 
identified as shock discontinuities.
\end{remark}

\subsection{Efficient implementation of the PA operator}
\label{sub:implementation} 

While the PA operator was defined on the interior of the reference element $\Oref=[-1,1]$ in \S \ref{sub:poly_anni}, 
the shock sensor proposed in \S \ref{sub:sensor} only relies on the values of the PA operator at the $p$ 
midpoints $\{ \xi_{k+\frac{1}{2}} \}_{k=0}^{p-1}$ of the collocation points $\{ \xi_k \}_{k=0}^p$. 
The same holds for the $\ell^1$ regularization term $H\left(v\right) = \norm{L_mv}_1$. 
The $\ell^1$-norm of the PA transformation is thus given by  
\begin{equation}\label{eq:1-norm}
    \norm{ L_m v }_1 
        = \norm{\mat{L_m} \vec{v}}_1
        = \sum_{k=0}^{p-1} \left| L_m[v]\left(\xi_{k+1/2}\right) \right|, 
\end{equation} 
where the vector $\vec{v}$ once more consists of nodal degrees of freedom. 
We now aim to provide an efficient implementation of the PA operator $L_m$ in form of a matrix representation 
$\mat{L_m}$, which maps the nodal values $\vec{v}$ to the values of the PA operator at the midpoints. 
Revisiting 
\eqref{eq:PA-operator}, this matrix representation is given by 
\begin{equation}
    \mat{L_m} = \mat{Q}^{-1} \mat{C} \in \R^{p \times (p+1)}
\end{equation} 
with 
\begin{align}
    \mat{Q} 
      & = \diag{ \left[ q_m\left(\xi_{1/2}\right), \dots, q_m\left(\xi_{p-1/2}\right) \right] }, \\ 
    q_m\left(\xi_{k+1/2}\right) 
      & = \sum_{x_j \in S_{\xi_{k+1/2}}^+} c_j\left(\xi_{k+1/2}\right) 
\end{align} 
and 
\begin{equation}
    \mat{C} = \left( c_j \left( \xi_{k + 1/2} \right) \right)_{k,j=0}^{p-1,p},
\end{equation}    
where 
\begin{align}
    c_j\left(\xi_{k+1/2}\right) & = 
    \left\{
    \begin{array}{ll}
        \frac{m!}{\omega_j\left(S_{\xi_{k+1/2}}\right)} & \text{if } x_j \in S_{\xi_{k+1/2}}, \\ 
        0 & \text{else},
    \end{array} \right. \\
    \omega_j\left(S_{\xi_{k+1/2}}\right) 
        & = \prod_{\begin{smallmatrix} x_i \in S_{\xi_{k+1/2}} \\i\neq j\end{smallmatrix}} \left(x_j- x_i\right).
\end{align}
Utilizing all prior matrix vector representations, we can now give the discretization of the $\ell^1$ regularized 
optimization problem \eqref{eq:op-problem2} by  
\begin{equation}\label{eq:op-problem-disc}
    \vec{u}^{\mathrm{spar}} 
        = \argmin_{\vec{v} \in \R^{p+1} } \left( \norm{ \mat{L_m} \vec{v} }_1 
        + \frac{\mu}{2} \norm{ \vec{v} - \vec{u} }_{\mat{M}}^2 \right). 
\end{equation} 
Thus, we are able to solve the optimization problem directly for the nodal degrees of freedom of 
the sparse reconstruction $u_p^{\mathrm{spar}}$. 
Alternatively, the fidelity term $\norm{v-u_p}$ can also be approximated as 
\begin{equation}
    \norm{ \vec{v}-\vec{u} }^2_2 
        = \sum_{k=0}^p \left| v(\xi_k) - u(\xi_k) \right|^2, 
\end{equation}
i.e., by the Euclidean norm, yielding 
\begin{equation}\label{eq:op-problem-disc2}
    \vec{u}^{\mathrm{spar}} 
        = \argmin_{\vec{v} \in \R^{p+1} } \left( \norm{ \mat{L_m} \vec{v} }_1 
        + \frac{\mu}{2} \norm{ \vec{v} - \vec{u} }_2^2 \right), 
\end{equation} 
instead of \eqref{eq:op-problem-disc}. 
Future works will investigate the influence of the choice of the discrete norm on the 
performance of the $\ell^1$ regularization. 
Here we decided to use the Euclidean norm and thus the minimization problem 
\eqref{eq:op-problem-disc2}, making the computations in \S \ref{sub:ADMM} more intelligible.

\subsection{The alternating direction method of multipliers}
\label{sub:ADMM} 

\begin{figure}[!htb]
  \centering
   \begin{subfigure}[b]{0.24\textwidth}
    \includegraphics[width=\textwidth]{%
      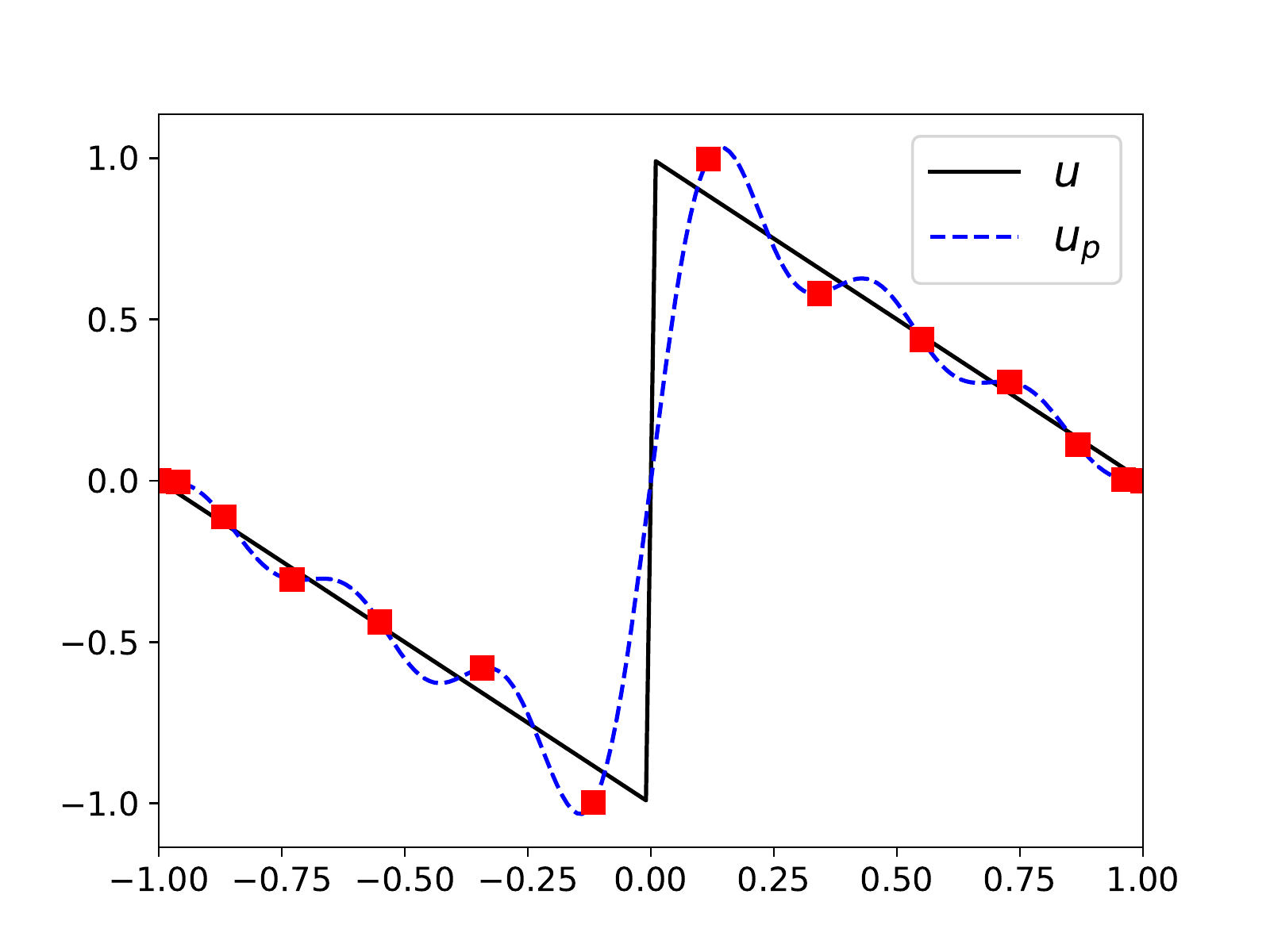}
    \caption{$u_p$}
    \label{fig:fun-values} 
  \end{subfigure}%
  ~
  \begin{subfigure}[b]{0.24\textwidth}
    \includegraphics[width=\textwidth]{%
      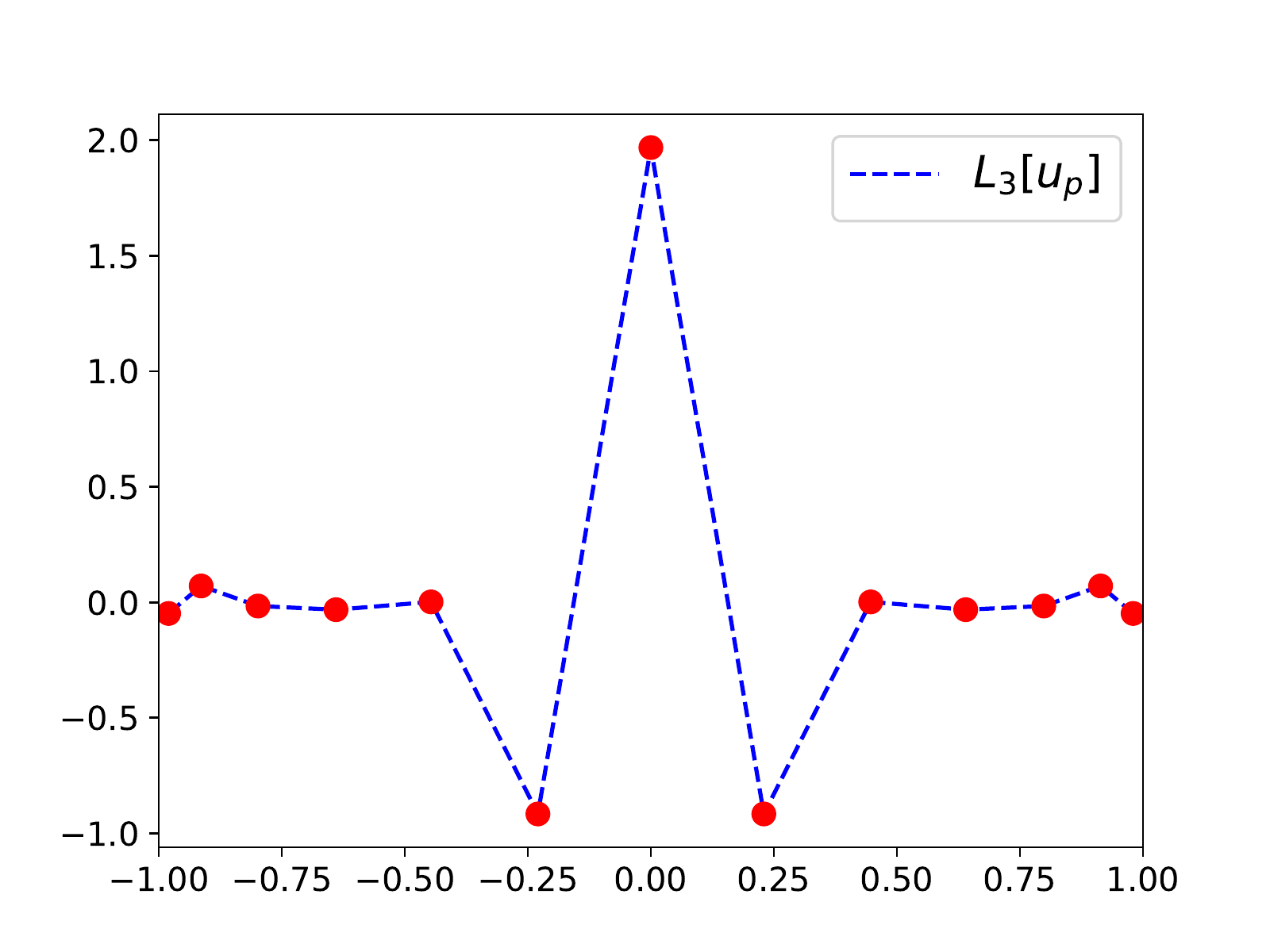} 
    \caption{$L_3[u_p]$}
    \label{fig:jump-values}
  \end{subfigure}%
  \begin{subfigure}[b]{0.24\textwidth}
    \includegraphics[width=\textwidth]{%
      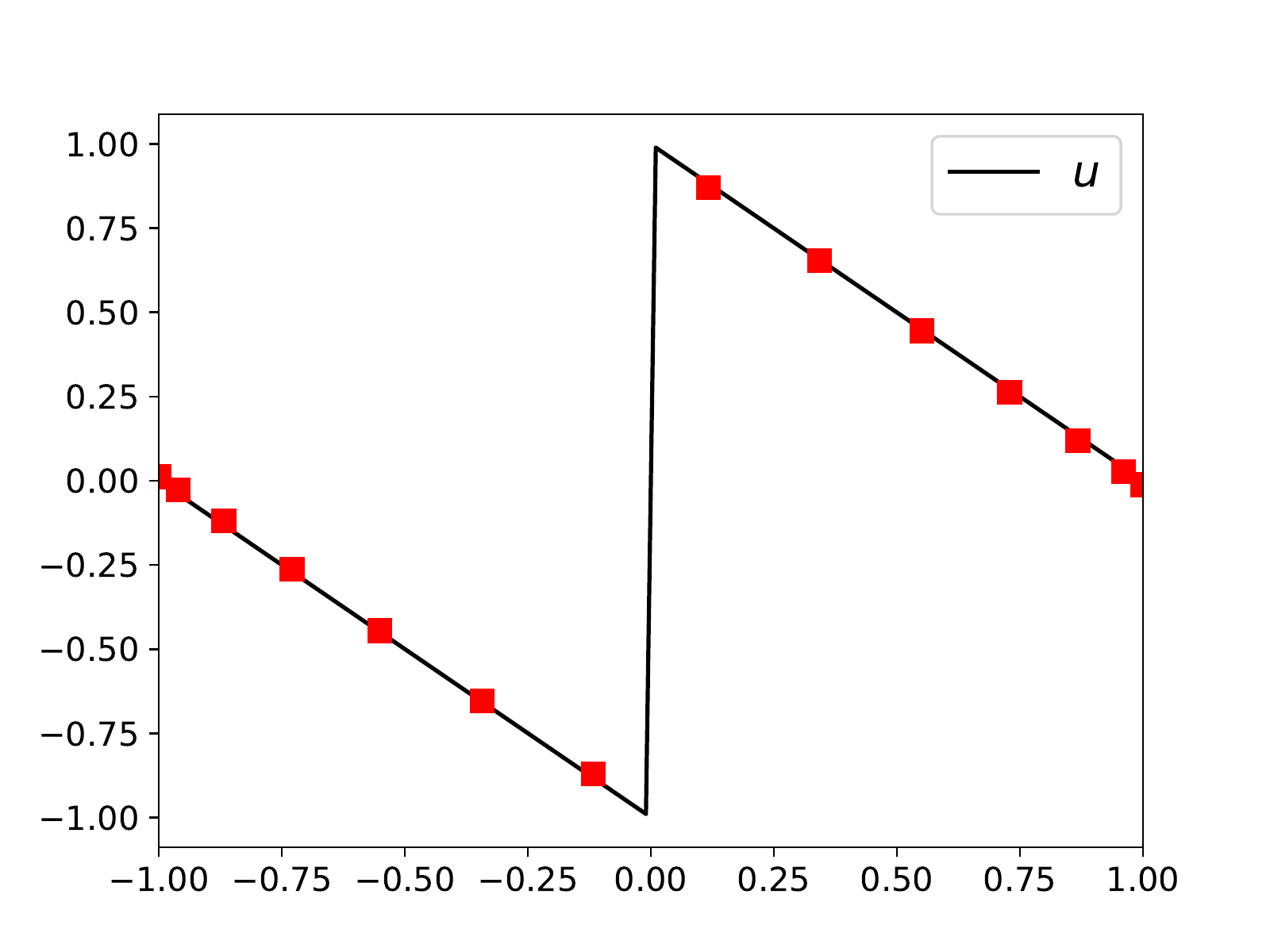} 
    \caption{$u_p^{spar}$}
    \label{fig:reg-fun-values}
  \end{subfigure}%
  ~
  \begin{subfigure}[b]{0.24\textwidth}
    \includegraphics[width=\textwidth]{%
      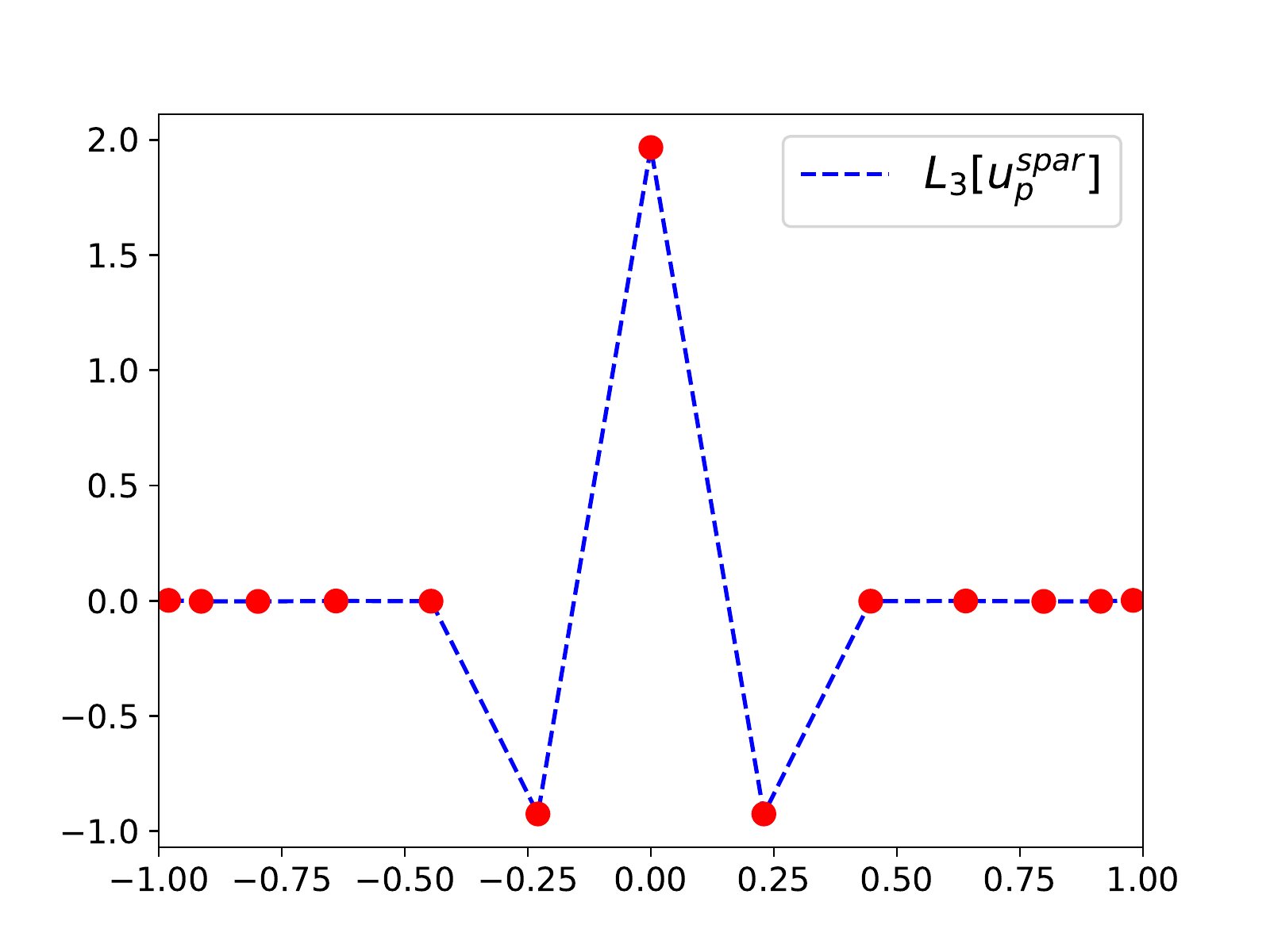} 
    \caption{$L_3[u_p^{spar}]$}
    \label{fig:reg-jump-values} 
  \end{subfigure}%
  \caption{Demonstration of $\ell^1$ regularization for the polynomial approximation (polluted by the Gibbs phenomenon) 
of a discontinuous function $u$.}
  \label{fig:l1-reg}
\end{figure}

Many techniques have been recently proposed to solve optimization problems in the 
form of \eqref{eq:op-problem-disc}. 
Following \cite{scarnati2017using}, we use the ADMM
\cite{li2013efficient,sanders2016matlab,sanders2017composite} in our implementation. 
The ADMM has its roots in \cite{glowinski1975approximation} 
and details of its convergence properties  can be found in 
\cite{gabay1976dual,glowinski1989augmented,eckstein1992douglas}.
In the context of $\ell^1$ regularization, ADMM is commonly implemented using the split Bregman method 
\cite{goldstein2009split}, which is known to be an efficient solver for a broad class of 
optimization problems. 
To implement the ADMM, it is first necessary to eliminate all nonlinear terms within the 
$\ell^1$-norm. 
We thus introduce a slack variable 
\begin{equation}
    \vec{g} = \mat{L_m} \vec{v} \in \R^p
\end{equation} 
and formulate \eqref{eq:op-problem-disc} equivalently as 
\begin{equation}\label{eq:op-problem3}
    \argmin_{\vec{v} \in \R^{p+1}, \vec{g} \in \R^p } \left( \norm{\vec{g}}_1
    + \frac{\mu}{2} \norm{ \vec{v} - \vec{u} }_{2}^2 
    \quad \text{s.t.} \quad \mat{L_m} \vec{v} = \vec{g} \right). 
\end{equation} 
To solve \eqref{eq:op-problem3}, we further introduce Lagrangian multipliers 
$\vec{\sigma} \in \R^p, \vec{\delta} \in \R^{p+1}$ and solve the 
unconstrained minimization problem given by
\begin{equation}\label{eq:op-problem4}
    \argmin_{\vec{v} \in \R^{p+1}, \vec{g} \in \R^p } J_{\vec{\sigma},\vec{\delta}}\left( \vec{v}, \vec{g} \right)
\end{equation} 
with objective function 
\begin{equation}
    J_{\vec{\sigma},\vec{\delta}}\left( \vec{v}, \vec{g} \right) 
        = \norm{\vec{g}}_1 
        + \frac{\mu}{2} \norm{ \vec{v} - \vec{u} }_{2}^2 
        + \frac{\beta}{2} \norm{ \mat{L_m} \vec{v} - \vec{g} }_{2}^2 
        - \scp{ \mat{L_m} \vec{v} - \vec{g} }{ \vec{\sigma} }_{2}
        - \scp{ \vec{v} - \vec{u} }{ \vec{\delta} }_{2}. 
\end{equation} 
Here, $\beta > 0$ is an additional positive regularization parameter
and recall that the data fidelity parameter $\mu$ is given by $\mu = \frac{2}{\lambda}$ for 
$\lambda > 0$; see \eqref{eq:op-problem} and \eqref{eq:op-problem2}.
Note that if the Lagragian multipliers $\vec{\sigma}, \vec{\delta}$ are updated a sufficient 
number of times, the solution of the unconstrained problem \eqref{eq:op-problem4} 
will converge to the solution of the constrained problem \eqref{eq:op-problem3}. 
In the ADMM, the solution is approximated by alternating between minimizations of $\vec{v}$ and $\vec{g}$. 
A crucial advantage of this method is that, given the current value of $\vec{v}$ as well as the Lagrangian multipliers, the optimal value of \vec{g} can be exactly determined by the 
\emph{shrinkage}-like formula \cite{goldstein2009split} 
\begin{equation}\label{eq:update-g}
    \left(\vec{g}_{k+1}\right)_i = shrink\left( \left(\mat{L_m} \vec{v}\right)_i - \frac{1}{\beta} \left(\vec{\sigma}_k\right)_i, \frac{1}{\beta} \right),
\end{equation} 
where 
\begin{equation}\label{eq:shrink}
    shrink\left( x, \gamma \right) = \frac{x}{|x|}\cdot \max\left( |x|-\gamma,0 \right).
\end{equation} 
Given the current value $\vec{g}_{k+1}$, on the other hand, the optimal value of $\vec{v}$ is computed by the gradient descent method as 
\begin{equation}\label{eq:update-v}
    \vec{v}_{k+1} = \vec{v}_k - \alpha \nabla_{\vec{v}} J_{\vec{\sigma},\vec{\delta}}\left( \vec{v}, \vec{g}_{k+1} \right), 
\end{equation} 
where 
\begin{equation}
    \nabla_{\vec{v}} J_{\vec{\sigma},\vec{\delta}}\left( \vec{v}, \vec{g}_{k+1} \right) 
        = \mu \left( \vec{v} - \vec{u} \right) 
        + \beta \left( \mat{L_m} \right)^T \left( \mat{L_m} \vec{v} - \vec{g}_{k+1} \right) 
        - \left( \mat{L_m} \right)^T \vec{\sigma}_k - \vec{\delta}_k, 
\end{equation}
and the step size $\alpha > 0$ is chosen to provide a sufficient descent in direction of the gradient. 
Finally, the Lagrangian multipliers are updated after each iteration by 
\begin{equation}\label{eq:update-L}
    \begin{aligned}
        \vec{\sigma}_{k+1} & = 
            \vec{\sigma}_k - \beta \left( \mat{L_m} \vec{v}_{k+1} - \vec{g}_{k+1} \right), \\ 
        \vec{\delta}_{k+1} & = 
            \vec{\delta}_k - \mu \left( \vec{v}_{k+1} - \vec{u} \right).
    \end{aligned}
\end{equation} 
Algorithm \ref{algo:l1} is borrowed from \cite{sanders2017composite,scarnati2017using} and compactly describes the above ADMM. 
\begin{algorithm}
\caption{ADMM}\label{algo:l1}
\begin{algorithmic}[1]
    \State{Determine parameters $\mu, \beta$, and $tol$} 
    \State{Initialize $\vec{v}_0, \vec{g}_0, \vec{\sigma}_0$, and $\vec{\delta}_0$} 
    \For{$k=0$ to $K$}
        \While{$\norm{\vec{v}_{k+1}-\vec{v}_k} > tol$} 
            \State{Minimize $J$ for $\vec{g}$ according to \eqref{eq:update-g}} 
            \State{Minimize $J$ for $\vec{v}$ according to \eqref{eq:update-v}} 
        \EndWhile 
        \State{Update the Lagrangian multipliers according to \eqref{eq:update-L}} 
    \EndFor 
\end{algorithmic}
\end{algorithm} 
Figure \ref{fig:l1-reg} demonstrates the effect of $\ell^1$ regularization to a polynomial approximation with degree 
$p=13$ of a discontinuous sawtooth function $u(x) = \sign(x)-x$ on $[-1,1]$. 

Due to the Gibbs phenomenon, the original polynomial approximation $u_p$ in Figure \ref{fig:fun-values} shows spurious oscillations \cite{hewitt1979gibbs}. 
As a consequence, the approximation of the corresponding jump function by the PA operator for $m=3$ in Figure \ref{fig:jump-values} also shows undesired oscillations away from the  detected discontinuity at $x=0$. 
Note that the neighboring undershoots around $x=0$ are inherent in the PA operator for $m=3$. 
Yet, the oscillations left and right from these undershoots stem from the Gibbs phenomenon and thus are parasitical. 
Applying $\ell^1$ regularization to enhance sparsity of the PA transform, however, is able to correct these spurious 
oscillations, resulting in a sparse representation of the PA transformation in Figure \ref{fig:reg-jump-values}. 
The corresponding nodal values, illustrated in Figure \ref{fig:reg-fun-values}, now approximate the nodal values of the true solution accurately.

For the numerical test in Figure \ref{fig:l1-reg}, we have chosen the parameters $K=400$, $\mu = 0.005$ ($\lambda=4\cdot10^2$), $\beta=20$, $\alpha = 0.0001$, and $tol = 0.001$. 
The Lagrangian multipliers have been initialized with $\vec{\sigma}_0 = \vec{0}$ and $\vec{\delta}_0 = \vec{0}$.

\subsection{Preservation of mass conservation}
\label{sub:mass-corr} 

\begin{figure}[!htb]
  \centering
   \begin{subfigure}[b]{0.6\textwidth}
    \includegraphics[width=\textwidth]{%
      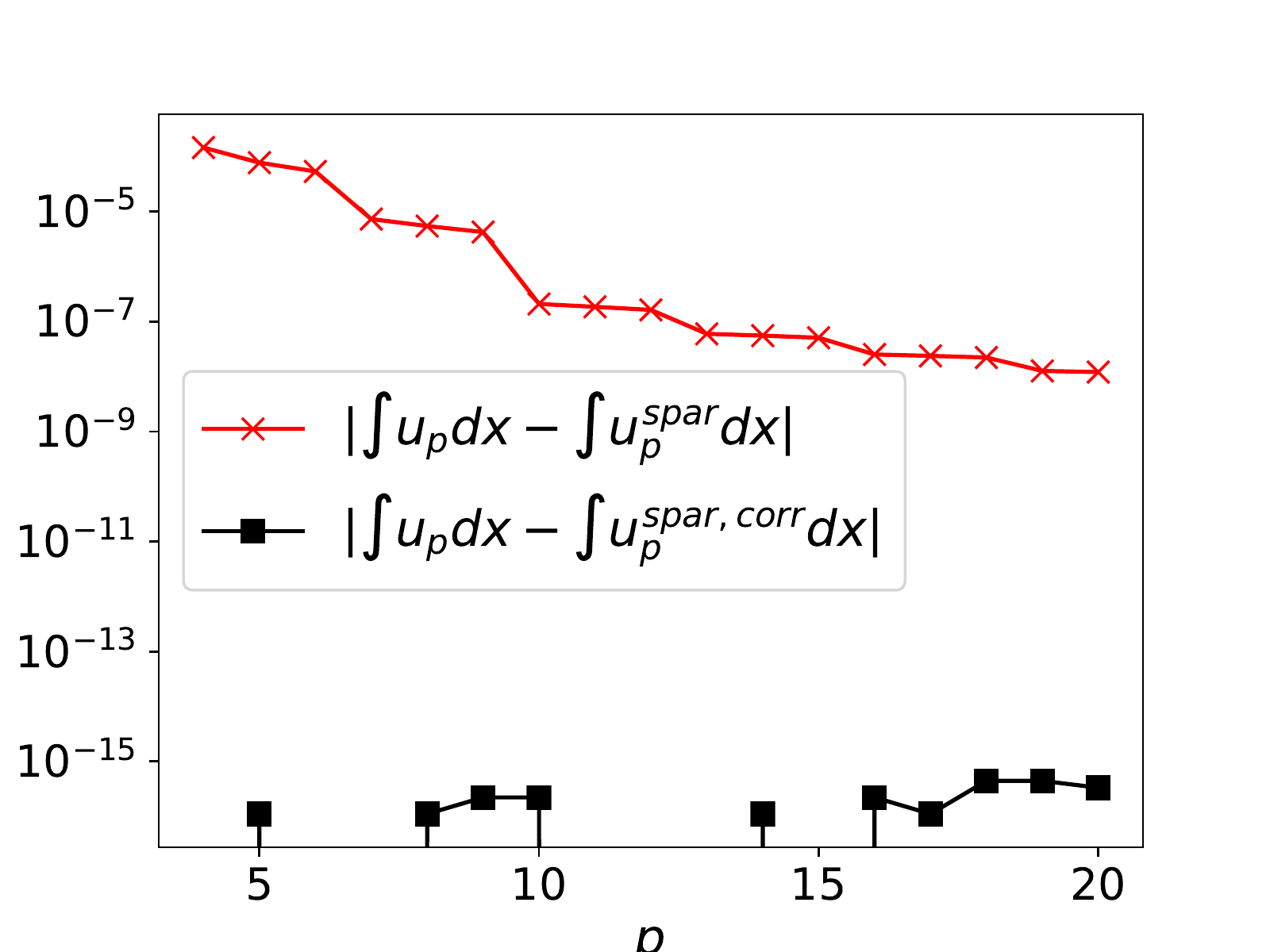}
  \end{subfigure}%
  \caption{Absolute difference between the mass of $u_p$ and $u_p^{\mathrm{spar}}$ (red crosses) 
and between the mass of $u_p$ and $u_p^{\mathrm{spar,corr}}$ (black squares) for $u(x) = 
\sign(x-0.5)+1$ on $\Omega_{\mathrm{ref}} = [-1,1]$ and increasing polynomial degree 
$p=4,\dots,20$. 
  In the cases where no (black) squares are visible the difference was below machine precision.}
  \label{fig:mass-corr}
\end{figure}

An essential property of standard (DG) finite element methods for hyperbolic conservation laws 
is that they are conservative, i.e., that 
\begin{equation}\label{eq:mass-cons}
  \int_{\Omega} u(t^{n+1},x) \d x 
    = \int_{\Omega} u(t^{n},x) \d x 
    + f\big|_{\partial \Omega}
\end{equation}
holds \cite{randall1992numerical,hesthaven2007nodal,gassner2013skew}. 
Any reasonable shock capturing procedure for hyperbolic conservation laws should preserve this 
property as well. 
In particular, in a troubled element $\Omega_i$ in which $\ell^1$ regularization is applied, 
\begin{equation}\label{eq:mass-cons-l1}
  \int_{\Omega_i} u_p^{\mathrm{spar}}(x) \d x 
    = \int_{\Omega_i} u_p(x) \d x 
\end{equation} 
should hold in order for $\ell^1$ regularization to preserve mass conservation of the underlying (DG) 
method. 
Unfortunately, \eqref{eq:mass-cons-l1} is violated when $\ell^1$ regularization is applied too 
naively. 
This is demonstrated in Figure \ref{fig:mass-corr} for a simple discontinuous example 
$u(x) = \sign(x-0.5)+1$ on $\Omega_{\mathrm{ref}} = [-1,1]$ with mass 
$\int_{\Omega_{\mathrm{ref}}} u \d x = 1$.

The (red) crosses illustrate the absolute difference between the mass of the polynomial 
interpolation $u_p$ and its sparse reconstruction $u_p^{\mathrm{spar}}$, i.e. 
\begin{equation}
  \left| \int_{\Omega_{\mathrm{ref}}} u_p(x) \d x 
    - \int_{\Omega_{\mathrm{ref}}} u_p^{\mathrm{spar}}(x) \d x \right|,
\end{equation}
for increasing polynomial degrees $p=4,\dots,20$. 
We observe that a naive application of $\ell^1$ regularization might destroy mass conservation. 
Thus, in the following, we present a simple fix for this problem.

\begin{remark}
A generic approach --- also to, for instance, ensure TVD and entropy 
conditions --- is to add additional constraints to the optimization problem \eqref{eq:op-problem}, 
e.g., 
\begin{equation}\label{eq:op-prob-mass-fix}
\begin{aligned}
  u_p^{\mathrm{spar}} 
    & = \argmin_{v \in \mathbb{P}_p\left(\Oref\right)} \left( \frac{1}{2}
\norm{v-u_p}_2^2 + \lambda  \norm{L_m v}_1 \right) \\
    \text{s.t.} & \quad 
    \int_{\Omega_{\mathrm{ref}}} v \d x 
    = \int_{\Omega_{\mathrm{ref}}} u_p(x) \d x 
\end{aligned}
\end{equation}
for conservation of mass to be preserved. 
In the present case,  where $u_p,v \in \mathbb{P}_p\left(\Oref\right)$ are 
expressed with respect to basis functions $\{\varphi_k\}_{k=0}^p$ with zero average for $k > 0$, the conditions in 
(\ref{eq:op-prob-mass-fix}) are easily met. 
Specifically, 
\begin{equation}
  u_p = \sum_{k=0}^p \hat{u}_k \varphi_k 
  \quad \text{and} \quad  
  v = \sum_{k=0}^p \hat{v}_k \varphi_k
\end{equation}
with 
\begin{equation}
  \int_{\Oref} \varphi_k \d x = \norm{\varphi_0}_2^2 \delta_{0k},
\end{equation}
implies that the additional constraint in \eqref{eq:op-prob-mass-fix} reduces to 
\begin{equation}
  \hat{u}_0 = \hat{v}_0. 
\end{equation}
Basis functions $\{\varphi_k\}_{k=0}^p$ with zero average for ${k > 0}$ are, for instance, given by the orthogonal 
basis (OGB) of Legendre polynomials.
\end{remark}

We now propose the following simple algorithm to repair mass conservation in the $\ell^1$ regularization:
\begin{algorithm}
\caption{Mass correction}\label{algo:mass-corr}
\begin{algorithmic}[1]
    \State{Compute $\hat{u}_0$} 
    \State{Compute $u_p^{\mathrm{spar}}$ according to 
\eqref{eq:op-problem}/\eqref{eq:op-problem-disc}} 
    \State{Represent $u_p^{\mathrm{spar}}$ w.r.t. an OGB: 
      $u_p^{\mathrm{spar}} = \hat{u}_0^{\mathrm{spar}} \varphi_0 + \dots + 
\hat{u}_p^{\mathrm{spar}} \varphi_p$} 
    \State{Replace $\hat{u}_0^{\mathrm{spar}}$ by $\hat{u}_0$}
\end{algorithmic}
\end{algorithm} 

\begin{figure}
  \centering
  \begin{subfigure}[b]{0.8\textwidth}
    \includegraphics[width=\textwidth]{%
      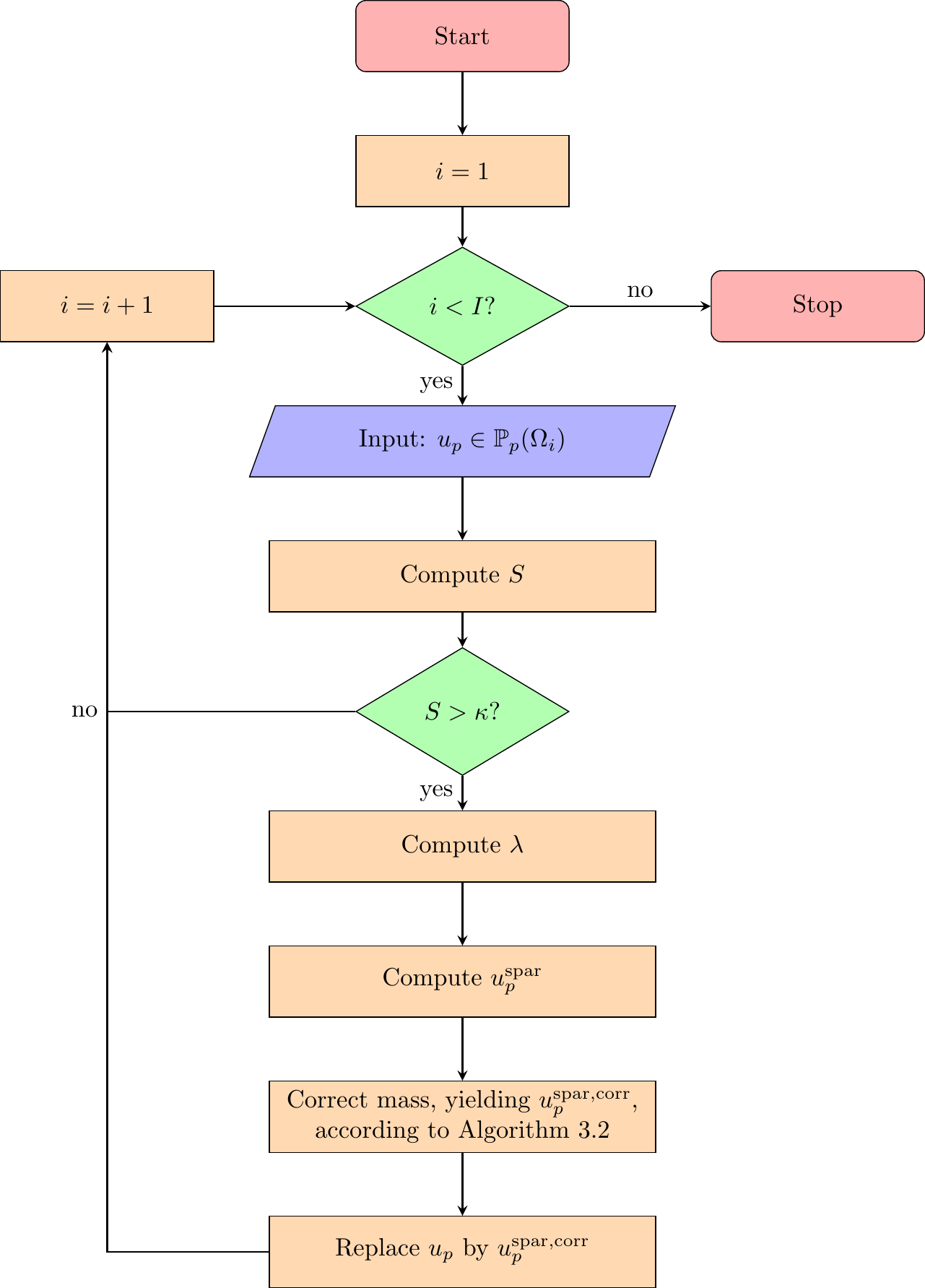}
  \end{subfigure}
  \caption{Flowchart describing $\ell^1$ regularization with mass correction. 
    The procedure is described for a fixed time and includes the loop over all $I$ elements.}
  \label{fig:flowchart}
\end{figure}

The advantage of this additional step is demonstrated in Figure \ref{fig:mass-corr} as well, where the 
absolute difference between the mass of $u_p$ and of its sparse reconstruction with additional mass 
correction, denoted by $u_p^{\mathrm{spar,corr}}$, is illustrated by (black) squares. 
In contrast to the sparse reconstruction without mass correction, illustrated by (red) crosses, 
$u_p^{\mathrm{spar,corr}}$ is demonstrated to preserve mass nearly up to machine precision 
($\approx 10^{-16}$). 
Finally, we note that for the test illustrated in Figure \ref{fig:reg-jump-values}, $\ell^1$ 
regularization with and without mass correction resulted in the same approximations, 
due to $u$ being an odd function. 
Thus, we omit those results.
We present a flowchart in Figure \ref{fig:flowchart} illustrating the proposed procedure for a fixed time $t < 
t_{\mathrm{end}}$.

\section{Numerical tests}
\label{sec:tests}

We now numerically demonstrate the application of $\ell^1$ regularization 
(with and without mass correction) to a nodal DG method 
for the inviscid Burgers' equation, the linear advection equation, and a nonlinear system of conservation laws. 
Our results show that $\ell^1$ regularization provides increased accuracy of the numerical solutions. 
In all numerical tests, we use a PA operator of third order and choose the same parameters as before, 
i.e., $\lambda_{\text{max}}=4\cdot10^2$,  $K=400$, $\beta=20$, $\alpha = 0.0001$, and $tol = 0.001$. 
We have made no effort to optimize these parameters.

\subsection{Inviscid Burgers' equation}
\label{sub:Burgers}

We start our numerical investigation by considering the inviscid Burgers' equation 
\begin{equation}\label{eq:Burgers}
    \partial_t u + \partial_x \left( \frac{u^2}{2} \right) = 0 
\end{equation} 
on $\Omega = [0,2]$ with initial condition 
\begin{equation}\label{eq:Burgers_IC}
    u(0,x) = u_0(x) = \sin\left( \pi x \right)
\end{equation} 
and periodic boundary conditions. 
For this test case, a shock discontinuity develops in the solution at $x=1$.

In the subsequent numerical tests, the usual local Lax--Friedrichs flux 
\begin{equation}\label{eq:LLF}
  \fnum(u_-,u_+) = \frac{1}{2} \left( f(u_+) + f(u_-) \right) 
    - \frac{{\alpha_{\mathrm{max}}}}{2} \left( u_+ - u_- \right)
\end{equation}
with maximum characteristic speed {$\alpha_{\mathrm{max}} = \max \{ |u_+|, |u_-| \}$}
is utilized. 
Further, we use the third order explicit strong stability preserving  
Runge--Kutta method with 
three stages (SSPRK(3,3)) for time integration (see \cite{gottlieb1998total}) and choose  
the ramp parameter ${\kappa = 0.8}$ in \eqref{eq:param-fun}.

\begin{figure}[!htb]
  \centering
  \begin{subfigure}[b]{0.32\textwidth}
    \includegraphics[width=\textwidth]{%
      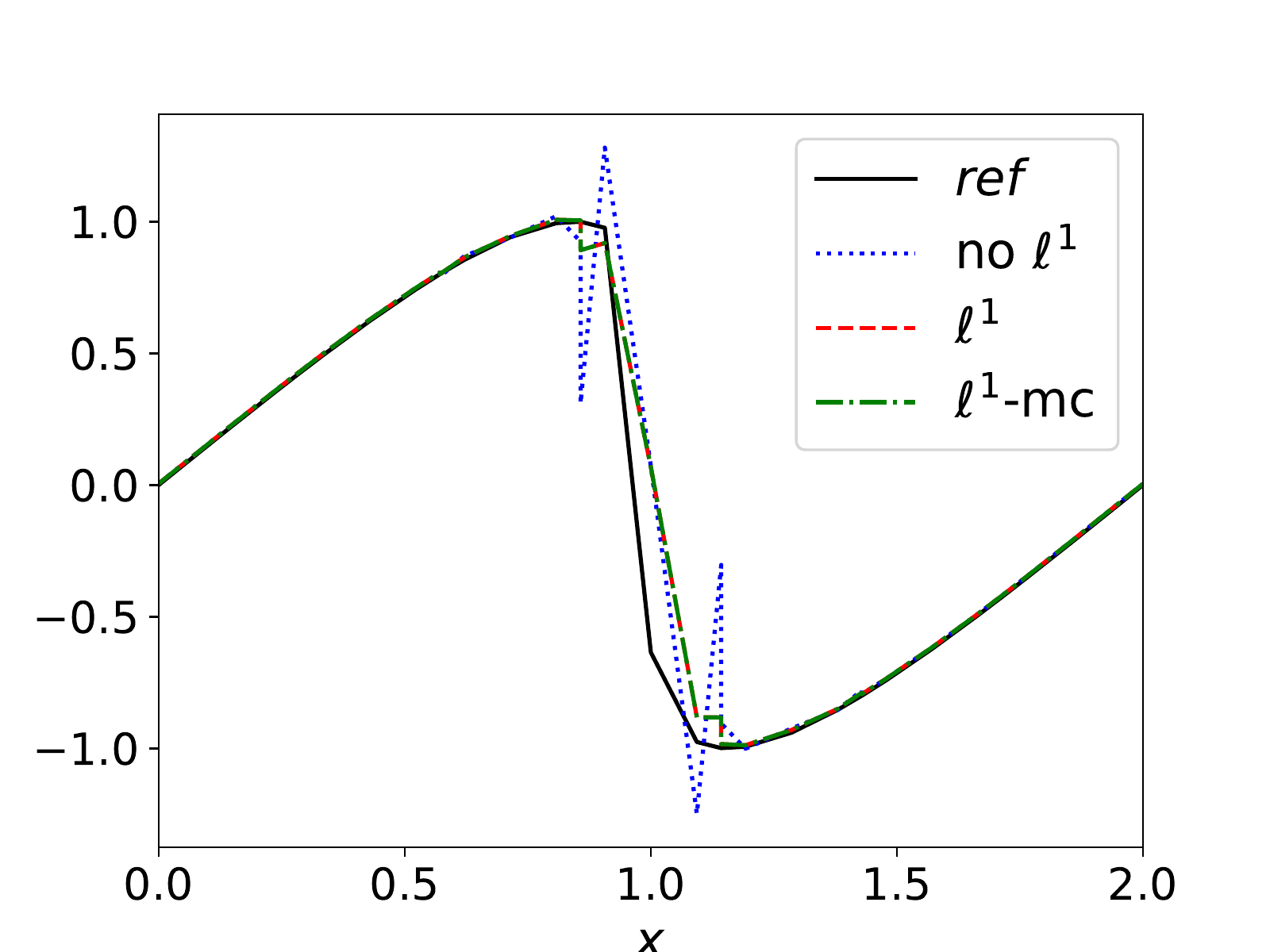}
    \caption{Solution for $p=4$.}
    \label{fig:Burgers_p04} 
  \end{subfigure}%
  ~ 
  \begin{subfigure}[b]{0.32\textwidth}
    \includegraphics[width=\textwidth]{%
      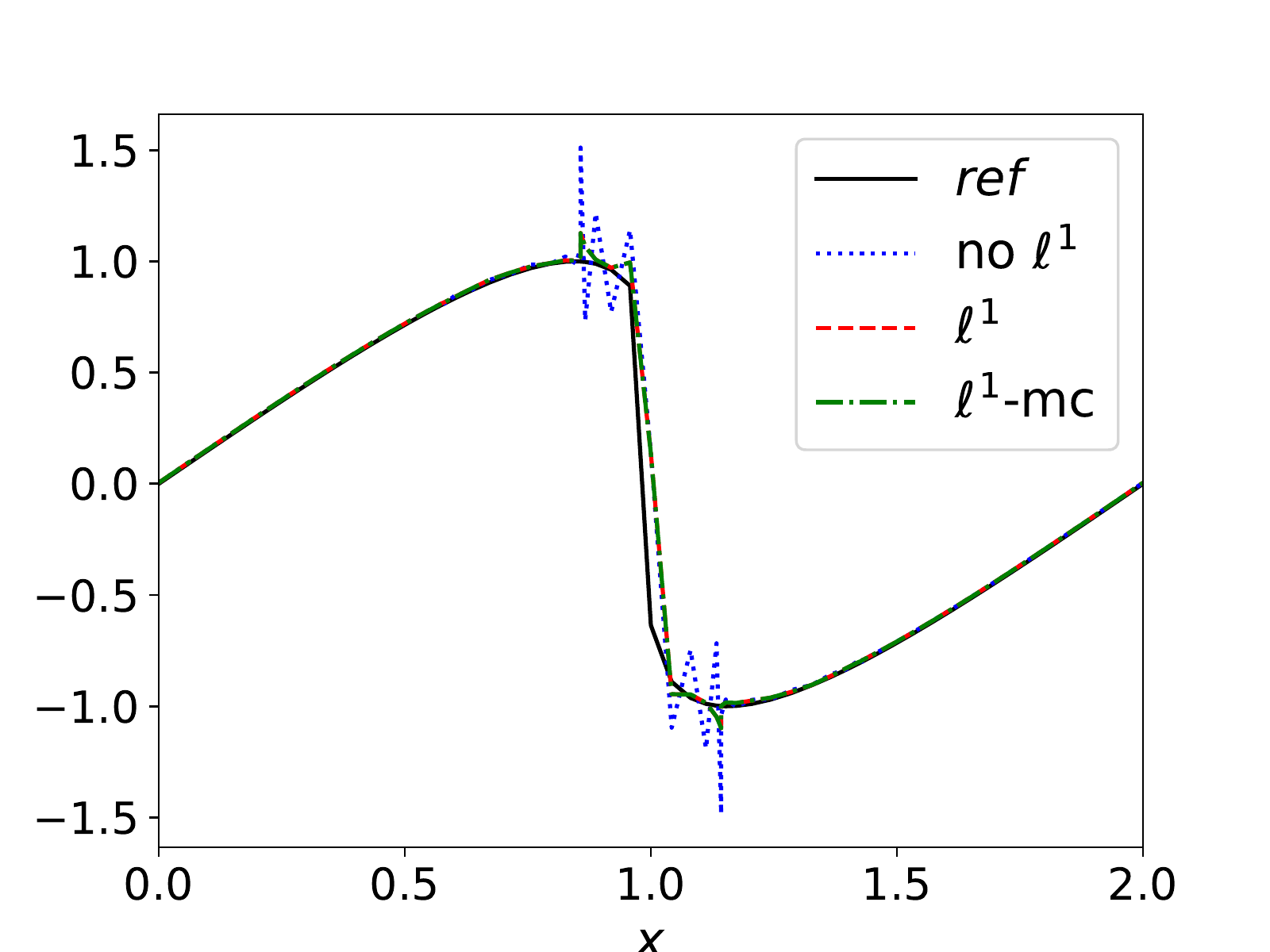}
    \caption{Solution for $p=10$.}
    \label{fig:Burgers_p10} 
  \end{subfigure}%
  ~ 
  \begin{subfigure}[b]{0.32\textwidth}
    \includegraphics[width=\textwidth]{%
      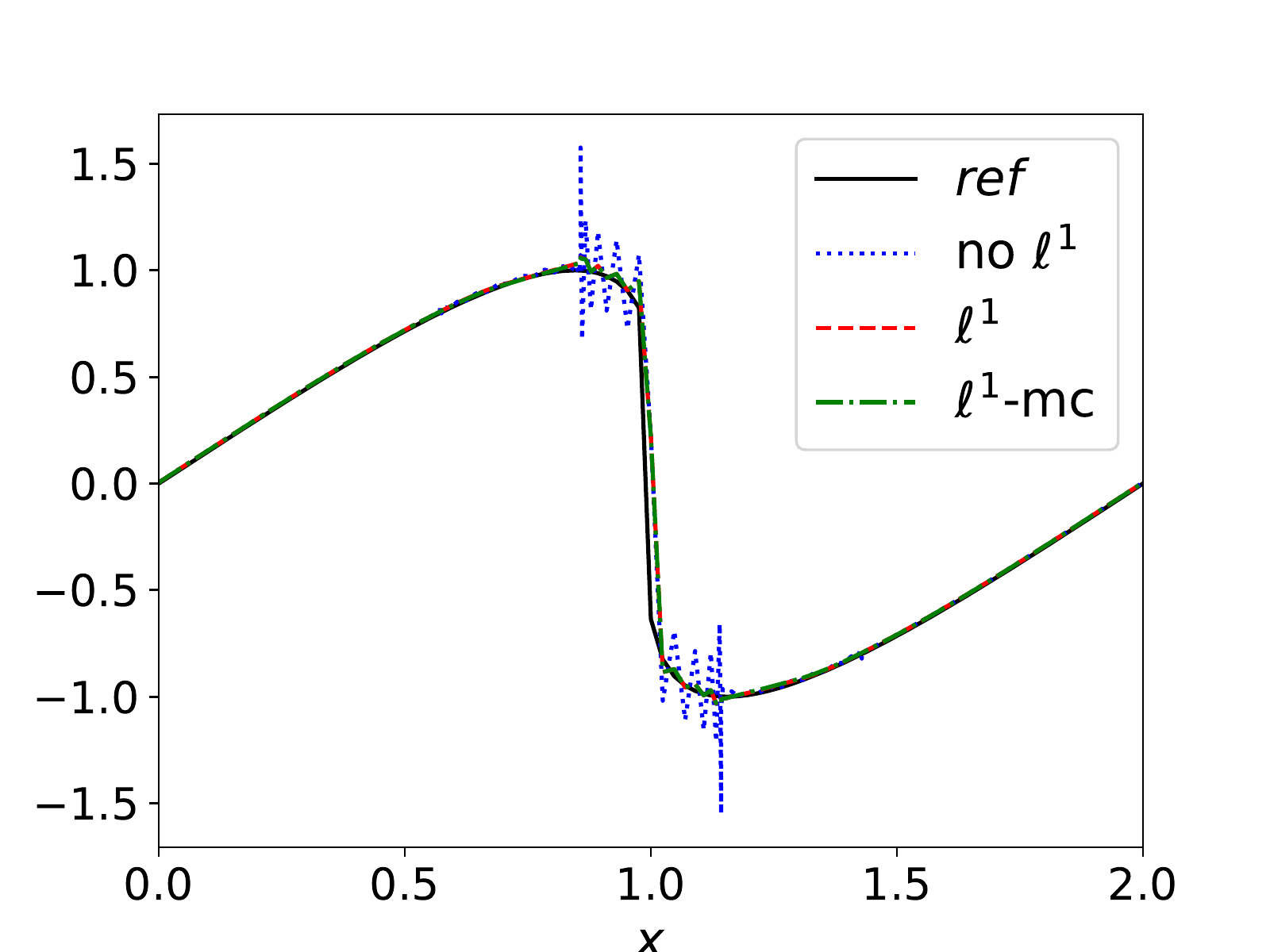}
    \caption{Solution for $p=18$.}
    \label{fig:Burgers_p18}
  \end{subfigure}%
  \\ 
  \begin{subfigure}[b]{0.32\textwidth}
    \includegraphics[width=\textwidth]{%
      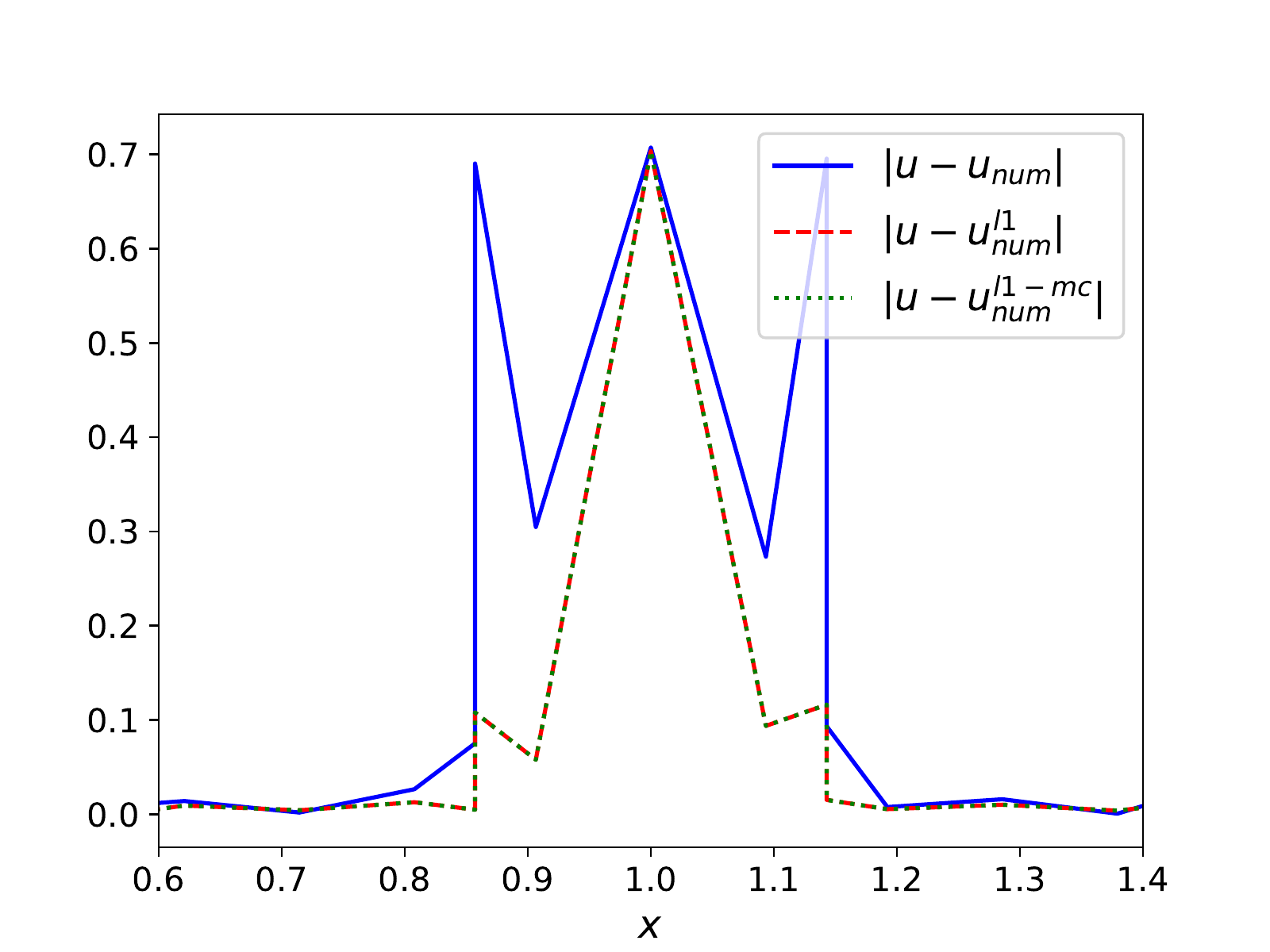}
    \caption{Errors for $p=4$.}
    \label{fig:pointwise_error_p04} 
  \end{subfigure}%
  ~
  \begin{subfigure}[b]{0.32\textwidth}
    \includegraphics[width=\textwidth]{%
      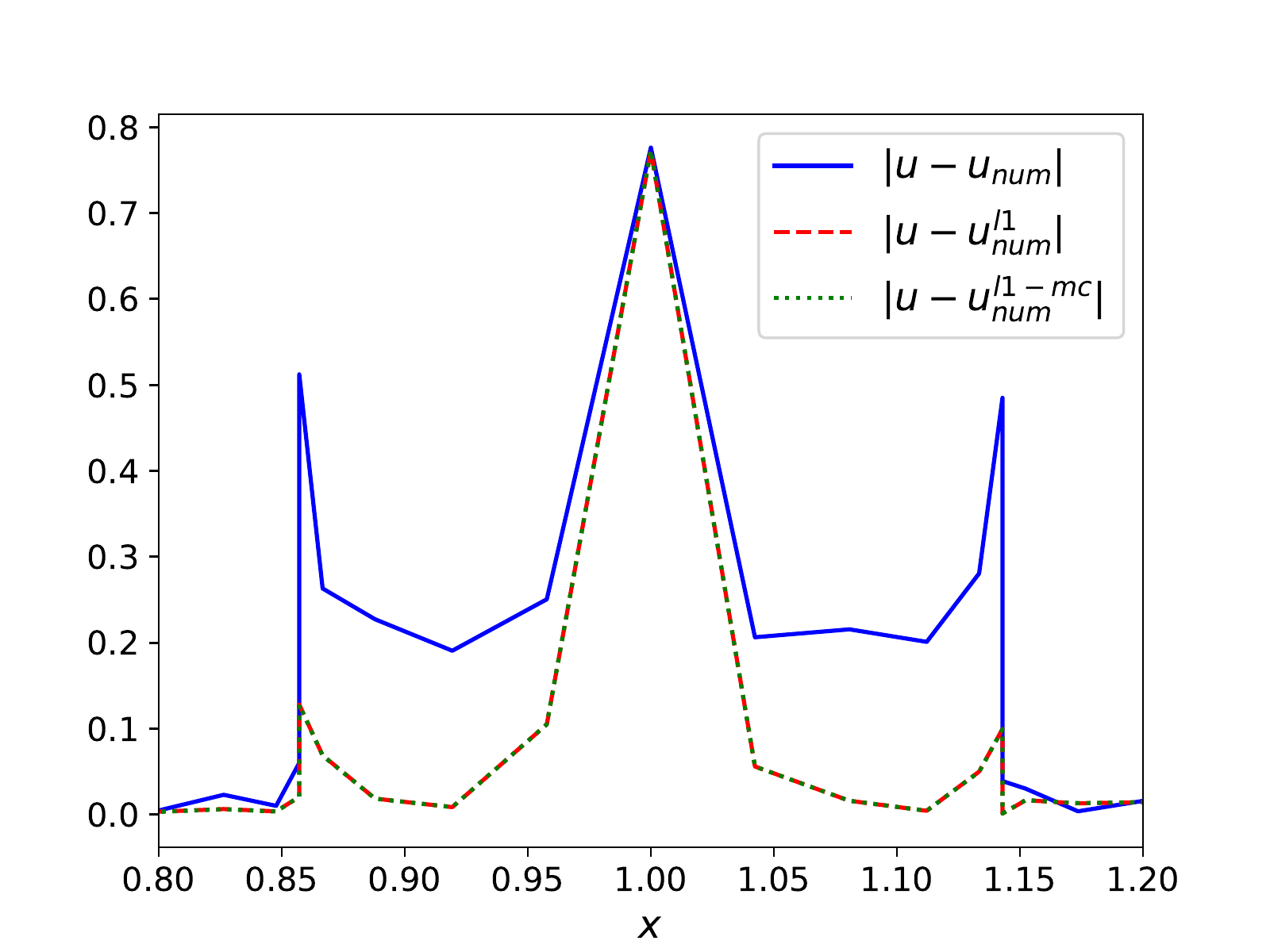}
    \caption{Errors for $p=10$.}
    \label{fig:pointwise_error_p10} 
  \end{subfigure}%
  ~
  \begin{subfigure}[b]{0.32\textwidth}
    \includegraphics[width=\textwidth]{%
      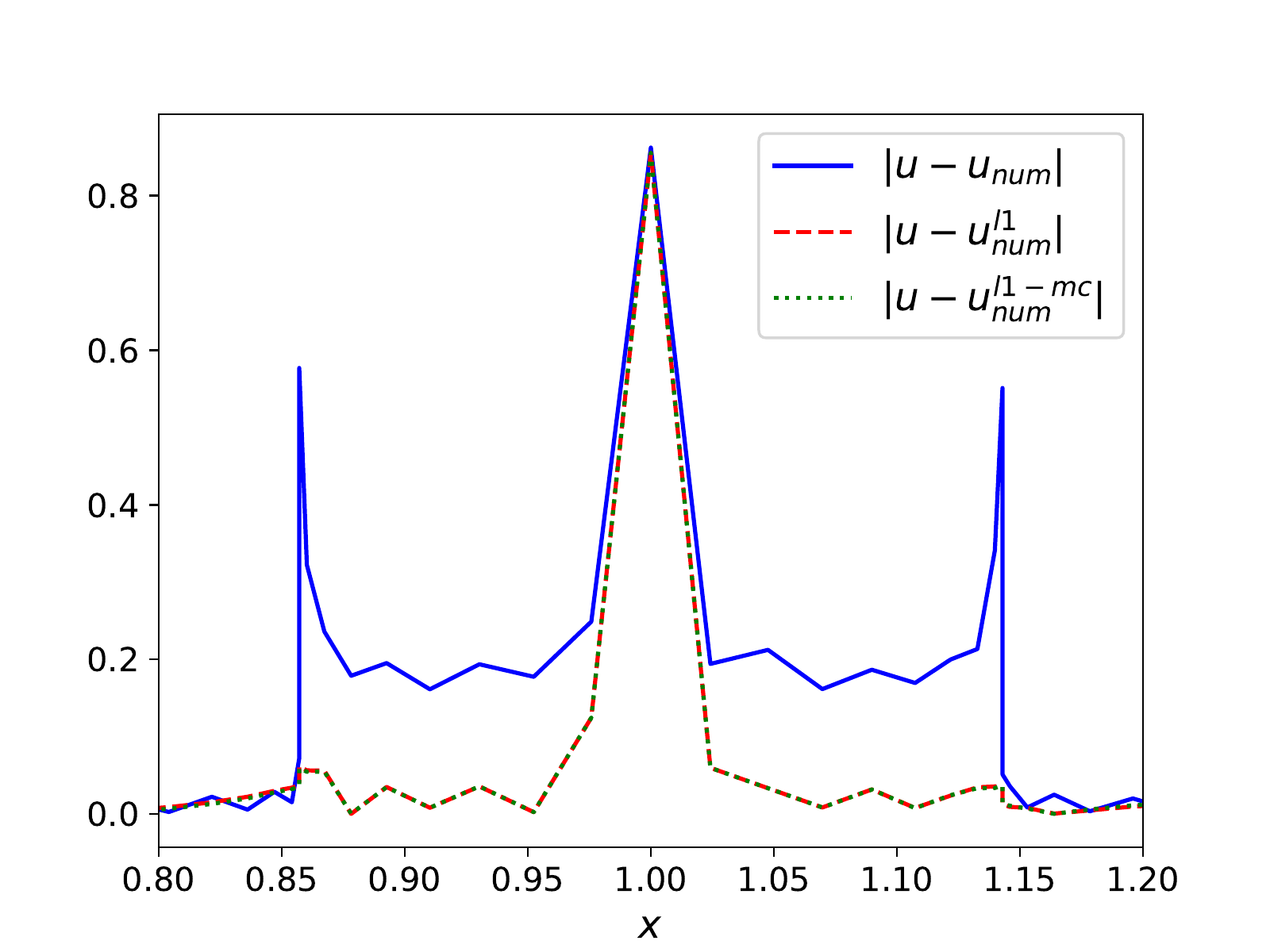} 
    \caption{Errors for $p=18$.}
    \label{fig:pointwise_error_p18}
  \end{subfigure}%
  \caption{Numerical solutions (and their pointwise errors) at $t=0.345$  by the DG method for 
$I=7$ elements 
and polynomial orders $p=4,10,18$. 
  Without $\ell^1$ regularization (straight blue line), with $\ell^1$ regularization (dashed red 
line), and with $\ell^1$ regularization and mass correction (dash-dotted green line).}
  \label{fig:l1-reg-Burgers}
\end{figure}

Figure \ref{fig:l1-reg-Burgers} illustrates the numerical solutions of the above test problem 
at time $t=0.345$. 
In Figure \ref{fig:Burgers_p04} a polynomial degree of $p=4$ was used on $I=7$ equidistant elements, 
while in Figure \ref{fig:Burgers_p10} a polynomial degree of $p=10$ was used on 
$I=7$ equidistant elements and in Figure \ref{fig:Burgers_p18} a polynomial degree of $p=18$ was used 
on $I=7$ equidistant elements.
All tests show spurious oscillations for the numerical solution without $\ell^1$ regularization. 
These oscillations are significantly reduced when  $\ell^1$ regularization is applied to the numerical 
solution after every time step.
As a consequence, Figures \ref{fig:pointwise_error_p04}, \ref{fig:pointwise_error_p10}, 
\ref{fig:pointwise_error_p18} illustrate how the pointwise error of the numerical solutions 
utilizing $\ell^1$ regularization (with and without mass correction) is reduced compared to 
the numerical solution without $\ell^1$ regularization. 
Again, considering an odd function, only very slight differences are observed between 
$\ell^1$ regularization with and without mass correction. 
The reference solution $u$ was computed using characteristic tracing.

%
%
We now extend the above error analysis for the $\ell^1$ regularization. 
Table \ref{tab:errors} lists the different common types of errors of the numerical solution for the above test problem \eqref{eq:Burgers}, \eqref{eq:Burgers_IC} and a varying polynomial degree $p$ as well as a varying number of equidistant elements $I$. 
Here we consider the global errors with respect to the discrete $\mat{M}$-norm, as defined in 
(\ref{eq:disc_inner_prod}), 
approximating the continuous $L^2$-norm and the discrete $1$-norm, given by
$\norm{u_p}_1 = \sum_{k=0}^p \omega_k \left|u_p(x_k)\right| = \scp{\vec{1}}{\left| \vec{u} 
\right|}_{\mat{M}}$ on the reference element.
The discrete $\infty$-norm is given by  $\norm{u_p}_\infty = \max_{k=0,\dots,p} |u_p(x_k)|$. 
For the discrete $\mat{M}$-norm and discrete $1$-norm, the global norms are defined by summing up over the 
weighted local norms, i.e.,  
\begin{equation}
    \norm{u_{\mathrm{num}}}^2 = \sum_{i=1}^I \ \frac{\left| \Omega_i \right|}{2} \ \norm{ u_p^{(i)} }^2,
\end{equation} 
where $u_p^{(i)}$ denotes the numerical solution (polynomial approximation) on the $i$th element 
$\Omega_i$.

\begin{table}[!htb] 
\resizebox{\columnwidth}{!}{%
    \centering 
    \begin{tabular}{ r | r || c | c | c || c | c | c || c | c | c }
      \hline
            & & \multicolumn{3}{c||}{$\norm{\cdot}_M$-error} 
	      & \multicolumn{3}{c||}{$\norm{\cdot}_1$-error} 
	      & \multicolumn{3}{c}{$\norm{\cdot}_\infty$-error} \\ 
        $p$ & $I$ & no $\ell^1$ & $\ell^1$ & $\ell^1$-mc 
		  & no $\ell^1$ & $\ell^1$ & $\ell^1$-mc 
		  & no $\ell^1$ & $\ell^1$ & $\ell^1$-mc \\ 
        \hline \hline
       
	  3 &  15 & 1.5e-2 & 3.3e-2 		& 3.3e-2 
		  & 1.2e-2 & 2.2e-1		& 1.2e-2 
		  & 1.7e-1 & 3.3e-1		& 3.3e-1 \\
            &  31 & 1.5e-2 & 2.0e-2		& 2.0e-2 
		  & 1.1e-2 & 1.2e-2		& 1.2e-2 
		  & 2.2e-1 & 3.2e-1		& 3.2e-1 \\
            &  63 & 1.9e-2 & 2.4e-2		& 2.4e-2 
		  & 1.1e-2 & 1.1e-2 \checkmark 	& 1.1e-2 
		  & 3.6e-1 & 5.6e-1		& 5.6e-1 \\
            & 127 & 3.3e-2 & 2.7e-2 \checkmark 	& 2.7e-2 
		  & 1.3e-2 & 1.1e-2 \checkmark 	& 1.1e-2 
		  & 1.4e-0 & 8.3e-1 \checkmark 	& 8.3e-1 \\
        \hline
       
          4 &  15 & 7.0e-2 & 5.9e-2 \checkmark & 5.9e-2 
		  & 3.9e-2 & 2.7e-2 \checkmark & 2.7e-2
		  & 7.6e-1 & 7.5e-1 \checkmark & 7.5e-1 \\
            &  31 & 5.5e-2 & 4.6e-2 \checkmark & 4.6e-2 
		  & 2.4e-2 & 1.7e-2 \checkmark & 1.7e-2 
		  & 8.7e-1 & 8.6e-1 \checkmark & 8.6e-1 \\
            &  63 & 4.3e-2 & 3.9e-2 \checkmark & 3.9e-2 
		  & 1.6e-2 & 1.3e-2 \checkmark & 1.3e-2 
		  & 1.0    & 1.0e-0 \checkmark & 1.0e-0 \\
            & 127 & 3.8e-2 & 3.6e-2 \checkmark & 3.6e-2 
		  & 1.2e-2 & 1.2e-2 \checkmark & 1.1e-2 
		  & 1.4    & 1.3e-0 \checkmark & 1.3e-0 \\
        \hline
       
          5 &  15 & 1.6e-2 & 1.5e-2 \checkmark & 1.5e-2 
		  & 1.2e-2 & 1.2e-2 \checkmark & 1.2e-2 
		  & 2.5e-1 & 1.9e-1 \checkmark & 1.7e-1 \\
            &  31 & 2.0e-2 & 1.4e-2 \checkmark & 1.3e-2 
		  & 1.3e-2 & 1.0e-2 \checkmark & 1.0e-2 
		  & 5.3e-1 & 2.6e-1 \checkmark & 2.5e-1 \\
            &  63 & 7.3e-2 & 1.9e-2 \checkmark & 1.6e-2 
		  & 2.2e-2 & 1.0e-2 \checkmark & 1.0e-2 
		  & 3.6e-0 & 5.3e-1 \checkmark & 4.3e-1 \\
            & 127 & NaN    & 2.7e-2 !          & 2.2e-2      
		  & NaN    & 1.0e-2 !          & 1.0e-2    
		  & NaN    & 1.1e-0 ! 	       & 9.1e-1 \\
        \hline
       
          6 &  15 & 5.9e-2 & 5.2e-2 \checkmark & 5.2e-2
		  & 3.2e-2 & 2.3e-2            & 2.3e-2      
		  & 8.1e-1 & 8.0e-1 \checkmark & 8.0e-1 \\
            &  31 & 4.7e-2 & 4.3e-2 \checkmark & 4.3e-2 
		  & 2.0e-2 & 1.6e-2 \checkmark & 1.6e-2 
		  & 9.6e-1 & 9.5e-1 \checkmark & 9.5e-1 \\
            &  63 & 3.9e-2 & 3.7e-2 \checkmark & 3.6e-2 
		  & 1.2e-2 & 1.2e-2 \checkmark & 1.2e-2 
		  & 1.2e-0 & 1.1e-0 \checkmark & 1.1e-0 \\
            & 127 & NaN    & 3.4e-2 ! 	       & 3.2e-2         
		  & NaN    & 1.1e-2 ! 	       & 1.1e-2         
		  & NaN    & 1.3e-0 ! 	       & 1.3e-0 \\
        \hline
       
          7 &  15 & 1.7e-2 & 1.8e-2    	       & 1.8e-2        
		  & 1.4e-2 & 1.2e-2 \checkmark & 1.2e-2 
		  & 2.2e-1 & 2.7e-1 	       & 2.8e-1 \\
            &  31 & 2.5e-2 & 1.7e-2 \checkmark & 1.8e-2 
		  & 1.5e-2 & 1.1e-2 \checkmark & 1.1e-2  
		  & 7.0e-1 & 3.8e-1 \checkmark & 4.0e-1 \\
            &  63 & 5.0e-1 & 2.1e-2 \checkmark & 2.2e-2  
		  & 1.4e-1 & 1.0e-2 \checkmark & 1.0e-2  
		  & 1.1e+1 & 6.9e-1 \checkmark & 7.5e-1 \\
            & 127 & 2.8e-2 & 2.5e-2 \checkmark & 2.5e-2  
		  & 1.0e-2 & 1.0e-2 \checkmark & 1.0e-2  
		  & 1.3e-0 & 1.2e-0 \checkmark & 1.2e-0 \\
        \hline
       
          8 &  15 & 5.9e-2 & 4.9e-2 \checkmark & 4.9e-2  
		  & 3.2e-2 & 2.1e-2 \checkmark & 2.1e-2  
		  & 9.9e-1 & 8.6e-1 \checkmark & 8.6e-1 \\
            &  31 & 4.5e-2 & 4.0e-2 \checkmark & 4.0e-2  
		  & 1.9e-2 & 1.4e-2 \checkmark & 1.4e-2  
		  & 1.0e-0 & 1.0e-0 \checkmark & 1.0e-0 \\
            &  63 & 3.8e-2 & 3.6e-2 \checkmark & 3.7e-2  
		  & 1.2e-2 & 1.2e-2 \checkmark & 1.2e-2  
		  & 1.4e-0 & 1.3e-0 \checkmark & 1.3e-0 \\
            & 127 & 1.5e-1 & 3.2e-2 \checkmark & 3.1e-2  
		  & 3.1e-2 & 1.1e-2 \checkmark & 1.1e-2  
		  & 6.9e-0 & 1.4e-0 \checkmark & 1.4e-0 \\
        \hline 
       
          9 &  15 & 2.0e-2 & 1.9e-2 \checkmark & 1.8e-2  
		  & 1.5e-2 & 1.3e-2 \checkmark & 1.3e-2  
		  & 4.5e-1 & 3.0e-1 \checkmark & 2.9e-1 \\
            &  31 & 7.8e-2 & 1.9e-2 \checkmark & 1.8e-2  
		  & 2.5e-2 & 1.1e-2 \checkmark & 1.1e-2  
		  & 4.5e-0 & 4.5e-1 \checkmark & 4.1e-1 \\
            &  63 & NaN    & 2.3e-2 ! 	       & 2.2e-2          
		  & NaN    & 1.0e-2 !  	       & 1.0e-2         
		  & NaN    & 9.0e-1 !	       & 8.3e-1 \\
            & 127 & NaN    & 2.7e-2 !	       & 2.8e-2           
		  & NaN    & 1.0e-2 !	       & 1.0e-2           
		  & NaN    & 1.4e-0 !	       & 1.4e-0 \\
        \hline
    
    \end{tabular} 
}
    \caption{
        Errors of the numerical solutions without and with $\ell^1$ regularization. 
        In cases where the error value is NaN, the numerical solution broke down before the final 
time was reached.
	``no $\ell^1$'' refers to the underlying DG method without $\ell^1$ regularization, 
        ``$\ell^1$'' refers to the DG method with $\ell^1$ regularization, and 
        ``$\ell^1$-mc'' refers to the DG method with $\ell^1$ regularization and additional mass 
correction.} 
    \label{tab:errors}
\end{table}

Table \ref{tab:errors} demonstrates that for almost all these norms as well as combinations of 
polynomial degrees $p=4,5,6,7,8,9$ and number of elements $I=15,32,63,127$, 
the numerical solution with $\ell^1$ regularization is more accurate than the numerical 
solution without $\ell^1$ regularization. 
Further, we observe just a slight difference in accuracy for $\ell^1$ regularization 
with and without mass correction. 
See, for instance, $p=5$ in Table \ref{tab:errors}.
We only utilize odd number of elements, so that the shock discontinuity arises in the interior of an element. 
By using an even number of elements, on the other hand, the shock would arise at the interface between two elements 
and the error analysis of the $\ell^1$ regularization would be blurred by the dissipation added by the numerical flux. 

In Table \ref{tab:errors}, all cases where the accuracy is increased or remains the same by applying 
$\ell^1$ regularization are flagged with a checkmark. 
It should be stressed that no effort was made to optimize the parameters in the $\ell^1$ regularization. 
In particular, the parameters in the $\ell^1$ regularization have not been adapted to the specific choice of the 
polynomial degree $p$, the number of elements $I$, or the number of time steps used in the solver.
We think that by further investigations of the parameters in the $\ell^1$ regularization, 
all above errors could be further reduced. 

Finally, special attention should be given to the cases flagged by an exclamation mark (!). 
In these cases, the numerical solver without $\ell^1$ regularization broke down completely. 
Yet, when $\ell^1$ regularization was applied, the same computations yielded fairly accurate 
numerical solutions.

\subsection{Linear advection equation}
\label{sub:linear}

The prior test case featured a shock discontinuity and might not fully reflect the behavior of 
the $\ell^1$ regularized DG method in smooth regions. 
Note that for the underlying spectral DG method used in this work, given $p+1$ nodes, 
the optimal order of convergence is $p+1$ in sufficiently smooth regions.
Especially in smooth regions, it is desirable that the convergence properties of the underlying 
method are preserved when using the $\ell^1$ regularization method. 
Convergence of the $\ell^1$ regularized DG method, however, also depends on the convergence of the PA 
operator. 
In all numerical tests presented in this work, a PA operator of order three, i.e., $L_3$, is used. 
Hence, if falsely activated (false positive; see Remark \ref{rem3}) by the discontinuity sensor 
\eqref{eq:sensor}, $\ell^1$ regularization might affect the order of convergence of the underlying 
method in smooth regions. 
In our numerical tests, we observed the sensor to be fairly reliable for sufficient resolution, 
and  $\ell^1$ regularization did not get activated in smooth regions.

To investigate this potential drawback and demonstrate the reliability of the $\ell^1$ regularized DG 
method in smooth regions, we now consider the linear advection equation 
\begin{equation}
  \partial_t u + \partial_x u = 0
\end{equation} 
on $\Omega=[0,2]$ with initial condition 
\begin{equation}
  u(0,x) = \sin\left( 2 \pi x \right)
\end{equation} 
and periodic boundary conditions, which provides a smooth solution for all times.
Table \ref{tab:errors-lin} lists comparative errors for DG with and without the $\ell^1$ 
regularization and with the additional mass conservation correction term at time $T = 2$.

\begin{table}[!htb] 
\resizebox{\columnwidth}{!}{%
    \centering 
    \begin{tabular}{ r | r || c | c | c || c | c | c || c | c | c } 
      \hline
            & & \multicolumn{3}{c||}{$\norm{\cdot}_M$-error} 
	      & \multicolumn{3}{c||}{$\norm{\cdot}_1$-error} 
	      & \multicolumn{3}{c}{$\norm{\cdot}_\infty$-error} \\ 
        $p$ & $I$ & no $\ell^1$ & $\ell^1$ & $\ell^1$-mc 
		  & no $\ell^1$ & $\ell^1$ & $\ell^1$-mc 
		  & no $\ell^1$ & $\ell^1$ & $\ell^1$-mc \\ 
        \hline \hline
       
	  3 &   2 & 1.2e-0 & 1.2e-0 & 1.2e-0
		  & 1.5e-0 & 1.6e-0 & 1.6e-0
		  & 9.5e-1 & 9.8e-1 & 9.8e-1 \\
	    &   4 & 1.3e-1 & 4.4e-1 & 4.4e-1
		  & 1.4e-1 & 5.8e-1 & 5.8e-1
		  & 1.3e-1 & 4.1e-1 & 4.1e-1 \\
	    &   8 & 6.3e-3 & 6.3e-3 & 6.3e-3
		  & 7.0e-3 & 7.0e-3 & 7.0e-3
		  & 1.2e-2 & 1.2e-2 & 1.2e-2 \\
	    &  16 & 3.8e-4 & 3.8e-4 & 3.8e-4
		  & 3.8e-4 & 3.8e-4 & 3.8e-4
		  & 9.9e-4 & 9.9e-4 & 9.9e-4 \\
        \hline
       
          4 &   2 & 3.4e-1 & 9.2e-1 & 9.2e-1 
		  & 4.2e-1 & 9.7e-1 & 9.6e-1
		  & 3.8e-1 & 9.0e-1 & 8.8e-1 \\
	    &   4 & 7.8e-3 & 7.8e-3 & 7.8e-3 
		  & 1.0e-2 & 1.0e-2 & 1.0e-2
		  & 1.2e-2 & 1.2e-2 & 1.2e-2 \\
	    &   8 & 4.2e-4 & 4.2e-4 & 4.2e-4
		  & 4.4e-4 & 4.4e-4 & 4.4e-4 
		  & 1.2e-3 & 1.2e-3 & 1.2e-3 \\
	    &  16 & 1.3e-5 & 1.3e-5 & 1.3e-5
		  & 1.3e-5 & 1.3e-5 & 1.3e-5
		  & 4.4e-5 & 4.4e-5 & 4.4e-5 \\
        \hline
       
          5 &   2 & 8.0e-2 & 1.0e-0 & 1.0e-0 
		  & 1.0e-1 & 1.3e-0 & 1.3e-0
		  & 1.3e-1 & 8.4e-1 & 7.8e-1 \\
	    &   4 & 2.1e-3 & 2.1e-3 & 2.1e-3
		  & 2.3e-3 & 2.3e-3 & 2.3e-3
		  & 5.3e-3 & 5.3e-3 & 5.3e-3 \\
	    &   8 & 2.8e-5 & 2.8e-5 & 2.8e-5 
		  & 2.9e-5 & 2.9e-5 & 2.9e-5
		  & 7.7e-5 & 7.7e-5 & 7.7e-5 \\
	    &  16 & 1.2e-6 & 1.2e-6 & 1.2e-6 
		  & 1.5e-6 & 1.5e-6 & 1.5e-6
		  & 1.6e-6 & 1.6e-6 & 1.6e-6 \\
        \hline
        
          6 &   2 & 2.1e-2 & 9.9e-1 & 9.9e-1 
		  & 2.5e-2 & 1.1e-0 & 1.1e-0
		  & 5.1e-2 & 9.9e-1 & 9.9e-1 \\
	    &   4 & 7.5e-5 & 7.5e-5 & 7.5e-5
		  & 8.0e-5 & 8.0e-5 & 8.0e-5
		  & 1.8e-4 & 1.8e-4 & 1.8e-4 \\
	    &   8 & 6.0e-6 & 6.0e-6 & 6.0e-6 
		  & 7.6e-6 & 7.6e-6 & 7.6e-6
		  & 7.5e-6 & 7.5e-6 & 7.5e-6 \\
	    &  16 & 7.3e-7 & 7.3e-7 & 7.3e-7 
		  & 9.4e-7 & 9.4e-7 & 9.4e-7
		  & 7.4e-7 & 7.4e-7 & 7.4e-7 \\
        \hline
        
          7 &   2 & 2.0e-3 & 2.0e-3 & 2.0e-3 
		  & 2.0e-3 & 2.0e-3 & 2.0e-3
		  & 6.4e-3 & 6.4e-3 & 6.4e-3 \\
	    &   4 & 3.9e-5 & 3.9e-5 & 3.9e-5 
		  & 4.9e-5 & 4.9e-5 & 4.9e-5
		  & 6.8e-5 & 6.8e-5 & 6.8e-5 \\
	    &   8 & 3.9e-6 & 3.9e-6 & 3.9e-6 
		  & 5.0e-6 & 5.0e-6 & 5.0e-6
		  & 4.1e-6 & 4.1e-6 & 4.1e-6 \\
	    &  16 & 4.9e-7 & 4.9e-7 & 4.9e-7 
		  & 6.3e-7 & 6.3e-7 & 6.3e-7
		  & 4.9e-7 & 4.9e-7 & 4.9e-7 \\
        \hline
    
    \end{tabular} 
}
    \caption{ 
	$u(0,x) = \sin\left( 2 \pi x \right)$.
        Errors of the numerical solutions without and with $\ell^1$ regularization. 
        ``no $\ell^1$'' refers to the underlying DG method without $\ell^1$ regularization, 
        ``$\ell^1$'' refers to the DG method with $\ell^1$ regularization, and 
        ``$\ell^1$-mc'' refers to the DG method with $\ell^1$ regularization and additional mass 
correction. 
        } 
    \label{tab:errors-lin}
\end{table}

Our results indicate that $\ell^1$ regularization is only activated, and thus affects convergence of 
the underlying method, if the solution is heavily underresolved. 
This can be noted in Table \ref{tab:errors-lin} from $p=3$ and $N=2,4$ as well as $p=4,5,6$ and 
$N=2$. 
For $p=7$, even $N=2$ provides sufficient resolution and the $\ell^1$ regularization is not activated.
For all numerical solutions, using at least $N=8$ elements, $\ell^1$ regularization does not affect 
accuracy and convergence of the underlying method in smooth regions. 
Finally, we note that in the cases where the numerical solution is heavily underresolved and $\ell^1$ 
regularization is activated, $\ell^1$ regularization with additional mass correction (l1-mc) provides 
slightly more accurate solutions than $\ell^1$ regularization without mass correction (l1). 
Finally, we note that because the $\ell^1$ regularization is only activated in elements 
containing discontinuities, the efficiency of our new method is comparable to the underlying DG 
method. 
Moreover, as in \cite{scarnati2017using}, we observed that for our specific set of test 
problems we were able to maintain stability for  time step sizes larger than the standard CFL 
constraints suggest. 
Theoretical justification for this will be part of future investigations.

\subsection{Systems of conservation laws}
\label{sub:systems}

We now extend our hybrid $\ell^1$ regularized DG method to the nonlinear system of conservation 
laws
\begin{equation}\label{eq:system}
  \partial_t 
  \begin{pmatrix} u_0 \\ u_1 \end{pmatrix} 
  + \frac{1}{2} \partial_x 
  \begin{pmatrix} 
    u_0^2 + u_1^2 \\ 
    2 u_0 u_1
    \end{pmatrix}
    = 0
\end{equation}
in the domain $\Omega = [0,2]$. 
System \eqref{eq:system} originates from a truncated polynomial chaos approach for Burgers' equation 
with uncertain initial condition 
\cite{pettersson2009numerical,pettersson2015polynomial,offner2018stability}. 
In this context, $u_0$ models the expected value of the numerical solution while $u_1^2$ 
approximates the variance. 
%
For the spatial semidiscretization of \eqref{eq:system} we follow 
\cite{offner2018stability,ranocha2018stability}, where a skew-symmetric formulation 
\begin{equation}\label{eq:PC_semidiscsystem}
\begin{aligned}
  \vec{\dot{u_k}}
  =&
  - \frac{1}{3} \sum_{i,j=0}^{1}  \langle \phi_i\phi_j\phi_k \rangle
    \left( \mat{D} \vec{u_i u_j} + \mat{u_j}^* \mat{D} \vec{u_i} \right)
  \\&
  - \mat{M}^{-1} \mat{R}^T \mat{B}\left(
    \vec{f_k}^\mathrm{num}
    - \sum_{i,j=0}^{1}  \langle \phi_i\phi_j\phi_k \rangle \left(
      \frac{1}{3} \mat{R} \vec{u_i u_j}
      + \frac{1}{6}\left(\mat{R} \vec{u_i}\right) \circ \left(\mat{R} 
\vec{u_j}\right)
      \right)
    \right)
\end{aligned}
\end{equation}
was proposed. 
Here $\langle \phi_i\phi_j\phi_k \rangle$ denotes the triple product
$\int  \phi_i(\xi)\phi_j(\xi)\phi_k(\xi)\omega(\xi) \dif \xi$, 
$\circ$ denotes the componentwise (Hadamard) product of two vectors, and $\{\phi_k\}$ is a set of 
orthogonal polynomials (typically Hermite polynomials are used). 

It was further proved in \cite{offner2018stability} that \eqref{eq:PC_semidiscsystem} yields an 
entropy  conservative semidiscretization when combined with the entropy
conservative flux $\vec{f_k}^\mathrm{num}$ presented in \cite{offner2018stability}. 
An entropy stable semidiscretization is thus obtained by adding a dissipative term $- \mat{Q} 
(\vec{u_k}_+ - \vec{u_k}_-)$ to the entropy conservative flux. 
Here we simply use a local Lax--Friedrichs type dissipation matrix
\begin{equation}
  \mat{Q} = \frac{{\alpha}}{2} \, \mat{\operatorname{I}}
  \quad \text{with} \quad
  {\alpha} = \max\set{ \abs{{\alpha}(-)}, \abs{{\alpha}(+)} }, 
\end{equation}
where $\abs{{\alpha}(\pm)}$ is the largest absolute value of all 
eigenvalues of the Jacobian $f'(u(\pm))$.

Even though the skew-symmetric formulation \eqref{eq:PC_semidiscsystem} combined with an 
appropriate numerical flux yields an entropy stable scheme, this test case demonstrates that 
additional regularization is still necessary to obtain reasonable numerical solutions. 
This was, for instance, stressed in \cite{ranocha2018stability}.
Here we demonstrate how $\ell^1$ regularization enhances the numerical solution for the nonlinear 
system \eqref{eq:system} with periodic boundary conditions and initial condition
\begin{align}
  u_0(x,0) & = 1 + 
  \left\{ 
  \begin{array}{lll}
    e \cdot \exp \left( -\frac{r^2}{r^2 - (x-0.5)^2} \right) 
      & , & \text{if } |x - 0.5| < r, \\
    0 & , & \text{if } |x - 0.5| \geq r, 
  \end{array} \right. \\
  u_1(x,0) & = 
  \left\{ 
  \begin{array}{lll}
    e \cdot \exp \left( -\frac{r^2}{r^2 - (x-0.5)^2} \right) 
      & , & \text{if } |x - 0.5| < r, \\
    0 & , & \text{if } |x - 0.5| \geq r, 
  \end{array} \right.
\end{align} 
where $r=0.5$.

For the more general case of systems of conservation laws, we propose a straightforward extension of 
our $\ell^1$ regularization technique. 
Specifically, the PA sensor \eqref{eq:sensor} is applied 
to every conserved variable $u_k$ separately and, once a discontinuity is detected, $\ell^1$ 
regularization \eqref{eq:op-problem2} is performed for the respective variable. 
For this test case, the ramp parameter in \eqref{eq:param-fun} has been chosen as $\kappa=0.9$.

\begin{figure}[!htb]
  \centering
  \begin{subfigure}[b]{0.32\textwidth}
    \includegraphics[width=\textwidth]{%
      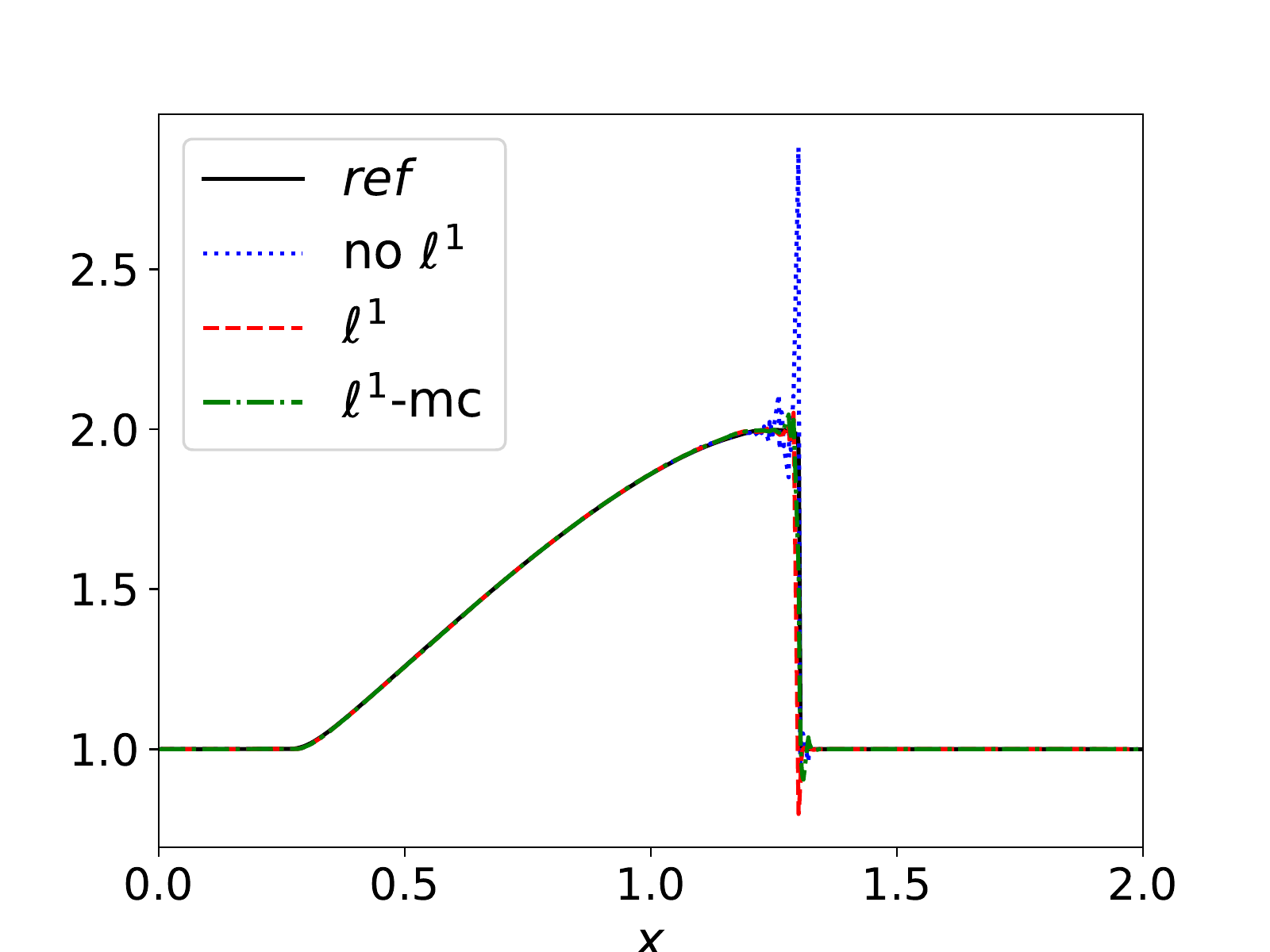}
    \caption{Coefficient $u_0$.}
    \label{fig:system_e-stable_u0} 
  \end{subfigure}%
  ~ 
  \begin{subfigure}[b]{0.32\textwidth}
    \includegraphics[width=\textwidth]{%
      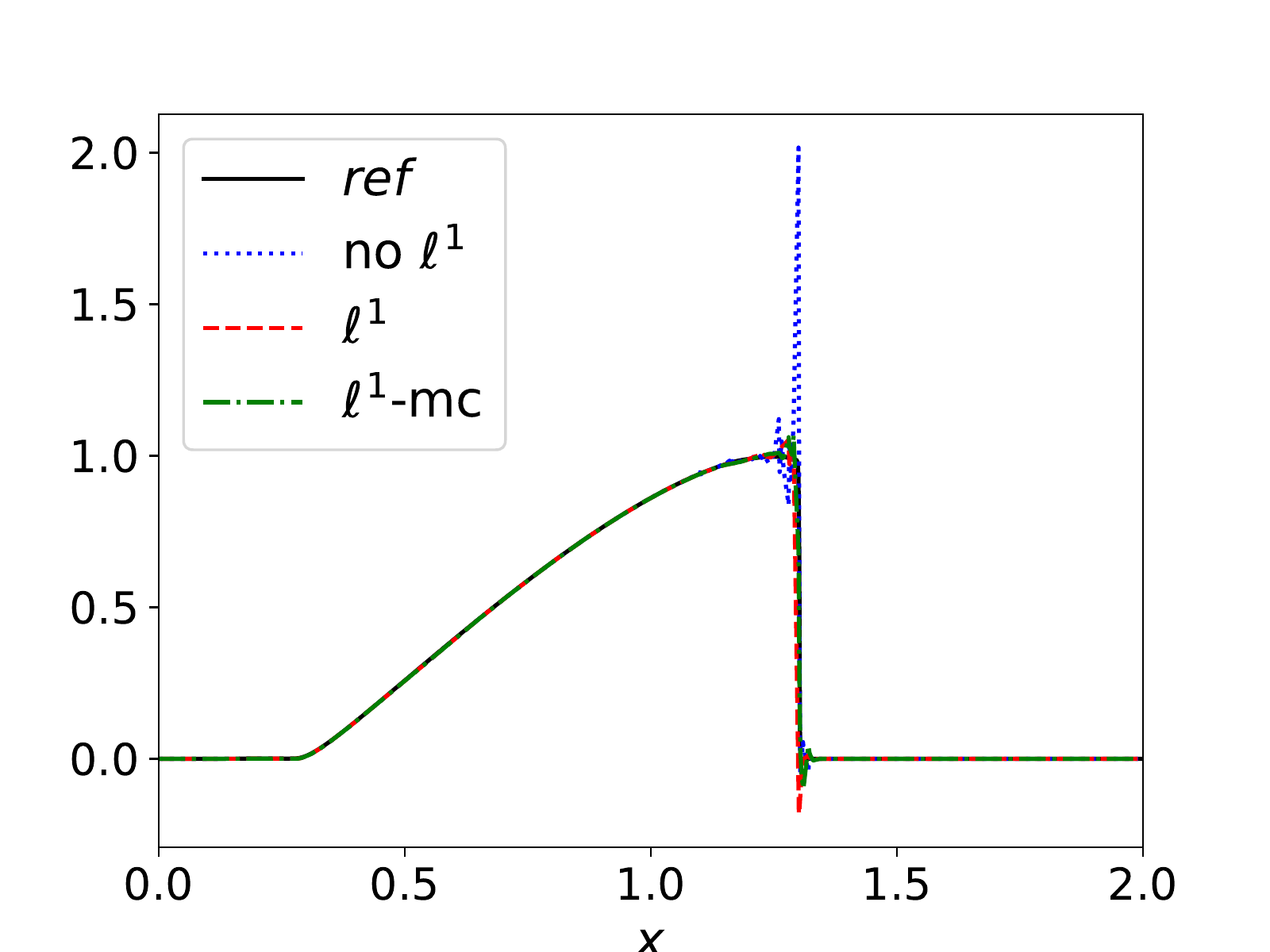}
    \caption{Coefficient $u_1$.}
    \label{fig:system_e-stable_u1} 
  \end{subfigure}%
  ~
  \begin{subfigure}[b]{0.32\textwidth}
    \includegraphics[width=\textwidth]{%
      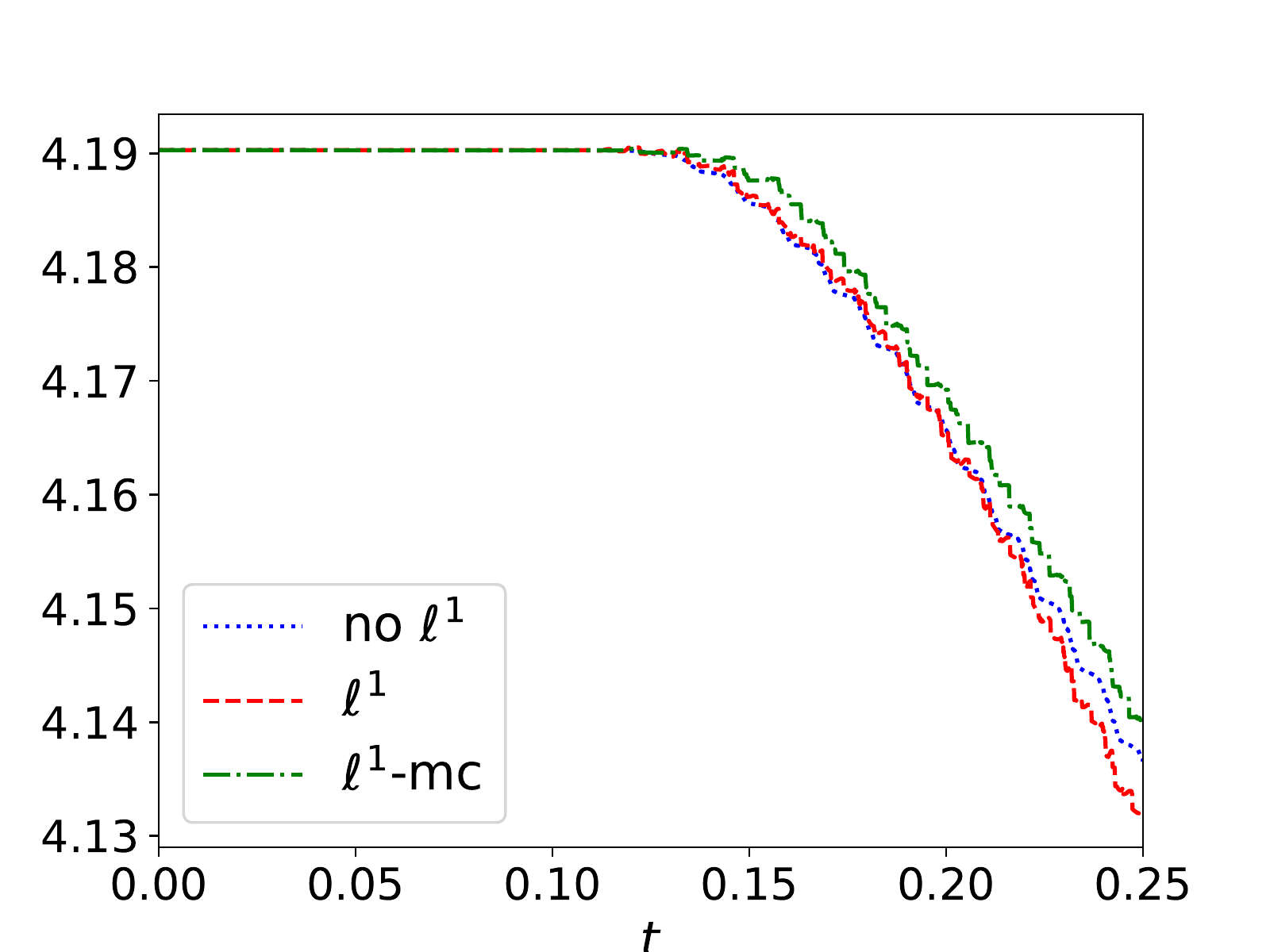}
    \caption{Energy over time.}
    \label{fig:system_e-stable_e} 
  \end{subfigure}%
  \\
  \begin{subfigure}[b]{0.32\textwidth}
    \includegraphics[width=\textwidth]{%
      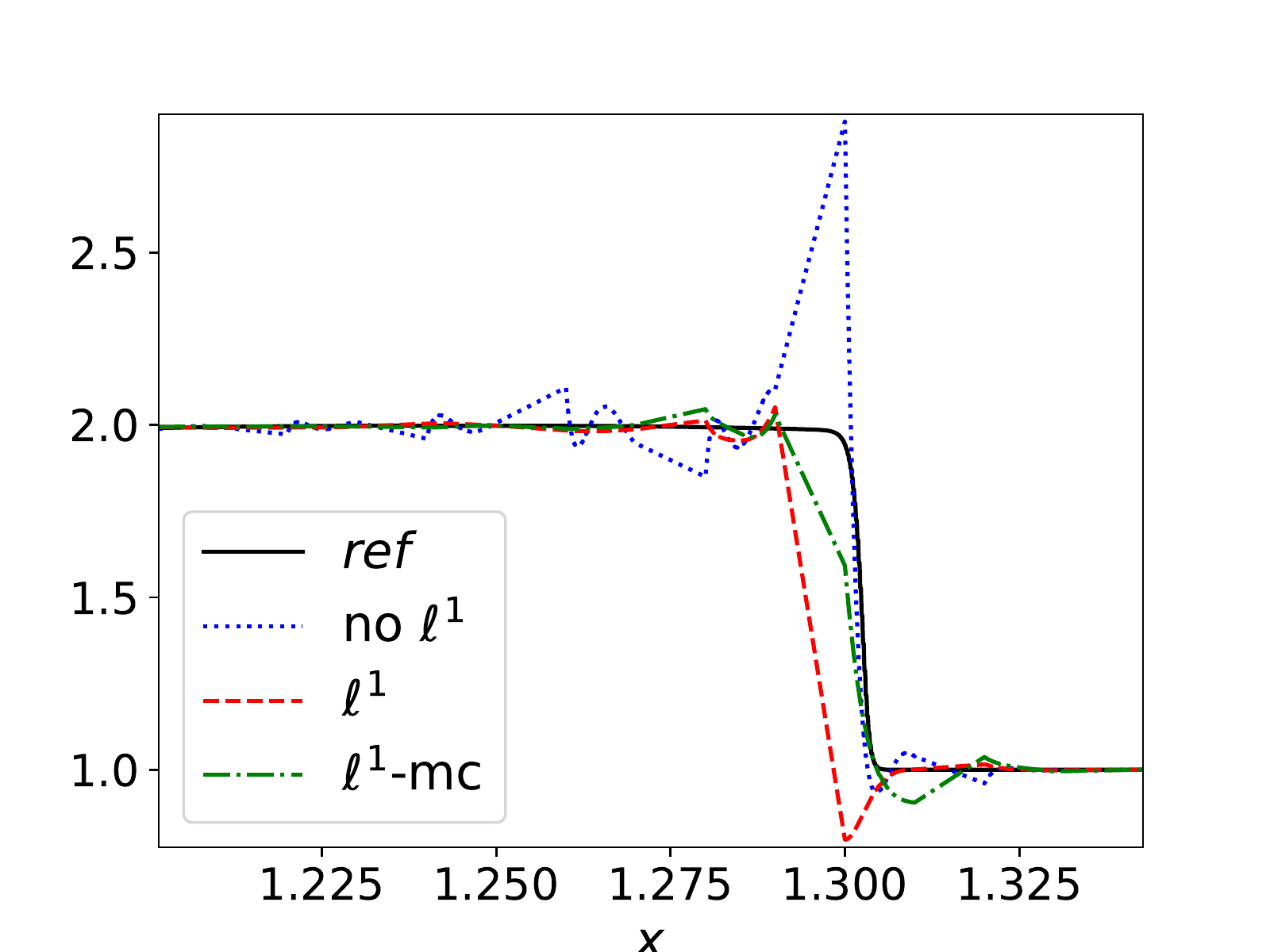}
    \caption{Coefficient $u_0$.}
    \label{fig:system_e-stable_u0_zoom} 
  \end{subfigure}%
  ~ 
  \begin{subfigure}[b]{0.32\textwidth}
    \includegraphics[width=\textwidth]{%
      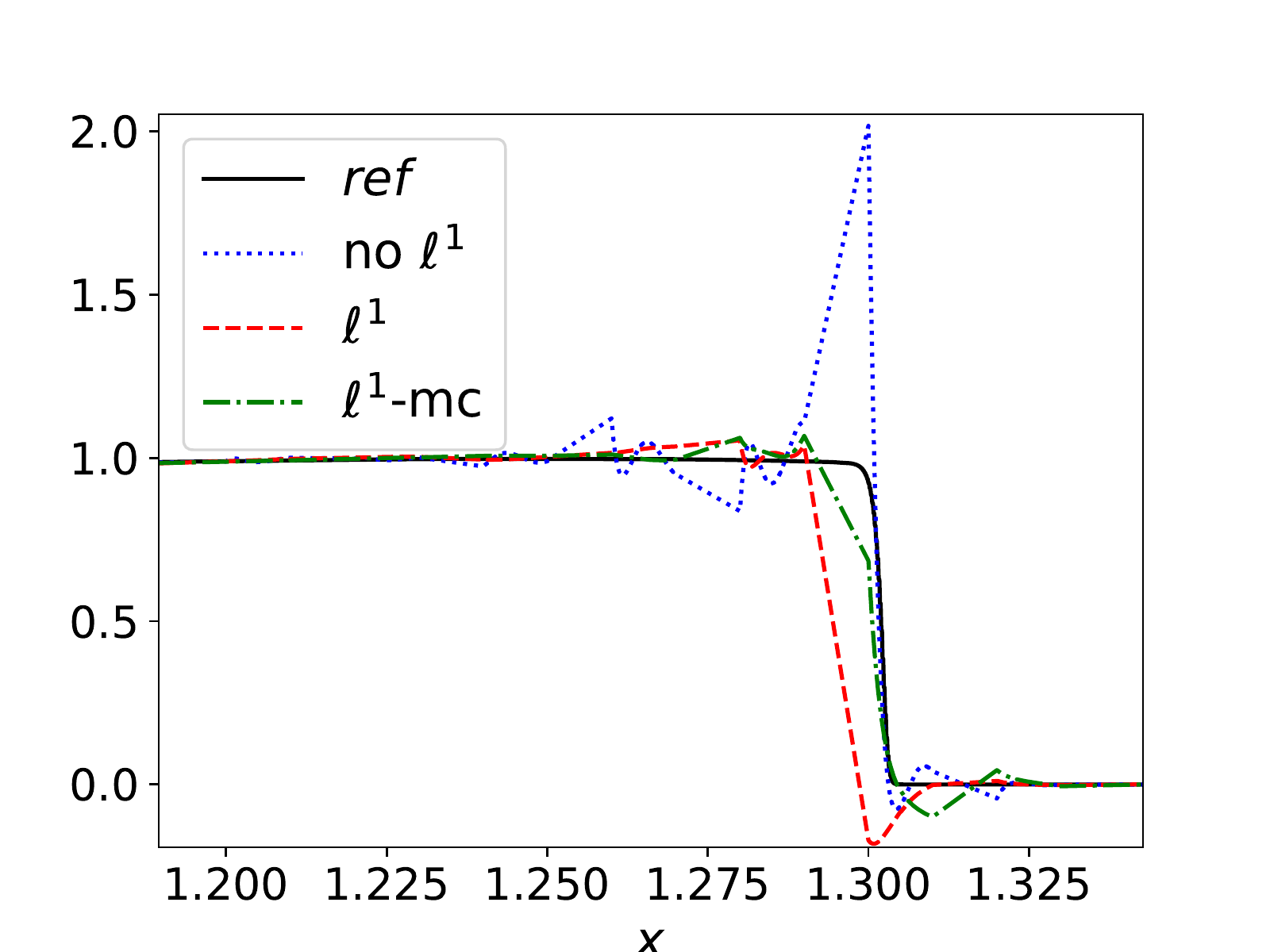}
    \caption{Coefficient $u_1$.}
    \label{fig:system_e-stable_u1_zoom} 
  \end{subfigure}%
  ~
  \begin{subfigure}[b]{0.32\textwidth}
    \includegraphics[width=\textwidth]{%
      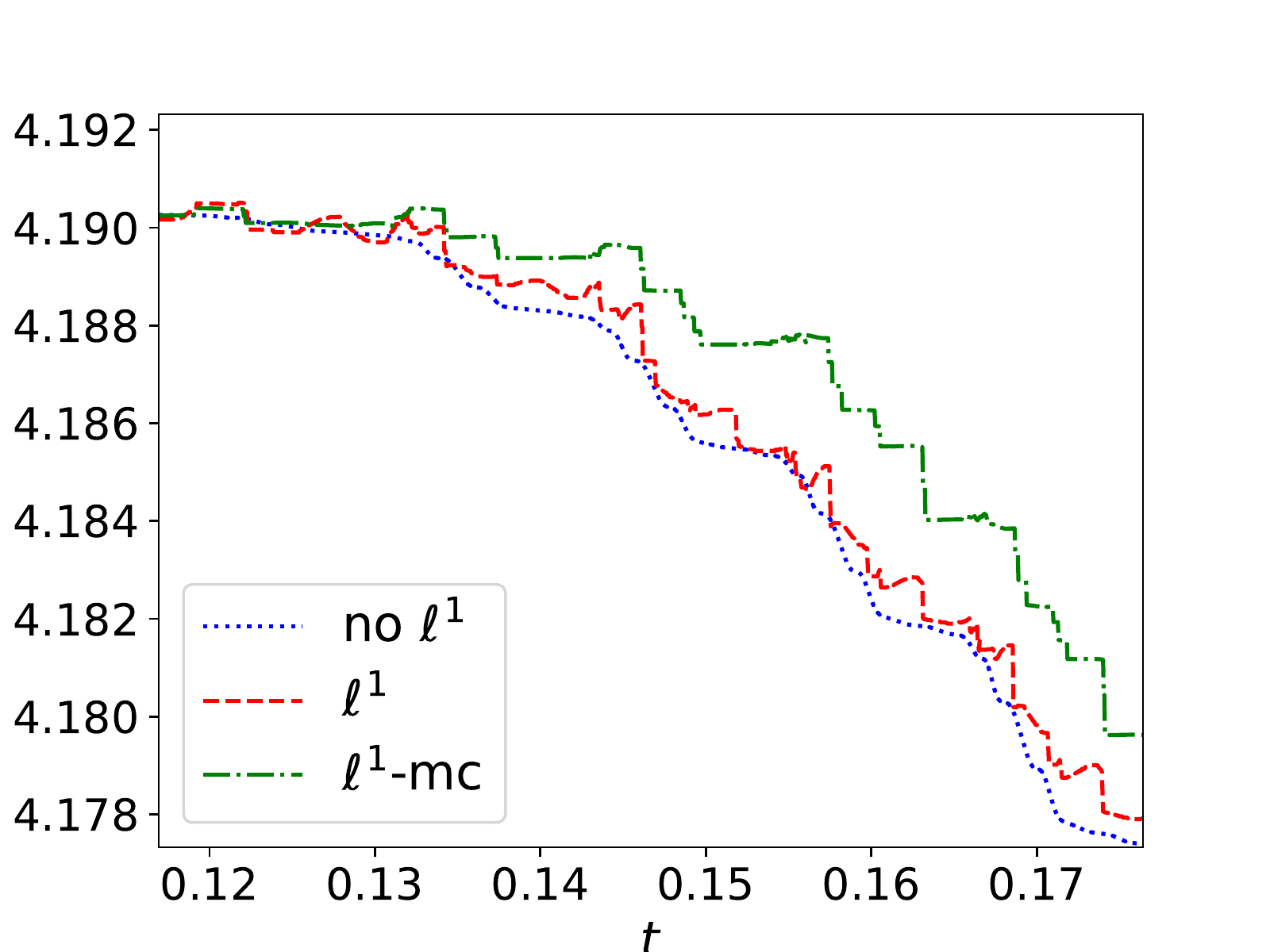}
    \caption{Energy over time.}
    \label{fig:system_e-stable_e_zoom} 
  \end{subfigure}%
  \caption{Components of the numerical solutions at $t=0.25$ by the DG 
method for $I=100$ elements and polynomial degree $p=6$. 
  Without $\ell^1$ regularization (straight blue line), with $\ell^1$ regularization (dashed red 
line), and with $\ell^1$ regularization and mass correction (dash-dotted green line).
  In all tests an entropy \emph{stable} numerical flux has been used.}
  \label{fig:system_e-stable}
\end{figure}

Figure \ref{fig:system_e-stable} illustrates the results of $\ell^1$ regularization (with 
and without mass correction) for the above described test case and for an entropy stable numerical 
flux. 
In all subsequent tests $I=100$ equidistant elements and a polynomial basis of degree $p=6$ have 
been used.
Further, for the $\ell^1$ regularization, the same parameters as before have been used, i.e.,  
$\lambda_{\text{max}}=4\cdot10^2$,  $K=400$, $\beta=20$, $\alpha = 0.0001$, and $tol = 0.001$. 

\begin{figure}[!htb]
  \centering
  \begin{subfigure}[b]{0.32\textwidth}
    \includegraphics[width=\textwidth]{%
      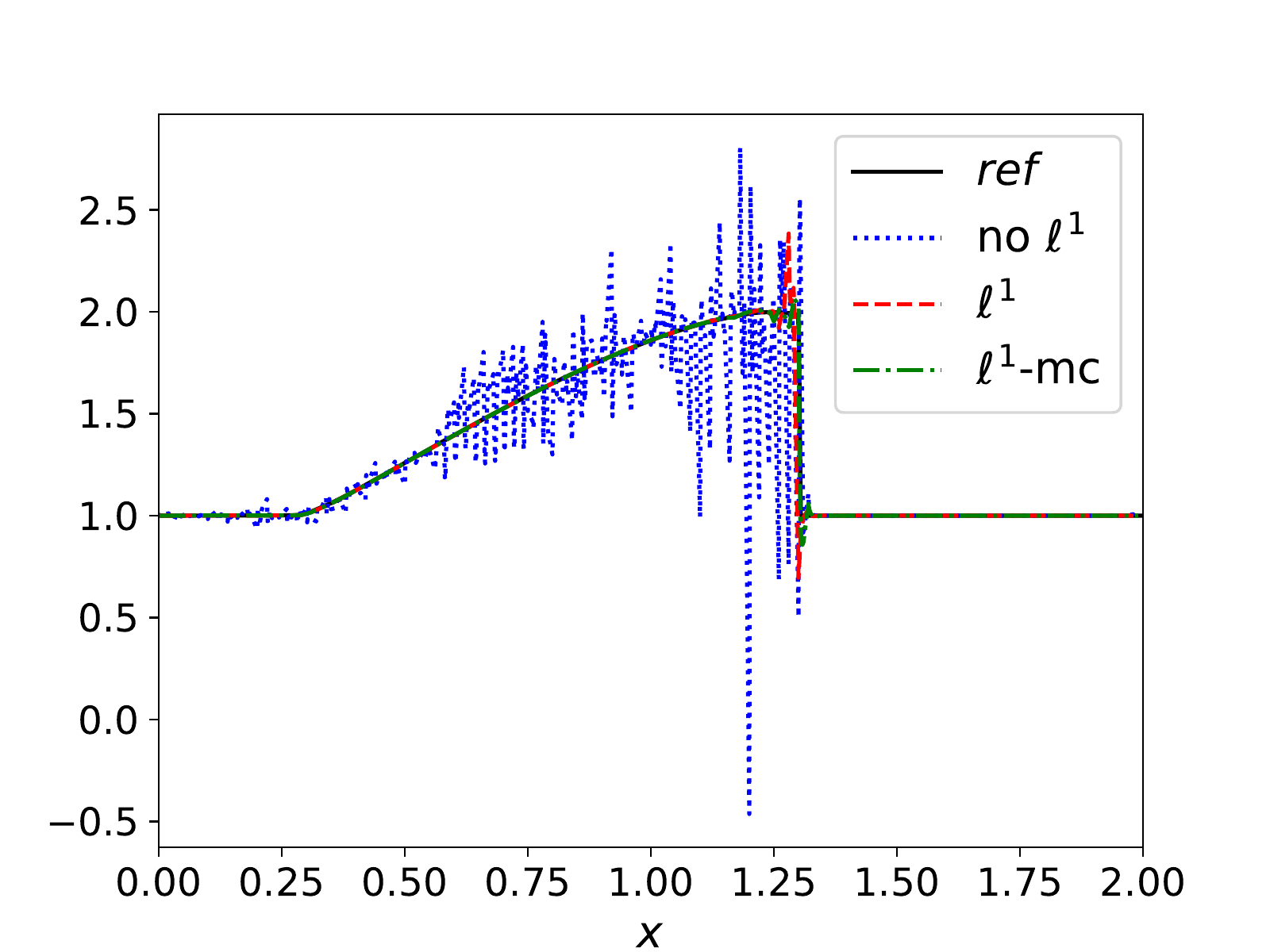}
    \caption{Coefficient $u_0$.}
    \label{fig:system_e-cons_u0} 
  \end{subfigure}%
  ~ 
  \begin{subfigure}[b]{0.32\textwidth}
    \includegraphics[width=\textwidth]{%
      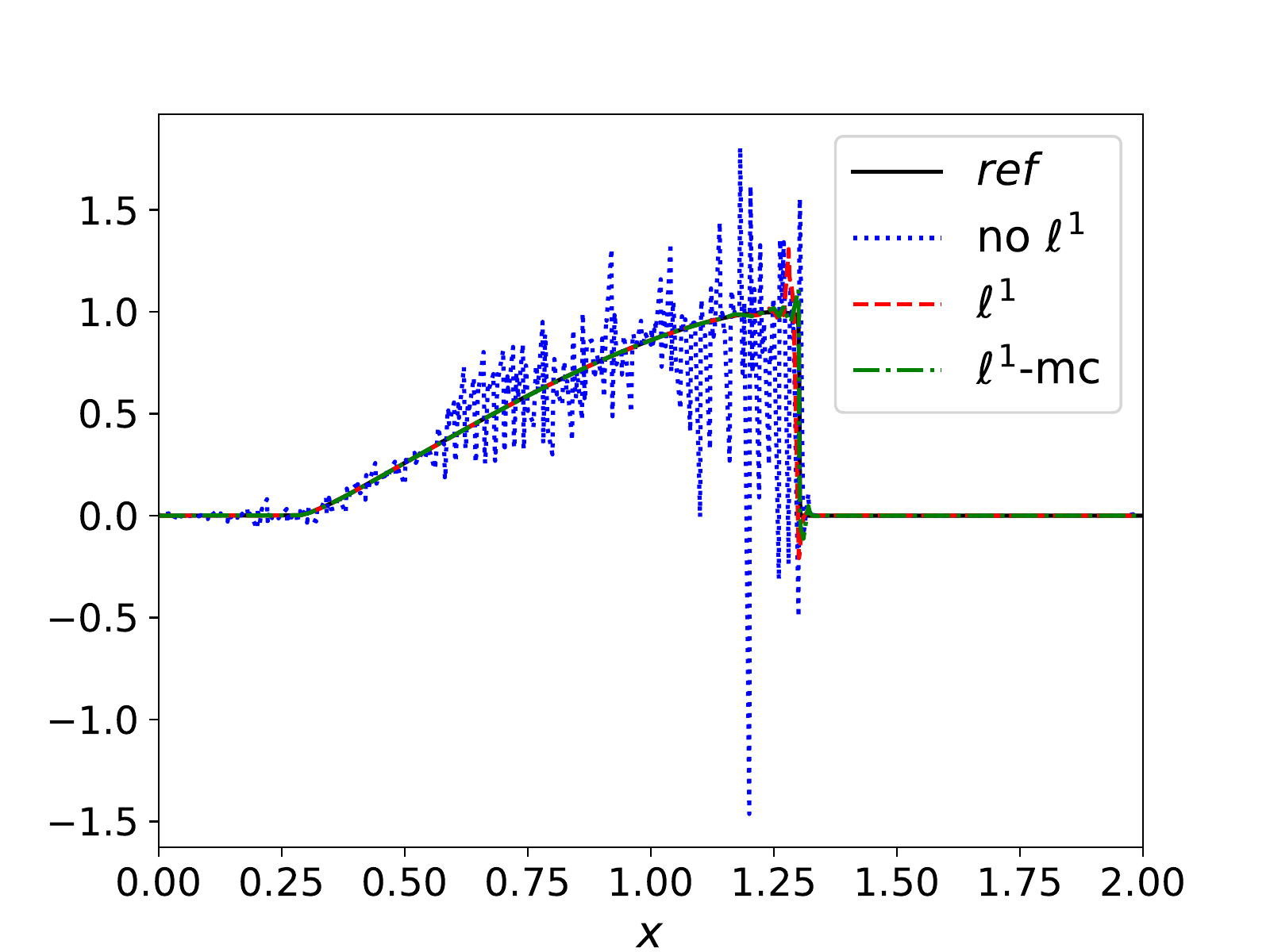}
    \caption{Coefficient $u_1$.}
    \label{fig:system_e-cons_u1} 
  \end{subfigure}%
  ~
  \begin{subfigure}[b]{0.32\textwidth}
    \includegraphics[width=\textwidth]{%
      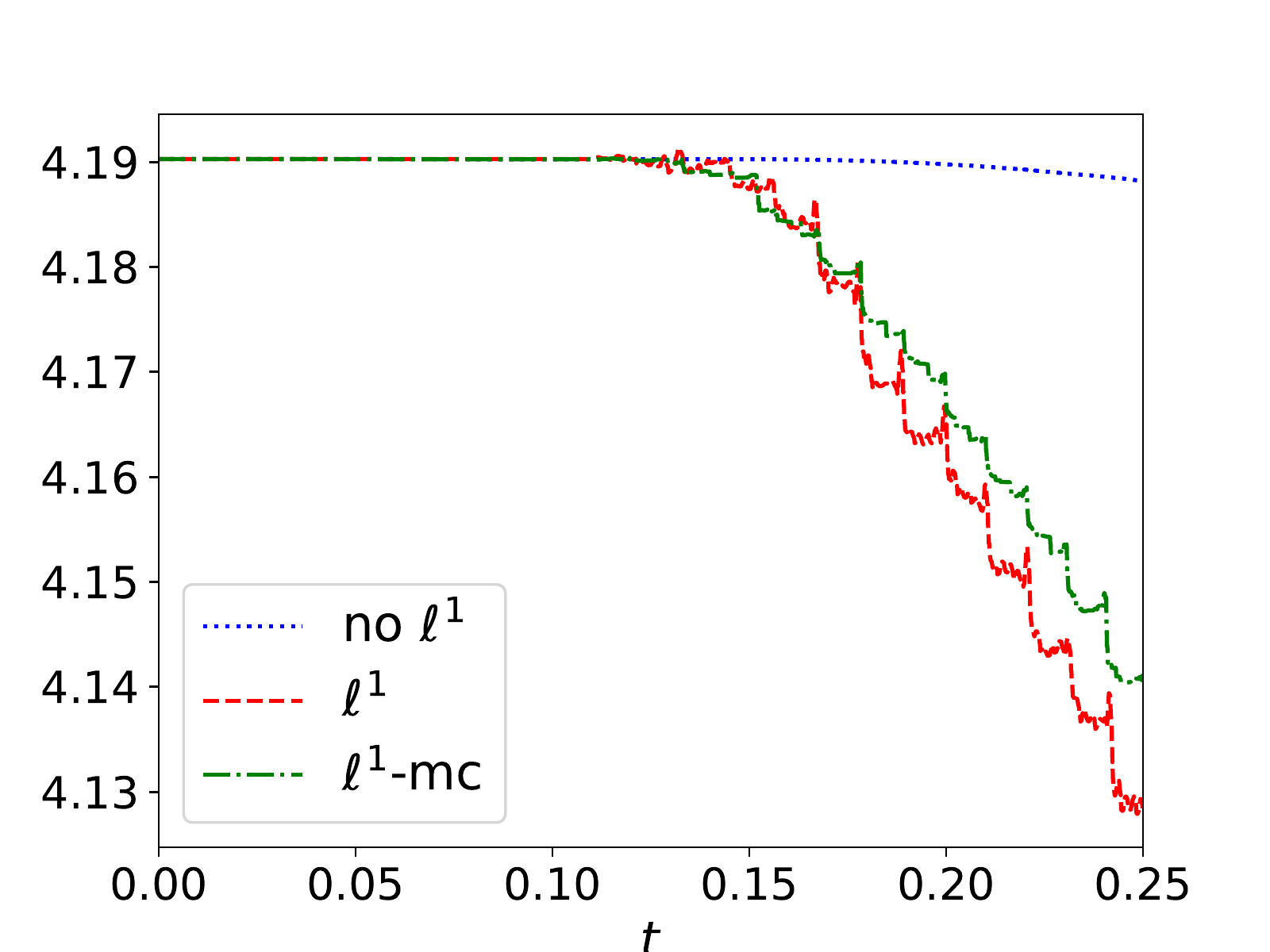}
    \caption{Energy over time.}
    \label{fig:system_e-cons_e} 
  \end{subfigure}%
  \\
  \begin{subfigure}[b]{0.32\textwidth}
    \includegraphics[width=\textwidth]{%
      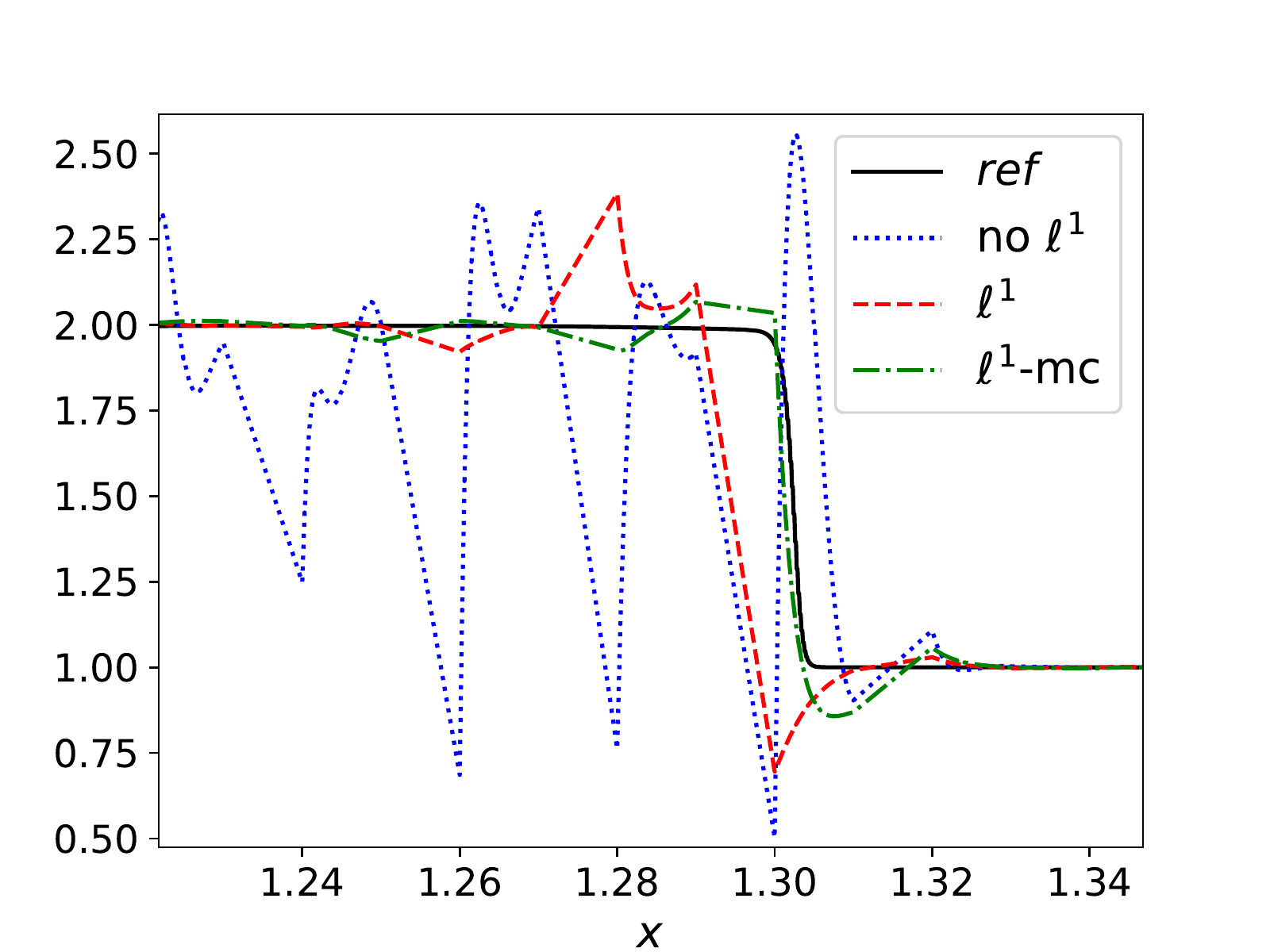}
    \caption{Coefficient $u_0$.}
    \label{fig:system_e-cons_u0_zoom} 
  \end{subfigure}%
  ~ 
  \begin{subfigure}[b]{0.32\textwidth}
    \includegraphics[width=\textwidth]{%
      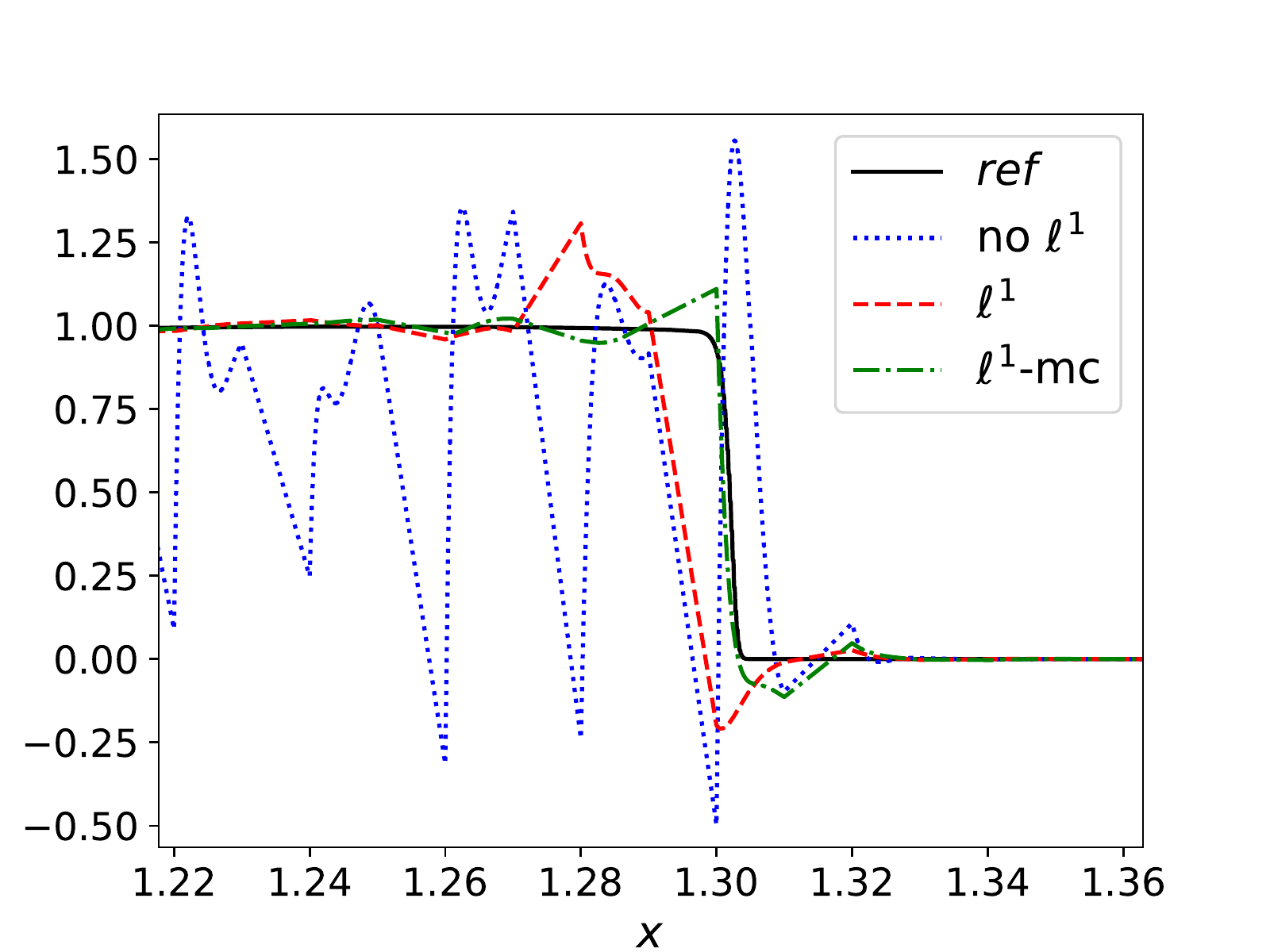}
    \caption{Coefficient $u_1$.}
    \label{fig:system_e-cons_u1_zoom} 
  \end{subfigure}%
  ~
  \begin{subfigure}[b]{0.32\textwidth}
    \includegraphics[width=\textwidth]{%
      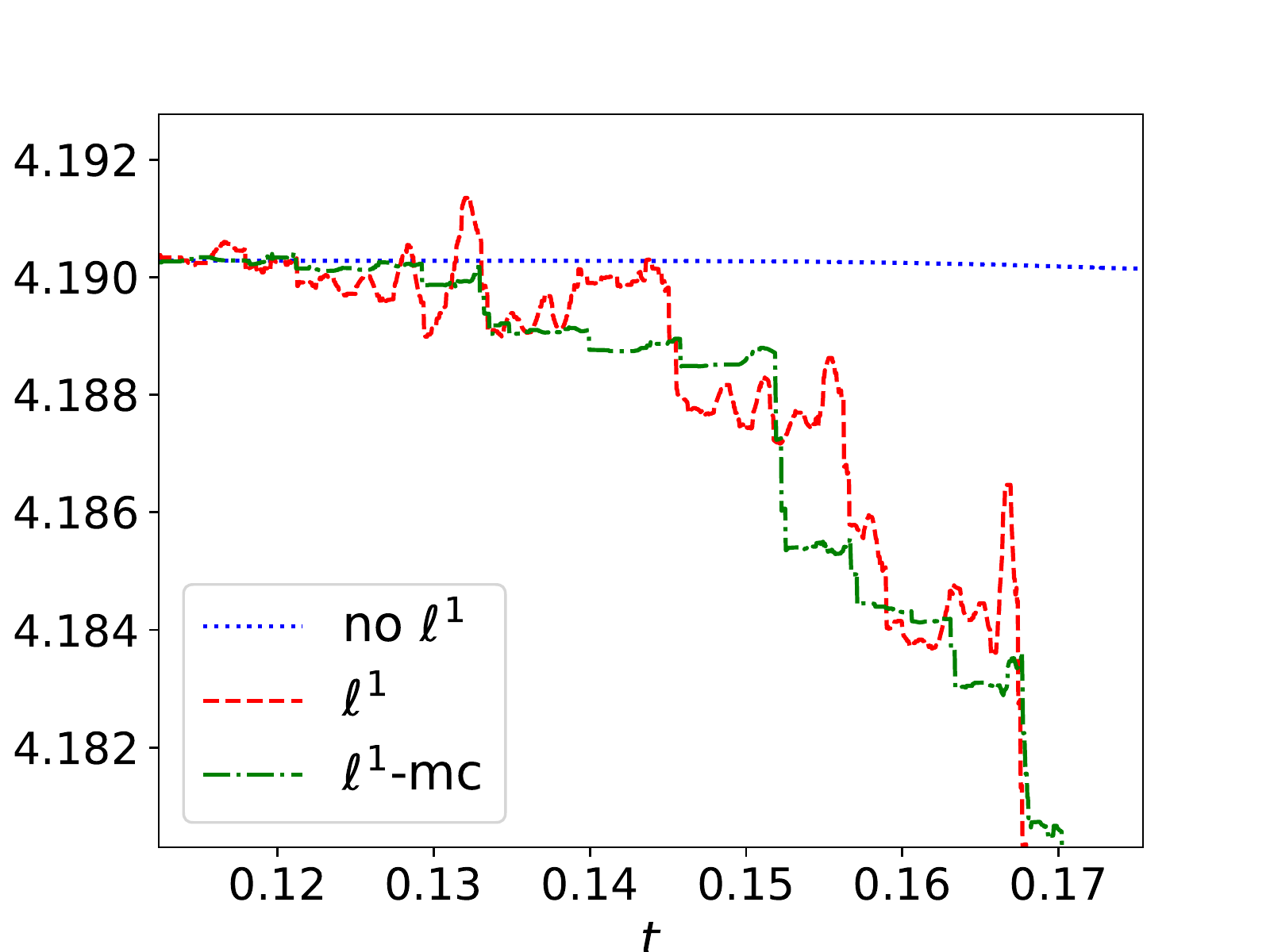}
    \caption{Energy over time.}
    \label{fig:system_e-cons_e_zoom} 
  \end{subfigure}%
  \caption{Components of the numerical solutions at $t=0.25$ by the DG 
method for $I=100$ elements and polynomial degree $p=6$. 
  Without $\ell^1$ regularization (straight blue line), with $\ell^1$ regularization (dashed red 
line), and with $\ell^1$ regularization and mass correction (dash-dotted green line).
  In all tests an entropy \emph{conservative} numerical flux has been used.}
  \label{fig:system_e-cons}
\end{figure}

Note that while the numerical solution without $\ell^1$ regularization shows heavy oscillations in 
both components, the numerical solution with $\ell^1$ regularization provides a significantly 
sharper profile.
Further, by consulting Figures \ref{fig:system_e-stable_u0_zoom} and  
\ref{fig:system_e-stable_u1_zoom}, it should be stressed that only $\ell^1$ regularization with 
additional mass correction is able to capture the exact shock location. 
Due to missing conservation, $\ell^1$ regularization without mass correction results in a slightly 
wrong location for the shock.
Finally, Figures \ref{fig:system_e-stable_e} and \ref{fig:system_e-stable_e_zoom} illustrate 
the energy of the different methods over time. 
We note from these figures that $\ell^1$ regularization (with and without mass correction) 
slightly increase the energy in this test case. 

In order to further emphasize the effect of $\ell^1$ regularization (with and without mass 
correction), similar results using an
entropy conservative numerical flux at the interfaces between elements are shown
in Figure \ref{fig:system_e-cons}. 
Once again, $\ell^1$ regularization is demonstrated to improve the numerical 
solution. 
The best results are obtained when $\ell^1$ regularization is combined with mass correction. 
In particular, only $\ell^1$ regularization with additional mass correction is able to accurately 
capture the shock at the right location.
Finally, consulting Figures \ref{fig:system_e-cons_e} and \ref{fig:system_e-cons_e_zoom}, 
$\ell^1$ regularization decreases the energy overall. 
Yet, it is demonstrated once more that an entropy inequality is not satisfied by $\ell^1$ 
regularization, i.e., the energy might increase as well as decrease by utilizing $\ell^1$ 
regularization. 
Thus, future work will focus on incorporating energy stability (as well as other properties like 
TVD or positivity) by additional constraints in the minimization problem \eqref{eq:op-problem}.
\section{Concluding remarks}
\label{sec:conclusion}

We have presented a novel approach to shock capturing by $\ell^1$ regularization using SE approximations. 
Our work not only is distinguished from previous studies 
\cite{schaeffer2013sparse,hou2015sparse+,lavery1989solution,lavery1991solution,guermond2008fast} 
by focusing on discontinuous solutions but further by promoting sparsity of the jump function instead 
of the numerical solution itself.   
By approximating the jump function with the (high order) PA operator,
we help to eliminate the staircase effect that arises for classical TV operators. 
Our results demonstrate that it is possible to efficiently implement a method that yields increased accuracy 
and better resolves (shock) discontinuities. 
In particular, no additional time step restrictions are introduced, 
in contrast to artificial viscosity methods when no care is taken in their construction. 
This approach for solving numerical conservation laws was first used in \cite{scarnati2017using}, 
where the Lax--Wendroff scheme and Chebyshev and Fourier spectral methods were used as the numerical 
PDE solver.
Our method improves upon the approach in \cite{scarnati2017using} in two ways. 
First, we employ the SE approximation for solving the conservation law, which allows 
element-to-element variations in the optimization problem.
In particular, $\ell^1$ regularization is only activated in troubled elements, which 
enhances accuracy and efficiency of the method. 
Second, in the process we proposed a novel discontinuity sensor based on PA operators of increasing orders, 
which is able to flag troubled elements as well as to steer the amount of regularization introduced by the sparse 
reconstruction. 

Numerical tests demonstrate the method using a nodal collocation-type discontinuous Galerkin 
method for the inviscid Burgers' equation, the linear advection equation, 
and a nonlinear system of conservation laws. 
Our results show that the method yields improved accuracy and robustness. 

No effort was made in our study to optimize any of the parameters involved in solving 
the optimization problem. 
This will be addressed in future work, along with the possibility to include additional constraints 
(e.g., for entropy, TVD, and positivity constraints), 
since preliminary results presented here are encouraging. 
The generalization of the approach itself to higher dimensions is straightforward and has already been demonstrated 
in \cite{scarnati2017using}. 
Of interest, however, would be the extension of the proposed approach to other classes of methods, such as finite 
volume methods. 
We believe $\ell^1$ regularization might be an important ingredient to make high order methods viable in several 
research applications.

\section*{Acknowledgements}
The authors would like to thank Chi-Wang Shu (Brown University) for helpful advice. 
Further, the authors would like to thank the anonymous referees for many helpful 
suggestions, resulting in an improved presentation of this work. 
Jan Glaubitz' work was supported by the German Research Foundation (DFG, Deutsche 
Forschungsgemeinschaft) under grant SO 363/15-1.
Anne Gelb's work was partially supported by AFOSR9550-18-1-0316 and NSF-DMS 1502640.

\bibliographystyle{elsarticle-num-names}
\bibliography{literature}

\end{document}